\documentclass[12pt,letterpaper]{article}

\usepackage{amsmath, amssymb}
\usepackage{lscape}
\usepackage{epsfig}
\usepackage{cite}

\usepackage{makeidx}         
\usepackage{subfig}
\usepackage{graphicx}        

\newcommand{\vect} \mathbf

\def\R{{\mathbb{R}}}

\begin{document}

\title{Robust artificial neural networks and outlier detection. \\  Technical report}
\author{Gleb Beliakov$^1$, Andrei Kelarev$^1$ and John Yearwood$^2$
\\[.2in]
$^1$School of Information Technology\\
Deakin University, 221 Burwood Hwy, Burwood 3125, Australia \\
\texttt{\{gleb,andrei\}@deakin.edu.au}\\[.2in]
$^2$School of Science, Information Technology and Engineering\\
University of Ballarat, P.O. Box 663, Ballarat, Victoria 3353, Australia\\
\texttt{j.yearwood@ballarat.edu.au}
}

\date{}
\maketitle

\begin{abstract}
Large outliers break down linear
and nonlinear regression models.
Robust regression methods allow one to
filter out the outliers when building a
model. By replacing the traditional
least squares criterion with the least
trimmed squares criterion, in which half
of data is treated as potential
outliers, one can fit accurate
regression models to strongly
contaminated data.
High-breakdown methods have become
very well established in linear
regression, but have started being
applied for non-linear regression
only recently.
In this work, we examine the problem
of fitting artificial neural networks
to contaminated data using least trimmed
squares criterion. We introduce a penalized least trimmed squares criterion which prevents unnecessary removal of valid data.
Training of ANNs leads
to a challenging non-smooth global
optimization problem. We compare the
efficiency of several derivative-free
optimization methods in solving it,
and show that our approach identifies
the outliers correctly when ANNs
are used for nonlinear regression.
\end{abstract}


\newpage

\baselineskip=\normalbaselineskip


\section{Introduction}\label{s:introduction}

We consider the problem of function
approximation with artificial neural
networks (ANNs). This generic task has
many applications in science and
engineering, such as signal processing,
pattern recognition, and control. The
goal is to model an unknown nonlinear
function based on observed input-output
pairs. ANNs are universal function
approximators \cite{Hornik1989}, and
usually they deliver good performance in
applications.

Error-free data are rarely provided in
applications. First, the data are
usually contaminated by noise, which
reflects inaccuracies in observations
and stochastic nature of the underlying
process. ANNs and other function
approximators deal with such noise quite
efficiently, by minimising the sum of
squared differences between the observed
and predicted values, or a more
sophisticated fitting criterion, such as
Huber-type functions, which are the
basis of M-estimators in
statistics~\cite{Huber1991_book}.

The second type of data contamination
has to do with either gross observation
errors (e.g., equipment malfunction, or
notorious replacing missing values with
zeroes, wrong decimal points and other
blunders), or the data reflecting a
mixture of different phenomena. These
data, which usually take aberrant
values, are called \emph{outliers}.
It has been noted  that typically
the occurrence of outliers in routine data
ranges from 1\% to 10\%. When fitting
a model to the data, outliers need to be
identified and eliminated, or,
alternatively, examined closely, as
they may be of the main interest
themselves. Notable examples are intrusion
and cyberattack detection, detection of
harmful chemicals and cancerous cells.

The methods of function approximation
based on the least squares (or, more
generally, maximum likelihood principle)
are not robust against outliers. In fact
just one aberrant value can make the
model's bias infinite (it is said that
the method breaks down). This phenomenon
is well known in linear regression
\cite{Rousseeuw1987_book,Maronna_book},
where a number of robust high-breakdown
methods have been developed. The popular
methods of least median of squares (LMS)
and least trimmed squares (LTS)
\cite{Rousseeuw1984_JASA} discard half
of the data as potential outliers and
fit a model to the remaining half. These
methods determine numerically which half
of the data should be discarded in order
to obtain the smallest value of the
respective objectives. That way, up to
half of the data can be outliers but
they do not break down the method. It is
said that their (asymptotic) breakdown
point is $\frac{1}{2}$. The outliers
themselves can be identified by their
large residuals, something that cannot
be achieved when using the least squares
estimators, or maximum likelihood
estimators (called M-estimators),
because of the masking effect (i.e., the
outliers affect the fitted model so much
that their residuals do not stand out).

Much less work has been devoted to
non-linear high-breakdown regression.
There are very few papers dealing with
the LTS method applied to ANNs, see
for example, \cite{Stromberg1992}. More
recently, Liano~\cite{Liano1996} used
M-estimators to study the mechanism by
which outliers affect the resulting ANNs.
Chen and et.
al~\cite{Chen1994}
also used M-estimators as a robust
estimator in the presence of outliers.
The Least Trimmed Squares estimator was
discussed
in~\cite{Jeng2011,Rusiecki2007,Rusiecki2008}.
A robust LTS backpropagation algorithm
based on simulated annealing was considered
in~\cite{Rusiecki2008}.

It is known that fitting the LTS or LMS
criterion is an NP-hard problem even in
linear regression. The objective
has a large number of local minima and
is non-smooth. The problem becomes even
more complicated for ANNs, because
\begin{itemize}
\item[(a)]
the dimensionality of the problem increases
even further with the number of hidden
neurons used, and
\item[(b)]
the ANN training is an NP-hard problem
even when using the traditional least
squares criterion.
\end{itemize}
Training ANNs with LTS or LMS
criterion is very challenging because of
a much higher number of local minima as
well as non-applicability of the
traditional fast backpropagation
algorithm because of non-smoothness of
the objective.

In this article we advance the methods
for robust fitting of ANNs using LTS and
related fitting criteria. Our first
contribution is to design a hybrid algorithm, which combines a derivative free optimisation method for initial training of ANN, removal of the detected outliers, and then fine tuning of ANN weights using clean data and  backpropagation. The second contribution is the design of an improved fitting criterion, called Penalised CLTS (PCLTS), which prevents unnecessary removal of valid data.  The LTS and LMS criteria have this undesirable effect, illustrated in our experiments. The PCLTS criterion prevents unnecessary removals by imposing a penalty on removal of every datum.

This article is structured as follows.
In Section~\ref{sec1}, we introduce the
problem of robust regression, recall
the definitions and the main features
of several existing high-breakdown
estimators, and discuss the associated
optimization problem.
In Section \ref{sec_criterion} we introduce the PCLTS criterion for ANN fitting.
In Section~\ref{sec_methods}, we outline
the existing approaches to solution of
the related optimization problem and present
three new methods we use in this study.
Section~\ref{sec_expeer} is devoted to
a comparative numerical study of the
optimization methods using several data sets.
Section~\ref{sec_concl} concludes the
article.

\section{High-breakdown robust estimators} \label{sec1}

In this section, we briefly introduce a
statistical problem which is in the
origin of the non-smooth global
optimization problem treated in this
paper. We start the discussion with
high-breakdown linear regression,
which will be followed by nonlinear regression.

\subsection{Robust linear regression}

Consider the standard multiple linear
regression model
\begin{equation} \label{model-linear}
y_i = x_{i1}\theta_1 + \ldots + x_{ip}\theta_p +\varepsilon_i ,
\quad\mbox{ for } i=1,\ldots,n,
\end{equation}
where $x_{ip}=1$ for regression with an
intercept term. Here $\{x_{ij}\}=X\in
\R^{n \times p}$ is the matrix of
explanatory variables of full column
rank, and $\varepsilon$ is a $n$-vector
of independent and identically
distributed  random errors with zero
mean and (unknown) variance $\sigma^2$.
The goal is to determine the vector of
unknown parameters $\theta \in \R^p$.
The goodness of fit is  expressed in
terms of the residuals $r_i(\vect{\theta})
= \vect{\theta}^t \vect{x}_i-\vect{y}_i$,
namely the  sum
$\sum_{i=1}^n r_i^2$ for the ordinary
least squares (OLS) and $\sum_{i=1}^n
|r_i|$  for the least absolute
deviations (LAD) methods.

As we mentioned in the introduction, the
OLS and LAD are very sensitive to large
outliers in the data. Just one grossly
atypical datum can affect the model. A
\emph{breakdown point} of a regression
estimator  is the smallest proportion
of contaminated data that can make the
estimator's bias arbitrarily large (see
\cite{Hampel1971}, \cite{Rousseeuw1987_book}, p.10).
The breakdown point of the  OLS and LAD
methods tends to zero as $\frac{1}{n}$ with the increasing
sample size $n$, and is said to be~0\%.

To overcome the lack of robustness of
OLS and LAD estimators, Rousseeuw
\cite{Rousseeuw1984_JASA} introduced the
least median of squares  and least
trimmed squares  estimators. The methods
of the least trimmed absolute deviations
(LTA) and the maximum trimmed likelihood
(MTL) were advocated in
\cite{Stromberg2000_JASA,Hawkins1999_CSDA,Hadi1997_CSDA}.
These methods are robust to leverage
points, and allow up to a half of the
data to be contaminated without
affecting the regression model. Atypical
data are then detected by their large
residuals. For recent accounts of the
state-of-the-art in high-breakdown
robust regression, see
\cite{Rousseeuw1987_book,Maronna_book}.

Essentially, the LTS, LTA and MTL
methods work in the following way. Half
of the sample is discarded, and a
regression model is built using the
other half. The sum of the residuals is
then evaluated. The objective is to find
the optimal partition of the data into
two halves, so that the sum of the
(squared or absolute) residuals is the
smallest. This is a combinatorial
formulation of the problem. Evidently
the solution is feasible only for small
data sets.

The same problem can be formulated as a
continuous non-smooth optimization
problem:
\begin{equation} \label{LTSOBJ}
    \mbox{Minimize }F(\theta)=\sum_{i=1}^h  (r_{(i)}(\theta))^2,
\end{equation}
where the residuals are ordered in the
increasing order
$|r_{(1)}|\leq |r_{(2)}|\leq \ldots \leq |r_{(n)}|$,
and $h=[(n+p+1)/2]$, where $[\cdot]$
is the floor function.
The variables in this model are
the regression coefficients $\theta$. For
small to moderate dimension but large
data sets this model offers significant
numerical advantages. It is the basis of
fast heuristic algorithms in
high-breakdown regression
\cite{Rousseeuw2006_DMKD}. Several
recent methods based on  this
formulation, including evolutionary and
semidefinite programming, were presented
in
\cite{Nunkesser2010_CSDA,Nguyen2010_CSDA,Cerioli2010_JCGS,Schyns2010_CSDA,Beliakov2012_OMS}.

The four methods mentioned above
achieve the highest attainable asymptotic breakdown
point of $1/2$, which means that up to
a half of the data can be contaminated
without affecting the estimator.
Consequently, the outliers can be
easily detected by their
large residuals, and either eliminated,
or alternatively, examined more closely
in the cases where the outliers
themselves are of the main interest.

It is shown in \cite{Bernholt2005_TR}
that computation of the high-breakdown
estimators is an NP-hard problem.
Indeed, the objective $F$ in the
methods mentioned above  has multiple
local minima. Consider, for instance,
the LTS estimator.
The problem~(\ref{LTSOBJ}) can be written as
 \begin{equation*}\label{prob_regressionOWAinc}
    \min_\theta \left( \min_{\pi} \sum_{i=1}^h r_{\pi(i)}^2(\theta)\right),
\end{equation*}
where $\pi$ is a permutation of the vector $(1,2,\ldots,n)$.
Subsequently we write it as
 \begin{equation*}\label{prob_regressionOWAinc1}
    \min_{\pi} \left(\min_{\theta } \sum_{i=1}^h  r_{\pi(i)}^2(\theta)\right).
\end{equation*}
The inner optimization problem is
convex, and has a unique minimum
(potentially, multiple minimizers). Then
problem (\ref{LTSOBJ}) will potentially
have as many as $ n \choose h$ local
minima (the number of permutations $\pi$
which result in distinct sums  for any
fixed~$\theta$).

In addition to determining the breakdown
point of an estimator, it is also essential
to ivestigate its efficiency. It is
customary to evaluate the efficiency by
comparing it to that of the OLS estimator.
A fully efficient estimator should deliver
the same accuracy as the maximum likelihood
based estimator (which is OLS when the
noise is Gaussian) when the data set
contains no outliers.
The relative efficiency of the LMS, LTS
and LTA methods for normally distributed
data is low
\cite{Maronna_book,Rousseeuw1987_book}.
However, fully efficient high-breakdown
estimators exist. The reweighted least
squares estimator (REWLSE) is one such estimator
presented in~\cite{Gervini2002_AS}.

To improve the efficiency, while preserving
the breakdown point, Gervini and
Yohai~\cite{Gervini2002_AS} use a
two-step process: an initial high-breakdown
estimator (like LMS, or LTS) provides a
robust estimate of scale $s$ used to
re-weigh the data.
The weights are given by the formula
\begin{equation}\label{weights}
w_i=\left\{
\begin{array}{cc}
                   1 & \mbox{ if } |r_i| < t_n, \\
                   0 & \mbox{ otherwise,}
                 \end{array}
\right.
\end{equation}
where $t_n$ is the adaptive cutoff value
beyond which the sample proportion of
absolute residuals exceeds the
theoretical proportion.
The weights $w_i$ given by
(\ref{weights}) effectively remove all
outliers. The adaptive estimate REWLSE
is then computed as a weighted OLS:
$\theta=(X^tWX)^{-1}X^tWY$ with
$W=diag(w_1,\ldots,w_n)$.

The two-step process in REWLSE computes
fully efficient estimators when the data
are normally distributed, as no data are
unnecessarily removed. REWLSE  inherits
the breakdown point of the initial
estimator and combines it with full
efficiency of the final LS estimator.

\subsection{Nonlinear robust regression}

We now consider a nonlinear regression model
\begin{equation} \label{model}
y_i = f(\vect x_i; \theta) +\varepsilon_i , \:\mbox{ for } i=1,\ldots,n,
\end{equation}
where $f$ is an arbitrary (nonlinear) function,
$\theta$ is a set of parameters
(it may vary depending on the particular
specification of a class of regression models),
$\vect x_i \in \R^{p-1}$ are the fixed data
points or inputs and $y_i$ are the outputs.
Regression neural networks give us examples
of such functions. For instance, each neural
network with one hidden layer and one output
defines a function of the form
$$
f(\vect x,\theta) =
\sum_{j=1}^{m_h} \theta^h_j \cdot
g \left( \sum_{k=1}^{p} \theta^i_{jk} x_k \right),
$$
where $\theta^h \in \R^{m_h}$
(and $\theta^i \in \R^{m_h \times p}$)
are the hidden (respectively, input)
layer weights, $m_h$ is the size of the
hidden layer, $p$ is the number of
inputs plus one (the bias term), and $g$
is a transfer function. Altogether, $f$
has $(p+1) \times m_h$ parameters
represented by $\theta = (\theta^h,\theta^i)$.

For regression ANNs, the OLS fitting
criterion is typically used, i.e., the
weights are found by minimising
\begin{equation} \label{LTSOBJANN}
  F(\theta)=\sum_{i=1}^n  (r_{i}(\theta))^2,
\end{equation}
with the residuals $r_i(\vect \theta) = f(\vect x_i; \theta) - y_i$.
Backpropagation is usually the algorithm
of choice for minimising $F$, although
Levenberg Marquardt  is also used
\cite{Hagan1994,Masters1995-advanced-book}.
In both cases, a randomly chosen initial
weighting vector $\theta_0$ is needed,
and often both methods are combined
with random start heuristic, because
$F$ has multiple local minima. The
number of local minima of $F$ grows
exponentially with the length of~$\theta$.

As in the case of linear regression, the
OLS criterion is not robust against
outliers. This can be clearly
demonstrated by replacing  one or more
data with very large or very small
values, see examples in Section
\ref{sec_expeer}. Unlike in the case of
linear regression, where the whole
regression model is shifted towards the
abnormal value, models provided by
regression ANNs exhibit wild
oscillations at the abnormal datum,
which significantly affect the rest of
the data.

There have been attempts to use more
robust Huber-type criterion (which is
used in M-estimators)
\cite{Chen1994,Liano1996}. Here the
squared residuals in (\ref{LTSOBJANN})
are replaced with $h(|r_i|)$, where $h$
is an non-negative monotone increasing
function with $h(0)=0$, whose growth
decreases with the size of the argument.
This way very large residuals have a
limited effect on the objective $F$. The
objective itself becomes more complex,
as even in linear regression, when $r_i$ depend  on $\theta$
linearly,  $F$
is the sum of quasi-convex terms, which is
not quasi-convex itself.
Huber-type functions $h$ also have another problem.
The scale of their
argument has to be either determined a
priori, or be data-dependent. In the
former case the estimator is not scale
invariant, whereas in the latter case
a robust estimation of the scale parameter $s$ is needed. When estimation of the scale $s$ is not robust (for example, taking $s$ as the mean of the absolute residuals and evaluating $h(\frac{|r_i|}{s})$),
the M-estimator will not be robust.

The LTS criterion (\ref{LTSOBJ}) was
also used in ANN training
\cite{Rusiecki2007}. It allows one to
discard up to half of the data as
potential outliers, and therefore make
ANN model robust against largely
abnormal data. Unlike in the case of
linear regression, however, the LTS
criterion may also discard good data
together with the outliers. When there
are no outliers in the data, it treats
good data as outliers, and builds wrong
regression models. We illustrate this on
several examples in Section~\ref{sec_expeer}.

\section{A new high-breakdown criterion for ANNs}
\label{sec_criterion}

In this section, we introduce a new
fitting criterion, which has allowed us
to overcome the deficiencies of the
M-type and LTS criteria mentioned in the
previous section. We propose this new
criterion, called Penalised CLTS
(PCLTS), with the following aims in mind

\begin{itemize}
\item[(A1)]
We need to discard data with unusually
large residuals as outliers.
\item[(A2)]
We need to penalise unnecessary removal
of data.
\item[(A3)]
We need to keep all data with residuals
which are comparatively small.
\item[(A4)]
For the purposes of optimising the
criterion, we need it to be based on
a Lipschitz-continuous function.
\end{itemize}

PCLTS is based on the CLTS criterion
\cite{Beliakov2012_OMS}, in which the
data are discarded if their absolute
residuals are $C$ times larger than the
median residual. The choice of $C=1$
corresponds to the LTS criterion, but
values of $C$ larger than one lead to
better efficiency of the estimator
compared to OLS in the absence of
outliers.

We propose the following objective,
which addresses the aims (A1) to (A4)
indicated above:
\begin{equation} \label{LTSOBJ1}
    F(\theta)=\sum_{i=1}^n  G(r_{(i)}(\theta)),
\end{equation}
where $r_i(\theta) = f(\vect x_i;\theta)-y_i$,
$s$ is the median of absolute
residuals $|r_1|,\ldots,|r_n|$, and
\begin{equation}\label{eq:G(t)=}
G(t) =
  \left\{
  \begin{array}{cl}
    t^2 & \mbox{ if } |t| \leq Cs,\\
    B & \mbox{ if } |t| \geq Cs(1+a),\\
    (|t|-Cs) \frac{B-(Cs)^2}{Csa} +(Cs)^2 & \mbox{ otherwise,}
  \end{array}
  \right.
\end{equation}
where $C\geq 1$ is a constant factor,
$a$ is a small positive number, and
$B \in [0,(Cs)^2]$ is a parameter
controlling the penalty for removal
of data with large residuals. Since
$$y = (x-Cs) \frac{B-(Cs)^2}{Csa} +(Cs)^2 = (Cs)^2\frac{x-Cs(1+a)}{Cs-
Cs(1+a)}+B\frac{x-Cs}{Cs(1+a)-Cs}$$
is a linear interpolation between the
points $(Cs,(Cs)^2)$ and $(Cs(1+a),B)$,
it is clear that the function $G$ is
Lipschitz-continuous.
When $C=1$, $B=0$ and $a \to 0$, we
obtain a function which coincides with
the LTS criterion almost everywhere.
When $B=(Cs)^2$, we obtain a
Huber-type function, which is
properly scaled, because the median
$s$ is a robust estimator of the scale
of the residuals (i.e., it is not
affected by very large residuals).

We will discuss the consequences of
one or another choice of the parameters
$B$ and $C$ when we discuss numerical
experiments in Section~\ref{sec_expeer}.
Before that, we address the issue of
numerical minimisation of~(\ref{LTSOBJ1}).

\section{Optimization methods} \label{sec_methods}

First, we will make a few general observations
about the objective in (\ref{LTSOBJ1}).
This function is non-convex and non-smooth,
although it is Lipschitz-continuous.
Compared to the LTS criterion (\ref{LTSOBJ}),
it has equally complex structure,
although our experiments with the CLTS
criterion in linear regression \cite{Beliakov2012_OMS}
indicated that the landscape of CLTS is
less rugged compared to that of LTS.
While the OLS criterion leads to a
complex objective when training the
standard ANNs (multiple local minima),
backpropagation seems to be quite efficient
in locating local minima of the OLS (\ref{LTSOBJANN}).
Neither backpropagation nor Levenberg Marquardt
will work with the new objective (\ref{LTSOBJ1}) directly.
Therefore we studied several alternative minimisation methods.

Specifically, we focused on the following
derivative free optimisation methods, which
have shown good potential when applied to
robust linear regression \cite{Beliakov2012_OMS}.
\begin{itemize}
\item[$\bullet$]
NEWUOA method of \cite{Powell2006},
combined with random start.
\item[$\bullet$]
derivative free bundle method
(DFBM) \cite{Bagirov2002_oms,Beliakov2007_Opt},
combined with random start.
\item[$\bullet$]
a derivative free method based on
dynamical systems
(DSO)~\cite{Mammadov2005_involume,Mammadov2007_OMS}.
\end{itemize}

Powell's NEWUOA is a derivative free
method for smooth functions, based on
quadratic model of the objective,
obtained by interpolation of function
values at a subset of $m$ previous trial
points.
While our objective is non-smooth
(although it is piecewise smooth),
we applied NEWUOA because this method
is faster than proper non-smooth optimisation
\cite{Makela_book} methods like DFBM.
Multiple local minima of the objective
mean that we have to make several starts
of the algorithm from different starting
points, where speed becomes an issue.
DSO was chosen because in previous
studies
\cite{Mammadov2005_involume,Mammadov2007_OMS,Mammadov2008}
it delivered good performance when the
objectives had a large number of variables.
We also tried other derivative free methods
(Nelder-Mead simplex method, pattern search
APPSPACK \cite{NelderMead, NR,More2009_SIAM}),
but they were not competitive and were
discarded in favour of those mentioned
above.

All mentioned methods could detect the
outliers effectively, but they were not
numerically efficient nor sufficiently
accurate when building regression
models. We explain this by a relatively
large number of variables -- the weights
of the ANN. Therefore we decided to
design a hybrid method by combining
detection and removal of outliers by
using PCLTS objective, and subsequent
training by backpropagation using
cleaned data. This approach is in line
with  Gervini and Yohai's REWLSE method
for linear regression
\cite{Gervini2002_AS}, in which  an
initial high-breakdown estimator is
improved by fully efficient OLS on
cleaned data.
This proved to deliver accurate
predictions for all data sets we have
considered.

Thus our ANN training algorithm has the following steps.
\begin{enumerate}
  \item[I.] Use one of the mentioned oprimisation methods with PCLTS objective (\ref{LTSOBJ1}).
  \item[II.] Clean the data by removing the data with the residuals larger than $Cs$.
  \item[III.] Use standard backpropagation with the OLS objective to train ANN on clean data.
\end{enumerate}

The most time consuming step here is
Step~I. It requires multiple evaluations
of the objective at a large number of
points (typically of order of $100,000$ in our
experiments). Step~III is standard in
ANN training, and it was executed using
standard ANN training library with
default parameters (in our studies we
used \emph{nnet} package in \emph{R}
software \cite{Rsoft}). It took
negligible CPU time compared to Step~I.
Step~II is trivial.

We used implementations of NEWUOA, DFBM
and DSO in C++ language (DFBM and DSO
taken form  GANSO library
\cite{Beliakov2007_Opt}, and a
translation of NEWUOA from Fortran to
C). All methods have few overheads, and
the main complexity was in evaluating
the objective. In order to achieve
competitive CPU time, we parallelised
calculation of the objective and
offloaded it to a graphics processing
unit (GPU), NVIDIA's Tesla C2050
\cite{NVIDIA_Tesla_C1060,NVIDIA_CUDA4}.
General Purpose GPUs have recently
become a valuable alternative to
traditional CPUs and CPU clusters. Tesla
C2050 has 448 GPU cores and 3 GB of RAM,
and can execute thousands of threads at
a time.  GPUs have limitations too, in
particular it is Single Instruction
Multiple Threads (SIMT) paradigm, which
means that instructions in different
threads (of the same thread block) need
to be almost identical. For parallel
calculation of the residuals this is not
an issue, as this task is trivially
parallel, and is executed using
\emph{for\_all} primitive
\cite{Thrust,S.Sengupta2007}.
Calculation of the median $s$ is done in
parallel using either GPU parallel
sorting
\cite{Sintorn2008_JPDC,Cederman2009_ACM}, or
a specific GPU selection algorithm
\cite{Beliakov2011_arxiv}. Summation is
performed in parallel using GPU
reduction
\cite{Thrust,S.Sengupta2007,NVIDIA_Reduction}.

\section{Numerical experiments} \label{sec_expeer}

\subsection{Data sets}
\label{sec_datasets}

We generated several artificial data
sets using test functions
considered in several previous papers
on robust ANNs, and also used several
real world data sets to which we
introduced outliers. The artificial data
sets are described below in Data Sets~1
through to~10. In all of these examples,
the independent variables $\vect x$ were
considered in $\R^m$, where the
dimension $m=p-1$ took on the values
$m=1,2,3,5,10$. The values of
$\vect x$ were chosen uniformly
at random in the segment $[\min_x,\max_x]$
indicated below in describing each data
set of the corresponding example. The
target variable was determined using
the formula
\begin{equation}\label{eq:testmodel}
y=h(\vect x)+\varepsilon,
\end{equation}
where
$\varepsilon \sim N(0,{\rm noiselevel})$.
In all data sets, the noise level adopted
the values $0, 0.1, 0.2$.
The sizes $n$ of the samples in each
data set were taken as $n\in \{100,500,5000\}$.
Then we replaced a proportion $\delta$
of the points with five subsets of outliers.
The outliers were divided
into 5 subsets of equal size. Each subset
was centered at a point with the first
$m$ coordinates chosen uniformly at
random in the same segment as
the explanatory variables, and
with the last coordinate equal to 10,000 + N(0,0.01).
For Data Sets~1 to~10 indicated
below, we examined all combinations of
the dimension $n = 1,2,3,5,10$, proportion
of outliers $\delta = 0,0.1,0.2$,
and noise level $\varepsilon=0, 0.1, 0.2$.

\textbf{Data Set 1} uses the
model~(\ref{eq:testmodel}) with the function
$h$ given by the equation
\begin{equation}\label{eq:dataset1}
h(\vect x) = ||\vect x||^{2/3},
\end{equation}
where $||\vect x||$ stands for the Euclidean
norm of $\vect x$ and where the explanatory variables are
again chosen uniformly at random in
the segment $[\min_x,\max_x]=[-2,2]$.
Earlier, the robustness
of ANN for this function and data with
outliers was investigated in
\cite{Chuang2004}, \cite{Shieh2009}
and \cite{Melegy2009} in the case
where $m=1$ and $n=501$.


\textbf{Data Set 2} deals with the
model~(\ref{eq:testmodel}) using the function
$h$ given by the equation
\begin{equation}\label{eq:dataset2}
h(\vect x) = x_1 e^{||\vect x||},
\end{equation}
where the explanatory variables are
again chosen uniformly at random in
the segment $[\min_x,\max_x]=[-2,2]$.
This function was earlier considered
in~\cite{Melegy2009}.

\textbf{Data Set 3} is defined by the
model~(\ref{eq:testmodel}) using the function
$h$ given by the formula
\begin{equation}\label{eq:dataset3}
\mbox{sinc}(\vect x) =
\frac{\sin(||\vect x||)}{||\vect x||},
\end{equation}
where the explanatory variables are
again chosen uniformly at random in
the segment $[\min_x,\max_x]=[-10,10]$.
The function~(\ref{eq:dataset3})
was considered in numerous articles.
For example, it was investigated in
\cite{Chen1994}, \cite{Chuang2004}
and~\cite{Smola1998}.

\textbf{Data Set 4} uses the
model~(\ref{eq:testmodel}) with the function
$h$ given by equation
\begin{equation}\label{eq:dataset4}
h(\vect x) =
\sin\left(\frac{5}{m}\sum^m_{i=1}x_i\right)
\,\mbox{acos}\left(\frac{1}{m}\sum^m_{i=1}x_i\right)
\cos\left(\frac{3}{m}\sum^m_{i=1}x_i-2/n\right),
\end{equation}
where $\vect x \in [-1,1]^m$
and $n$ is the sample size.
Earlier, the robustness of ANN for
this function was investigated in
\cite{Chuang2000} and~\cite{Chuong2007}.

\textbf{Data Set 5} is given by
the model~(\ref{eq:testmodel}) with the
following function
\begin{equation}\label{eq:dataset5}
h(\vect x)=\sin(10\pi||\vect x||)
+\sin(20\pi||\vect x||),
\end{equation}
where the explanatory variables
are uniformly chosen in the segment
$[\min_x,\max_x]=[0, 0.3]$. This
function was investigated in
\cite{Chuong2007}, \cite{Jeng2000}
and~\cite{Pati1993}.

\textbf{Data Set 6} deals with
the model~(\ref{eq:testmodel}) given by
the following function
\begin{equation}\label{eq:dataset6}
h(\vect x)=(x_1^2-x_2^2+x_3^2-x_4^2+\cdots)
\sin(0.5(x_1+x_3+\cdots)),
\end{equation}
where the explanatory variables
are uniformly chosen in the segment
$[\min_x,\max_x]=[-2,2]$.
This function was investigated in
\cite{Jeng2000} and~\cite{Chuong2007}.

\textbf{Data Set 7} is defined by
the model~(\ref{eq:testmodel}) with the
function $h$ given by the formula
\begin{equation}\label{eq:dataset7}
h(\vect x) = \frac{\sin(x_1+x_3+\cdots)}{x_1+x_3+\cdots}
\frac{\sin(x_2+x_4+\cdots)}{x_2+x_4+\cdots},
\end{equation}
where $\vect x \in [-5,5]^m$. Earlier,
this function was considered in
\cite{Chuong2007}, for only one value
of~$m=2$.

\textbf{Data Set 8} is given by
the model~(\ref{eq:testmodel}) with the
function $h$ defined by the equation
\begin{equation}\label{eq:dataset8}
h(\vect x) = 0.2(x_1x_2\cdots x_m)+
1.2\sin(||\vect x||^2),
\end{equation}
where $\vect x \in [-1,3]^m$. Earlier,
this function was investigated
in~\cite{Chuong2007}.

\textbf{Data Set 9} uses the function
$h$ for the model~(\ref{eq:testmodel})
defined by the equation
\begin{equation}\label{eq:dataset9}
h(\vect x) = \max\{ \exp(-10 x_1^2),
\exp(-50x_2^2), 1.25\exp(-5(||\vect x||^2))\},
\end{equation}
where $\vect x \in [-2, 2]^m$.

\begin{figure}
\centering
 \includegraphics[width=0.5\textwidth]{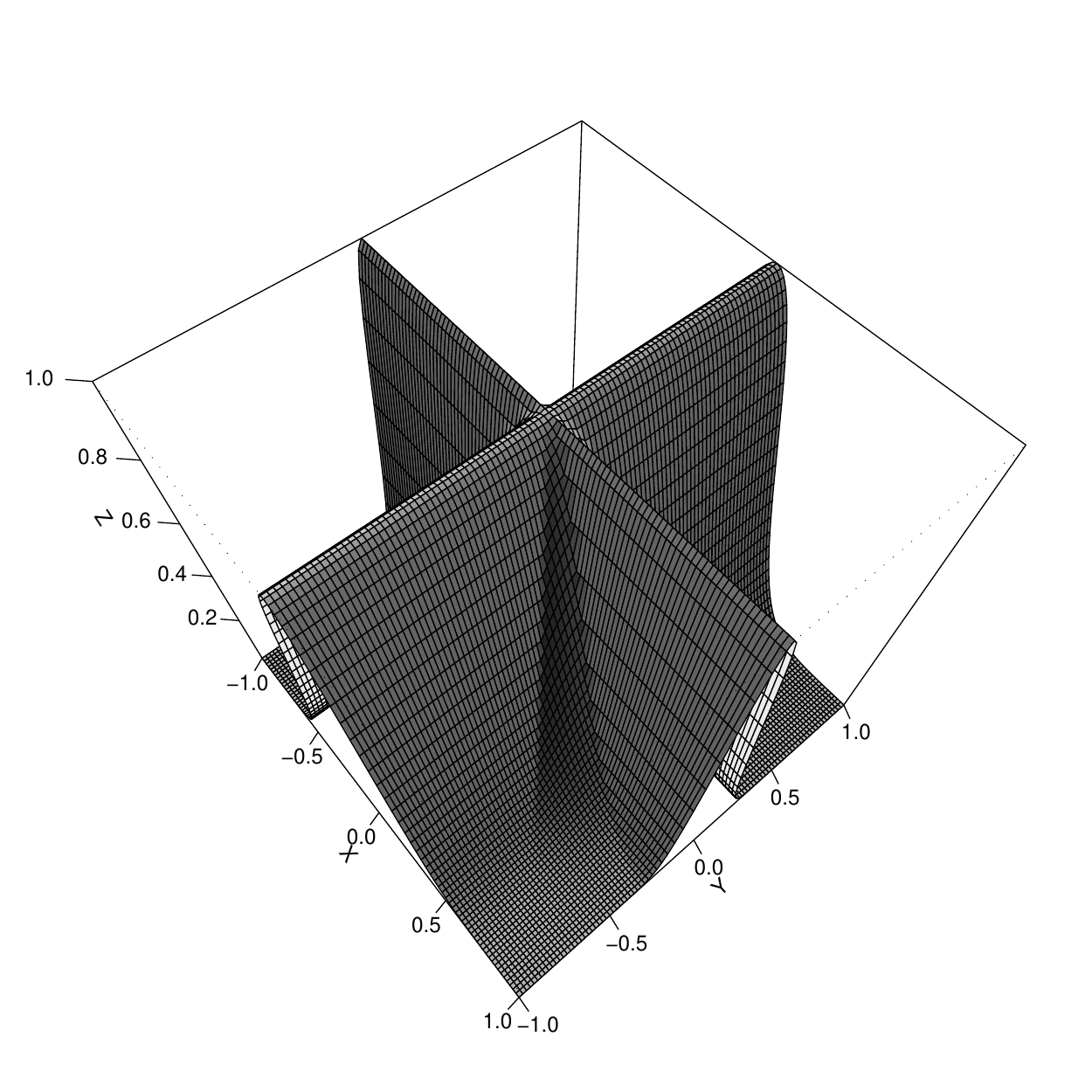}
\caption{Plot of the surface of the function $h$ used in Data set 9.}
\label{fg:dim2_model9}
\end{figure}

This function has an interesting surface
plot in two variables, illustrated in
Figure~\ref{fg:dim2_model9}.

\textbf{Data Set 10} is defined by
the model~(\ref{eq:testmodel}) with the
following function
\begin{equation}\label{eq:dataset10}
h(\vect x) = 0.5 ||\vect x|| \sin(||\vect x||) + \cos^2(||\vect x||),
\end{equation}
where $\vect x \in [-6, 6]^m$.
It was considered in \cite{Lee1999}.

\textbf{Data Set 11} is a real
world data set from the standardized
benchmark collection for neural networks
{\sc Proben~1} publicly available from
\cite{proben1}. {\sc Proben~1} was
considered in many articles, for example,
\cite{Igel2003,Ilonen2003,Prechelt94,Prechelt98,Tomandl2001}.
The Data Set~11 is the set `building'
with 4208 instances, where the value
of hourly electricity consumption is
regarded as a function of 14 input
variables and has to be predicted. Outliers have been introduced by replacing several values of the dependent variable with very large values (10,000).


\textbf{Data Set 12} is also a real
world data set from {\sc Proben~1}.
It uses the values of output variable
`cold water consumption' viewed as a
function of 14 input variables, see
\cite{proben1,Prechelt94}.

\textbf{Data Set 13} is a real world
data set from {\sc Proben~1} using
the values of output variable
`hot water consumption' viewed as a
function of 14 input variables, see
\cite{proben1,Prechelt94}.

\textbf{Data Set 14} is a real
world data set `hearta' from
{\sc Proben~1}, see \cite{proben1,Prechelt94}.

\subsection{Parameters of the algorithms}

We varied two parameters of the PCLTS
objective: $C$ and $B$. We took $C \in
\{1,8\}$. The value $C=1$ (with $B=0$)
corresponds to the LTS criterion, and it
was interesting to observe the
difference larger values of $C$ made.
The value of the second parameter $B$
were chosen in $\{0,1,8\}$. Here we
studied the influence of the penalty for
removing the data when cleaning the data
set.

As far as the parameters of optimisation
algorithms are concerned, we fixed them
in order to give each method
approximately the same CPU budget as the
others. We fixed the number of random
starts  of NEWUOA and DFBM at 200.
DFBM ran with the maximum number of iterations set to 1000.

We used \emph{nnet} package in \emph{R}
to perform the final stage of ANN
training with cleaned data with default
parameters.

\subsection{Benchmarking criteria}

As it is customary in data analysis
literature, we used the Root Mean
Squared Error (RMSE) to measure the
quality of approximation. The data were
split into training and test subsets.
The training subset contained noisy data
with or without outliers. The test
subset contained noiseless data, i.e.
the accurate values of the test
functions $h$ at uniformly
distributed data within the domain
of each function we considered.
The size of each test sample was
equal to that of the corresponding
training sample. The test data set
was generated separately and was not
provided to the training algorithm.
We report the RMSE on the test data
set, which is the most important
characteristic, reflecting ANN's
generalisation ability.

We also report the average CPU time
taken by the PCLTS and the nnet
procedures. The CPU time did not depend
much on the test function, nor the
Gaussian noise level or proportion of
outliers, but rather on the dimension of
the problem. Therefore we averaged the
CPU time over different
experiments.

\subsection{Results and discussion}

In the Appendix we presented detailed results for each of the test function we considered. In Tables \ref{t:rmse1}--\ref{t:rmse10} we compare RMSE values for PCLTS and standard backpropagation (nnet). The values in bold indicate the cases the method  broke down.

\begin{figure}
\centering
\subfloat[noise level 0\%]{\label{fg:c1-d1_020_500_p2_f8_b8_m4_i200_n0_h10_nnet_test_predictions}
 \includegraphics[width=0.3\textwidth]{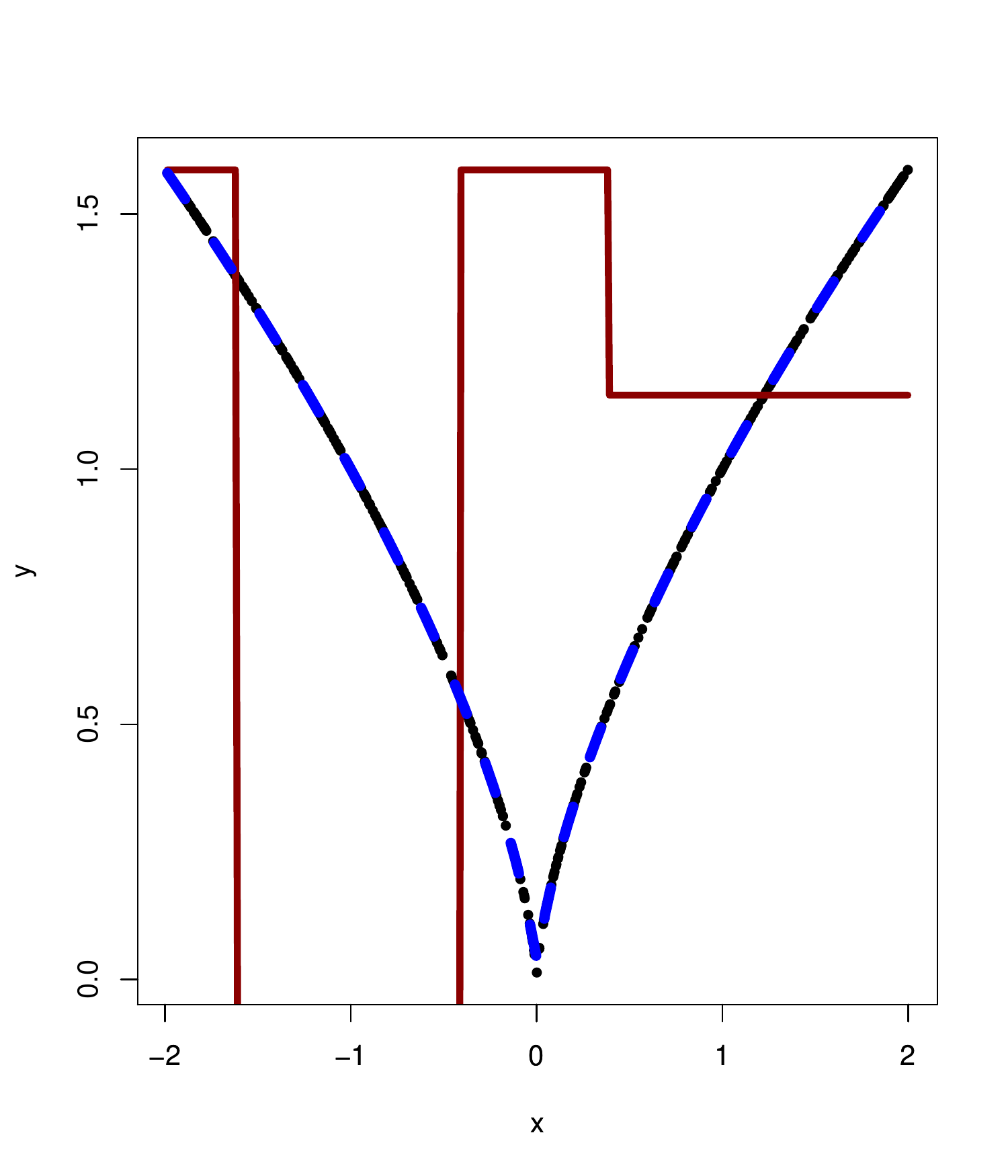}}
\subfloat[noise level 10\% ]{\label{fg:c1-d1_020_500_p2_f8_b8_m4_i200_n10_h10_nnet_test_predictions}
 \includegraphics[width=0.3\textwidth]{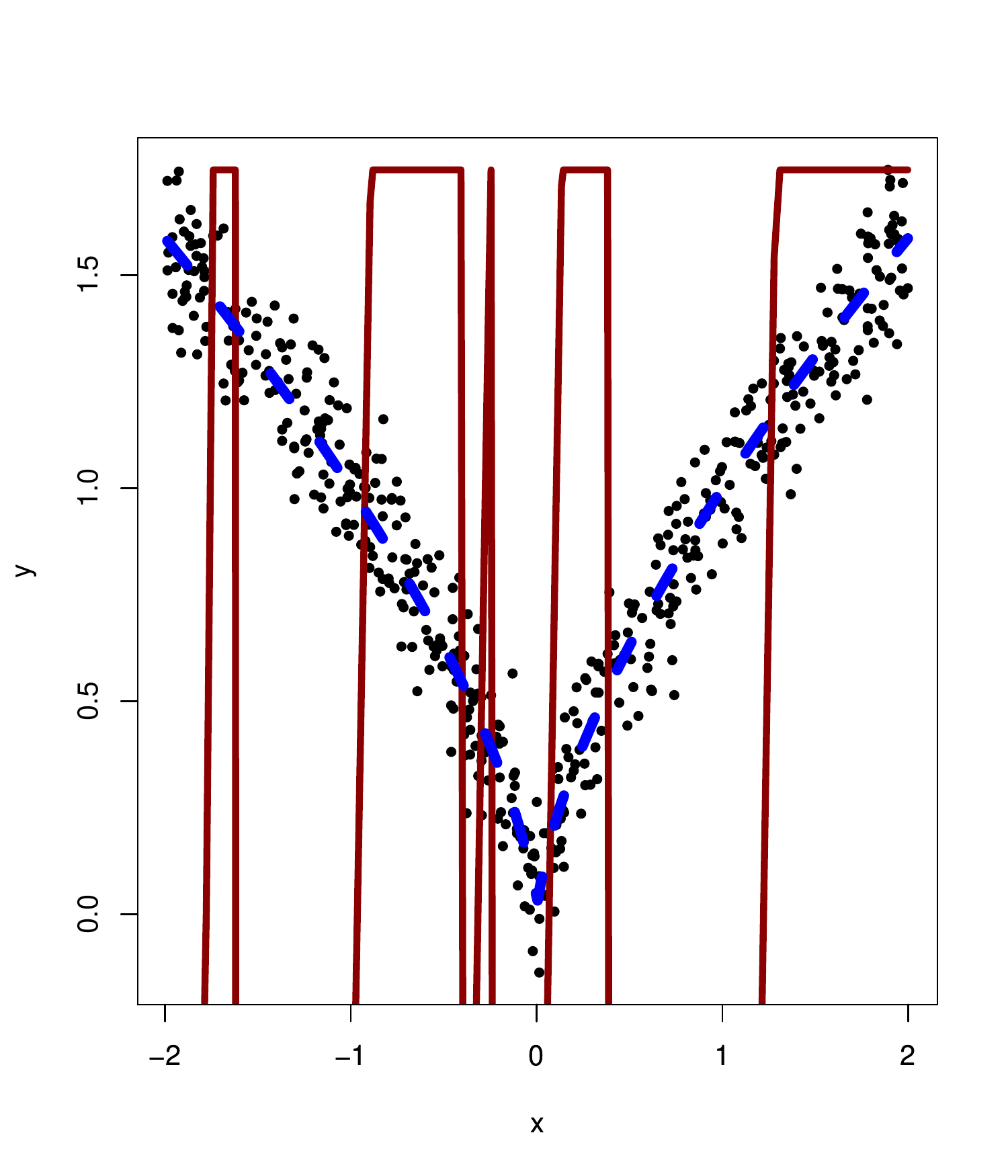}}
\subfloat[noise level 20\%]{\label{fg:c1-d1_020_500_p2_f8_b8_m4_i200_n20_h10_nnet_test_predictions}
 \includegraphics[width=0.3\textwidth]{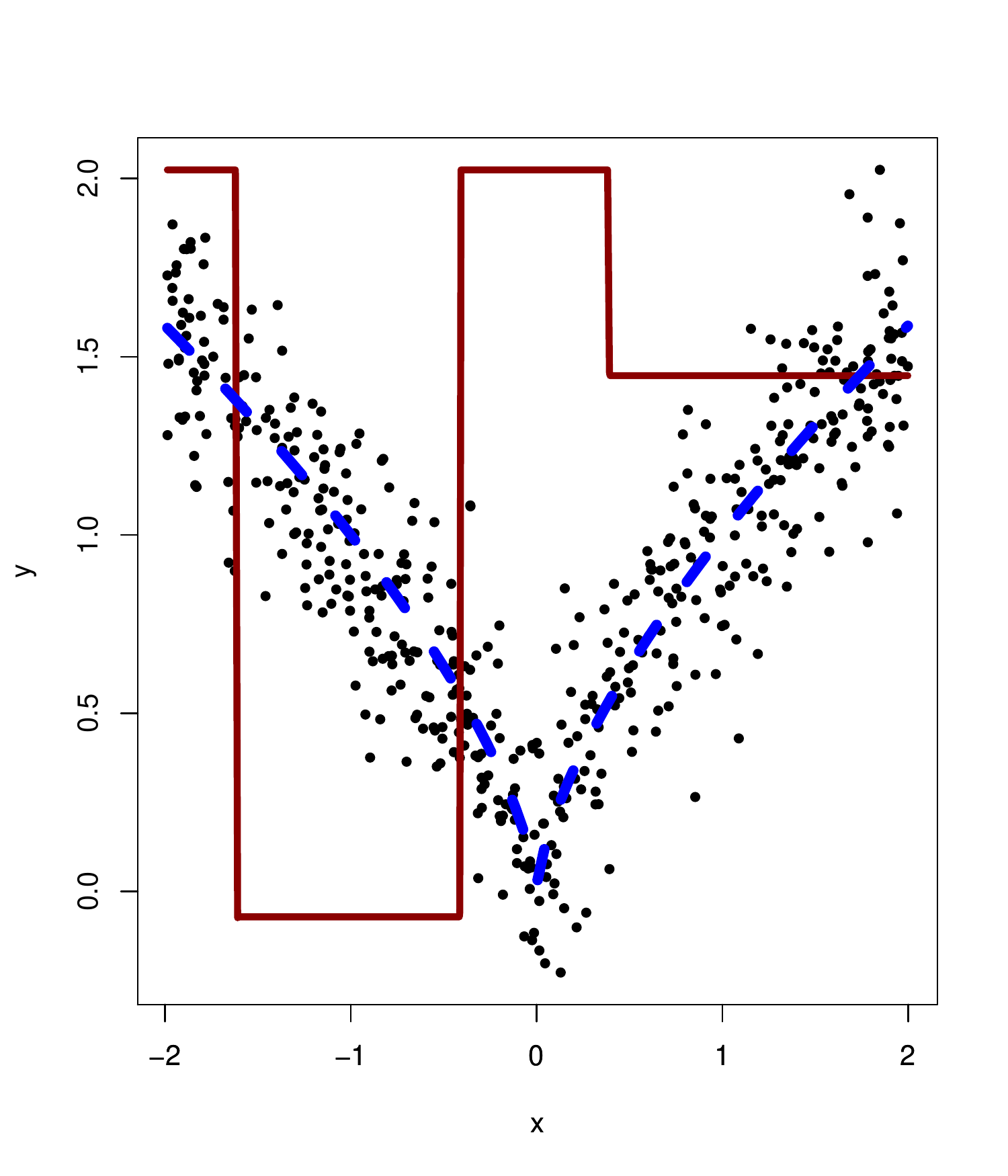}}
\\
\subfloat[noise level 0\%]{\label{fg:c1-020_500_f1_b0_m4_i180_n0_h10}
 \includegraphics[width=0.3\textwidth]{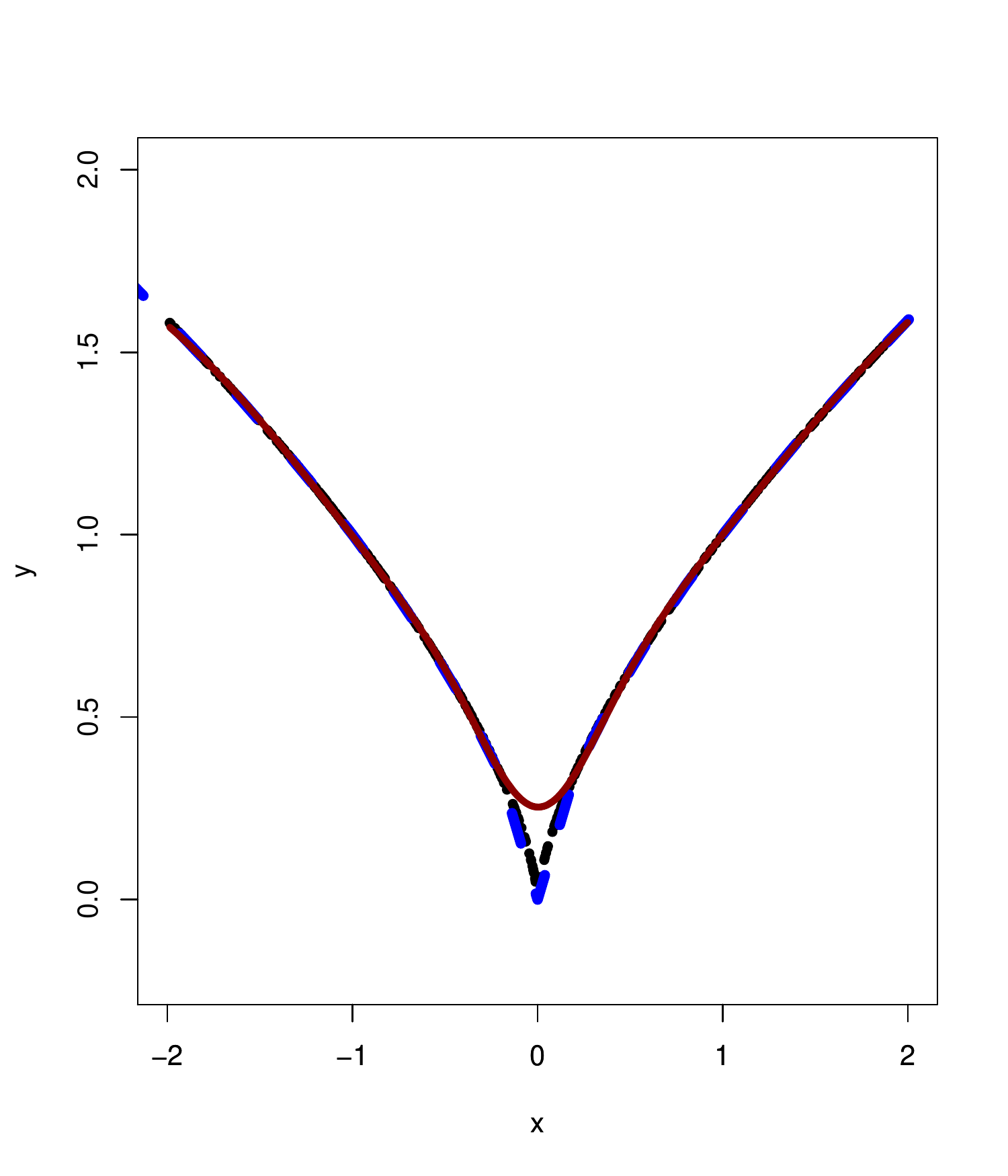}}
\subfloat[noise level 10\%]{\label{fg:c1-020_500_f1_b0_m4_i180_n10_h10}
 \includegraphics[width=0.3\textwidth]{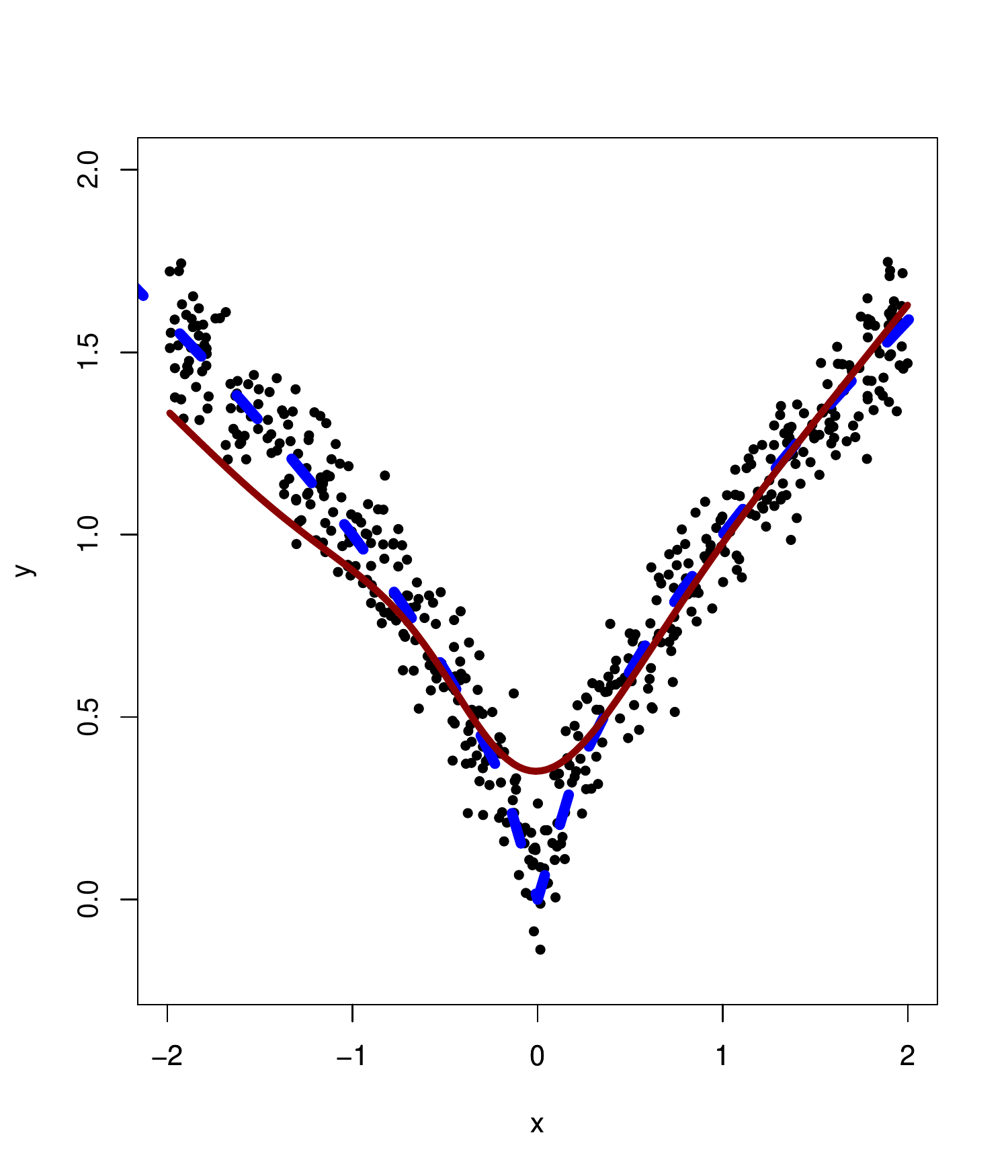}}
\subfloat[noise level 20\%]{\label{fg:c1-020_500_f1_b0_m4_i180_n20_h10}
 \includegraphics[width=0.3\textwidth]{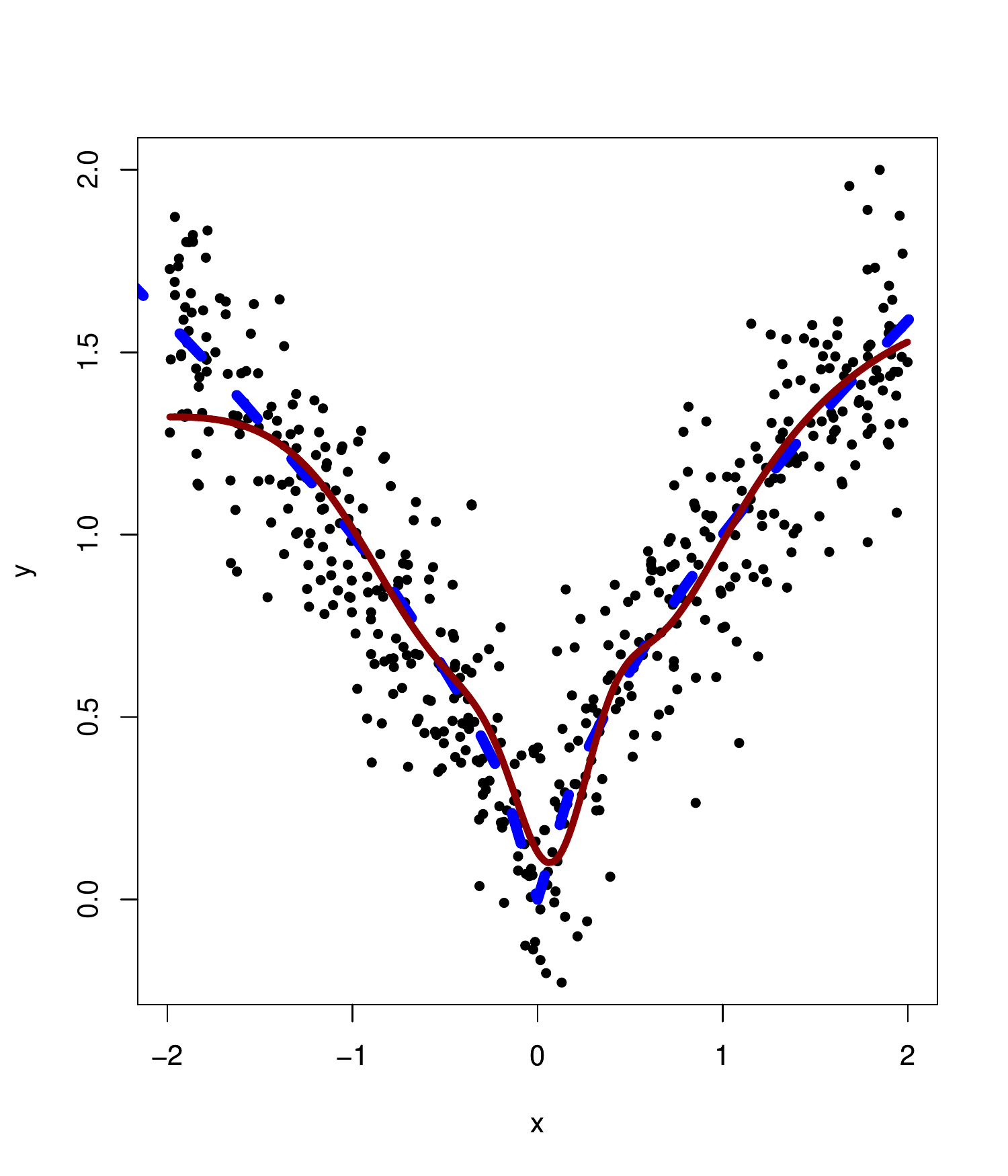}}

\caption{Data set 1, standard ANN (above) versus LTS (below).
The proportion of outliers is $\delta=0.2$, and $n=500$. Gaussian noise in the data is 0, 0.1 and 0.2 respectively. Red solid curve is the ANN prediction, the blue dashed curve is the (noiseless) test data, and black dots are training data. The outliers are not shown.}
\label{fg:ds1nnet}

\end{figure}

\begin{figure}
\centering
\subfloat[noise level 10\%]{\label{fg:c8-d1_020_500_p2_f8_b8_m4_i200_n10_h10_cuda_nnet_short_train_predictions}
 \includegraphics[width=0.3\textwidth]{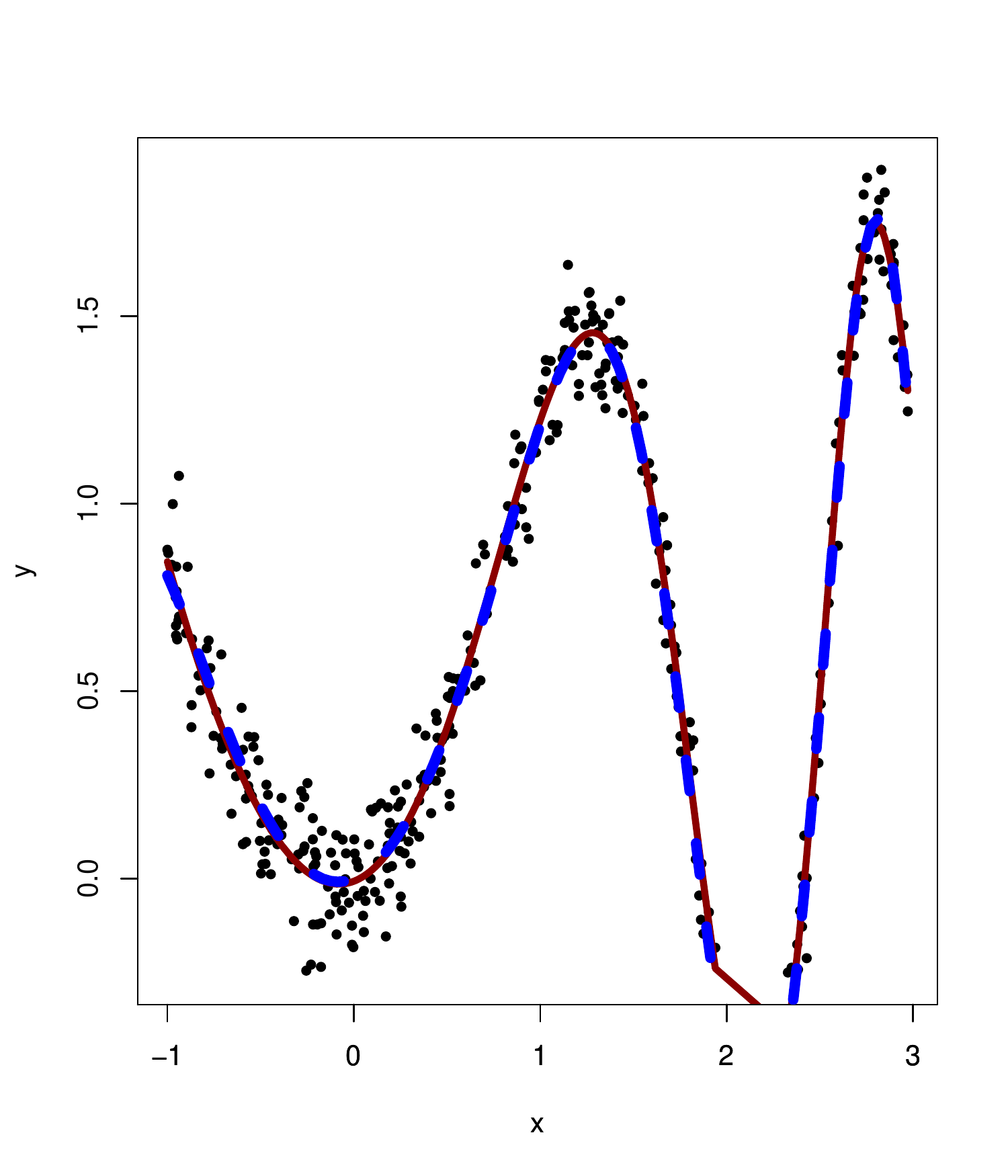}}
\subfloat[noise level 20\%]{\label{fg:c8-d1_020_500_p2_f8_b8_m4_i200_n20_h10_cuda_nnet_test_predictions}
 \includegraphics[width=0.3\textwidth]{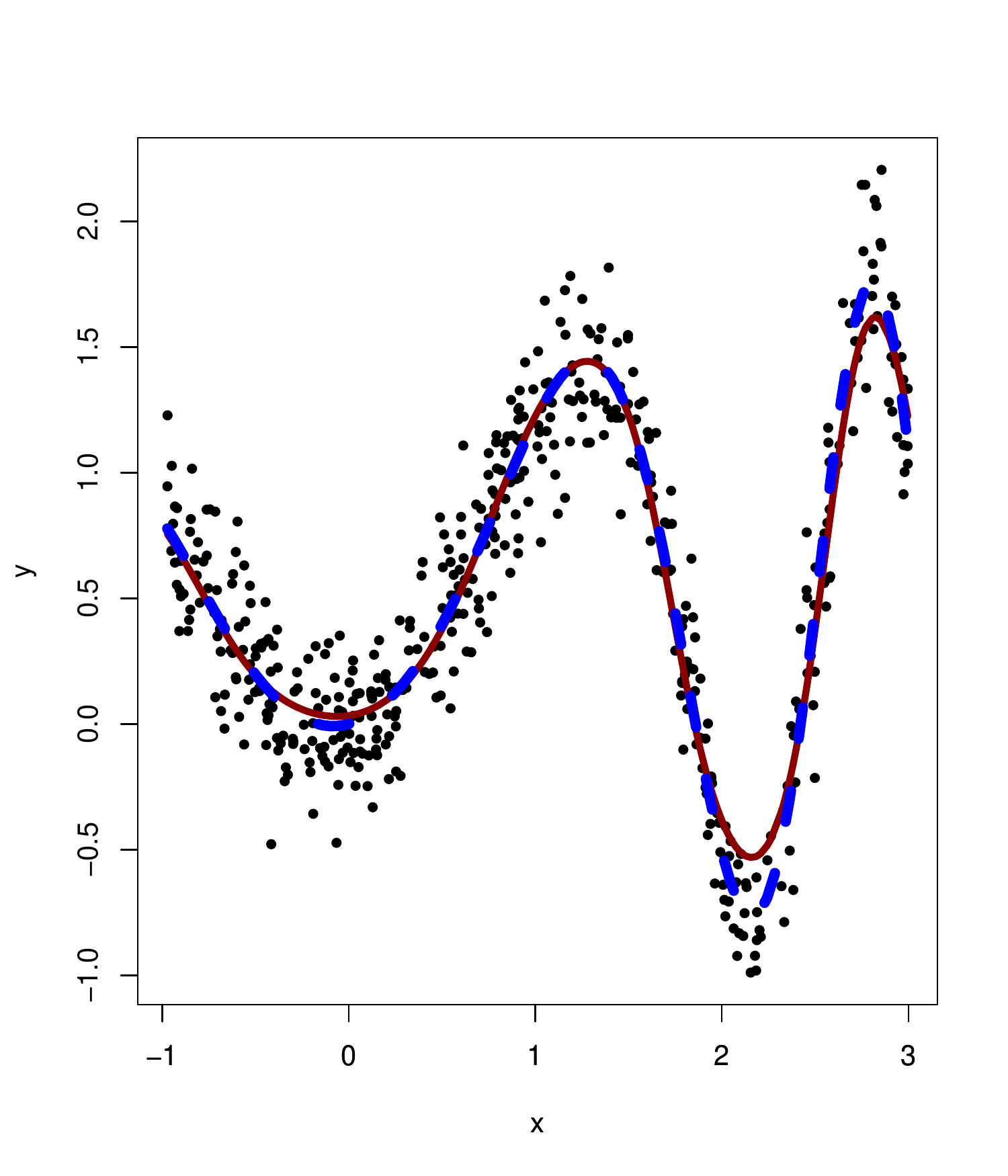}}
\subfloat[prediction by nnet]{\label{fg:c8-d1_020_500_p2_f8_b8_m4_i200_n20_h10_nnet_test_predictions}
 \includegraphics[width=0.3\textwidth]{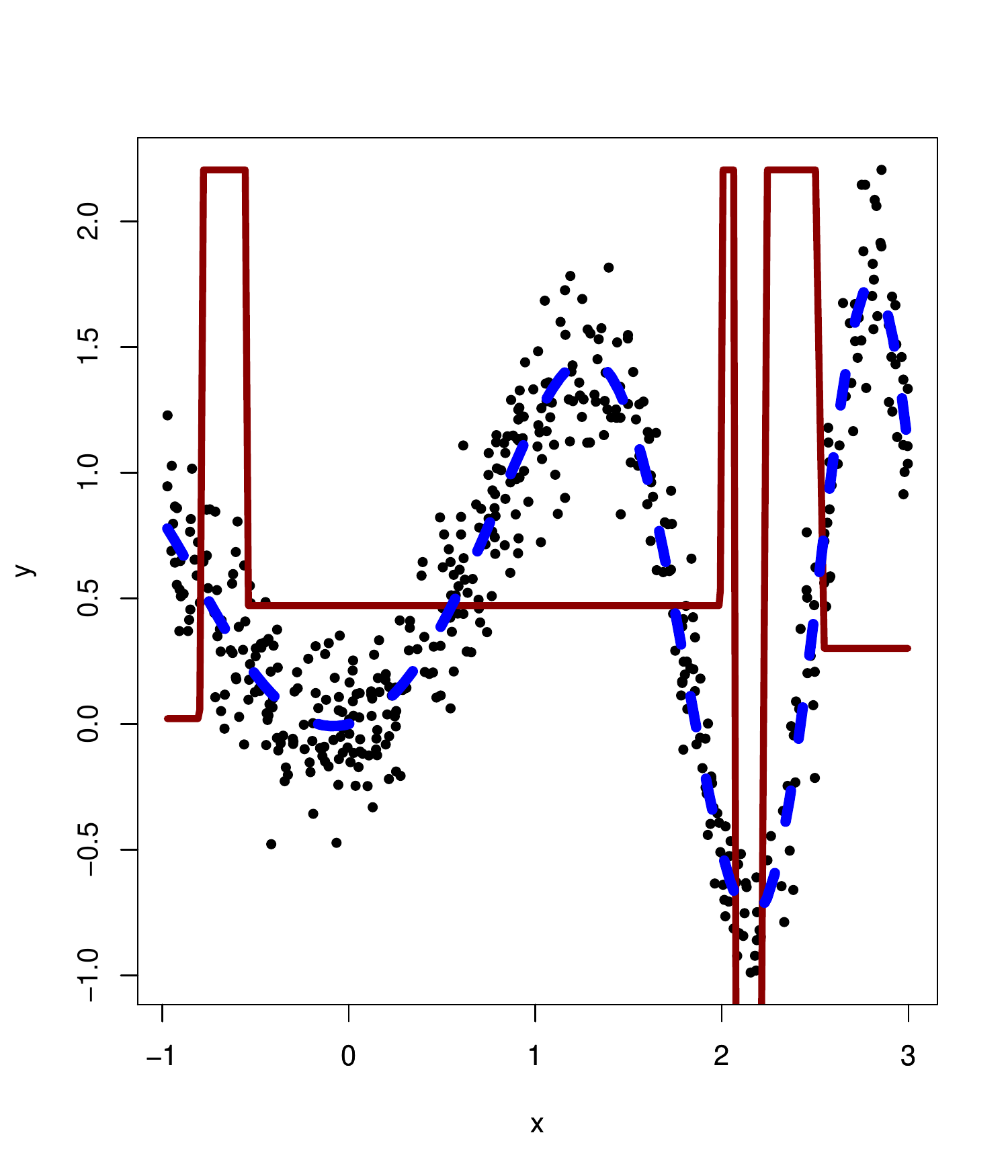}}
\caption{Data set 8, $\delta=0.2$, $n=500$. PCLTS ($C=8$, $B=8$) correctly predicts the model, but backpropagation fails.}
\label{fg:ds8}
\end{figure}

Two facts are apparent from these tables. First, it is the failure of the traditional ANN training to predict the test data sets. This can be viewed in Figure \ref{fg:ds1nnet}, top row, where the predictions are not even close to the majority of the train data, nor the test data. Similar pictures were obtained for other data sets, e.g. Figure \ref{fg:ds8},c). This is reflected in very large RMSE in the tables, which is consistent across the tables. Second, it is the ability of our method based on the PCLTS objective to filter out the outliers. We note that the RMSE are practically the same as RMSE of backpropagation \emph{in the absence of outliers} (the values in nnet column not in boldface). This indicates that all outliers were filtered out, and only few, if any, good data were removed. So our method achieves the same accuracy when using contaminated data, as nnet when using only clean data, as if the outliers were not there. This is repeated across all data sets we used.

Furthermore, when there are no outliers in the data, our method delivers the same RMSE as backpropagation, which means that no (or almost no) data were unnecessarily removed.

Table \ref{t:all-cpu-large1} show the CPU time of PCLTS and backpropagation. We see that  PCLTS takes 1,000-10,000 times longer to train the ANN, compared to backpropagation. We do not observe a clear pattern of rising CPU cost with the increased size of the data set. This can be explained, on one hand, by a smaller number of iterations of the optimisation algorithm for larger data, as the objective becomes less rough, and by the use of GPU to offload calculation of the residuals. We have used 448 cores for these calculations, and for 5000 data GPU calculations did not reach the saturation point. While as expected, PCLTS takes much longer, this time is compensated by the quality of the result. On the other hand, the CPU appears not to be excessive.

Tables \ref{t:rmse11}--\ref{t:rmse14} contain RMSE and CPU time for the real world data sets we considered. Again, we observe that standard backpropagation breaks down, while our method based on PCLTS objective provides a very good fit, comparable to backpropagation in the absence of outliers.
CPU times are consistent with those reached in the examples with artificial data sets.

\begin{figure}
\centering
\subfloat[noise level 0\%]{\label{fg:c1-020_500_f1_b0_m4_i180_n0_h10}
 \includegraphics[width=0.3\textwidth]{plot-c1-cuda-020_500_f1_b0_m4_i180_n0_h10}}
\subfloat[noise level 10\%]{\label{fg:c1-020_500_f1_b0_m4_i180_n10_h10}
 \includegraphics[width=0.3\textwidth]{plot-c1-cuda-020_500_f1_b0_m4_i180_n10_h10}}
\subfloat[noise level 20\%]{\label{fg:c1-020_500_f1_b0_m4_i180_n20_h10}
 \includegraphics[width=0.3\textwidth]{plot-c1-cuda-020_500_f1_b0_m4_i180_n20_h10}}
\\
\subfloat[noise level 0\%]{\label{fg:c1-020_500_f8_b8_m4_i200_n0_h10}
 \includegraphics[width=0.3\textwidth]{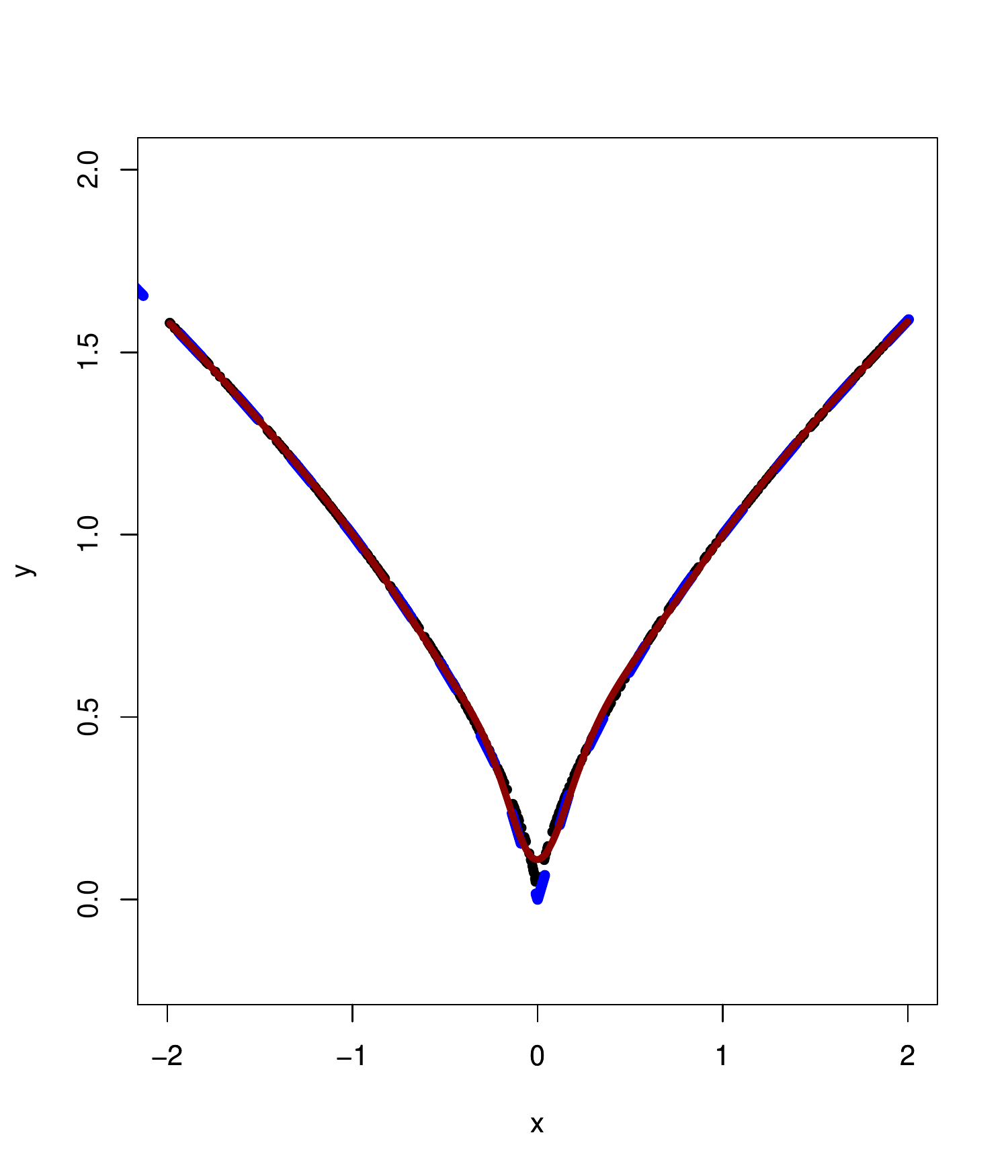}}
\subfloat[noise level 10\%]{\label{fg:-c1-020_500_f8_b8_m4_i200_n10_h10}
 \includegraphics[width=0.3\textwidth]{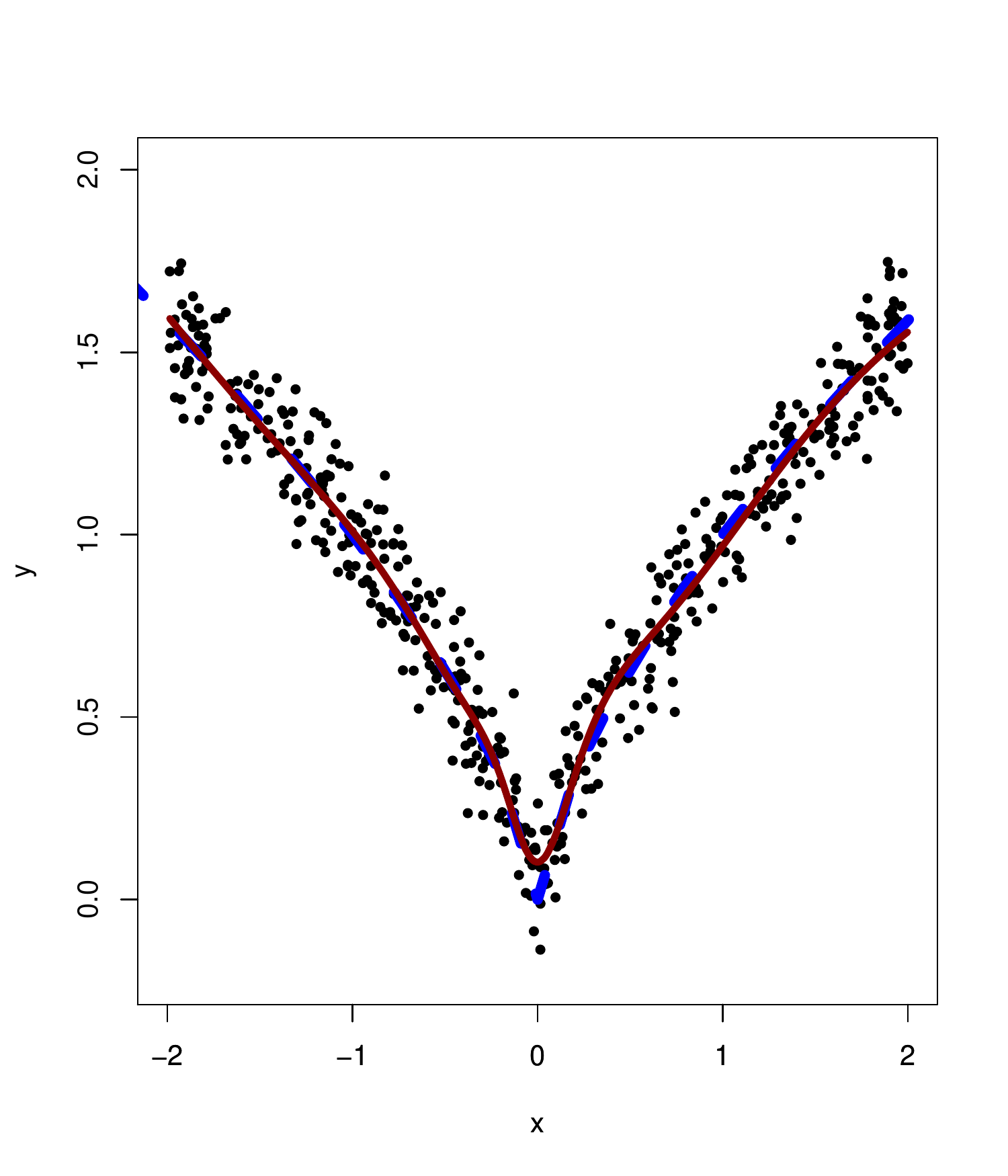}}
\subfloat[noise level 20\%]{\label{fg:c1-020_500_f8_b8_m4_i200_n20_h10}
 \includegraphics[width=0.3\textwidth]{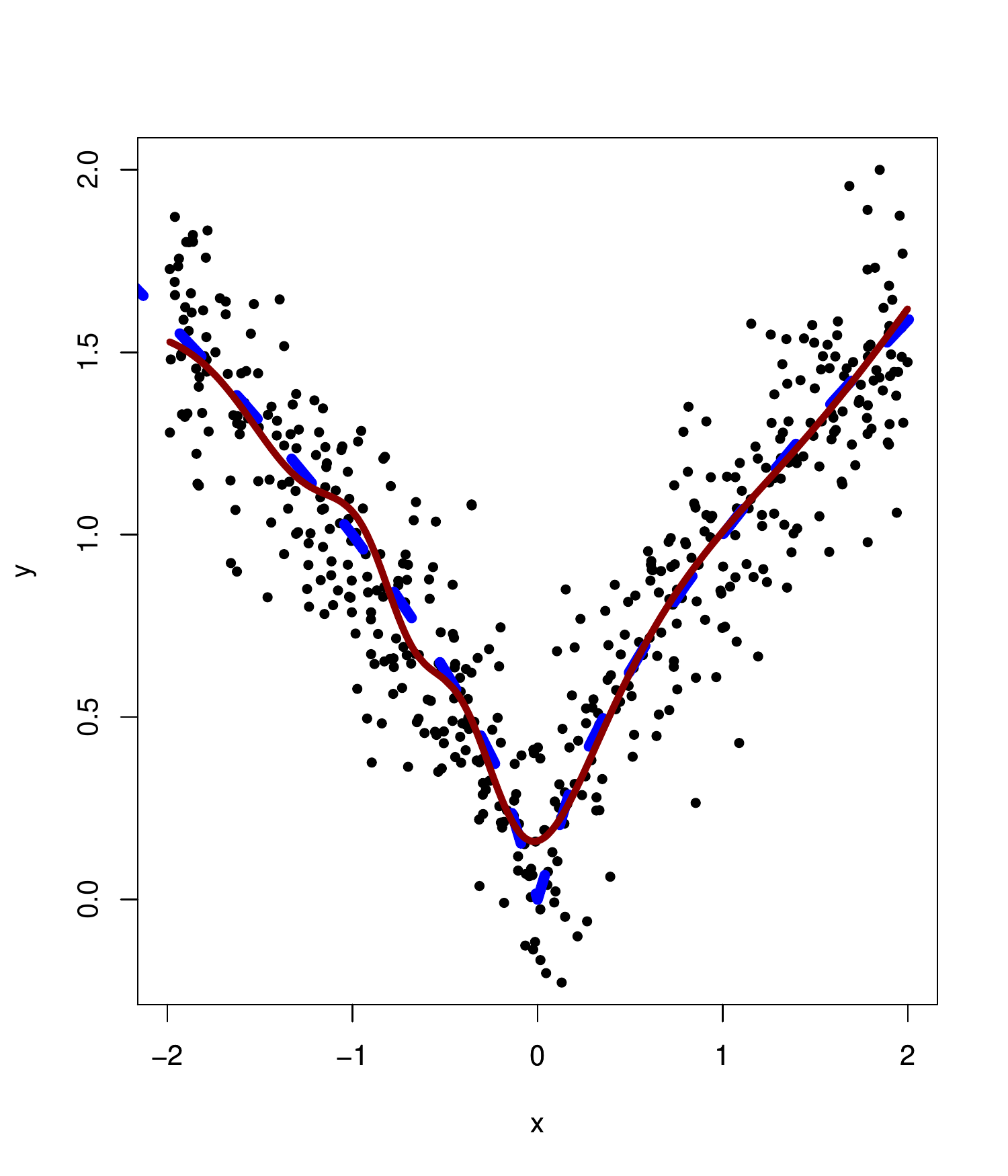}}
\caption{Data set 1, The effect of using PCLTS objective with parameters $C=8$ and $B=8$ (below) versus  LTS objective (above). The proportion of outliers is $\delta=0.2$, and $n=500$. With the LTS objective, some data are unnecessarily removed as outliers, and the ANN model is fitted incorrectly in these regions, whereas PCLTS objective avoids this.}
\label{fg:ds1trainPCLTS}
\end{figure}

\begin{figure}
\centering
\subfloat[noise level 0\%]{\label{fg:c3-cuda-test-000_500_f1_b0_m4_i180_n0_h10}
 \includegraphics[width=0.3\textwidth]{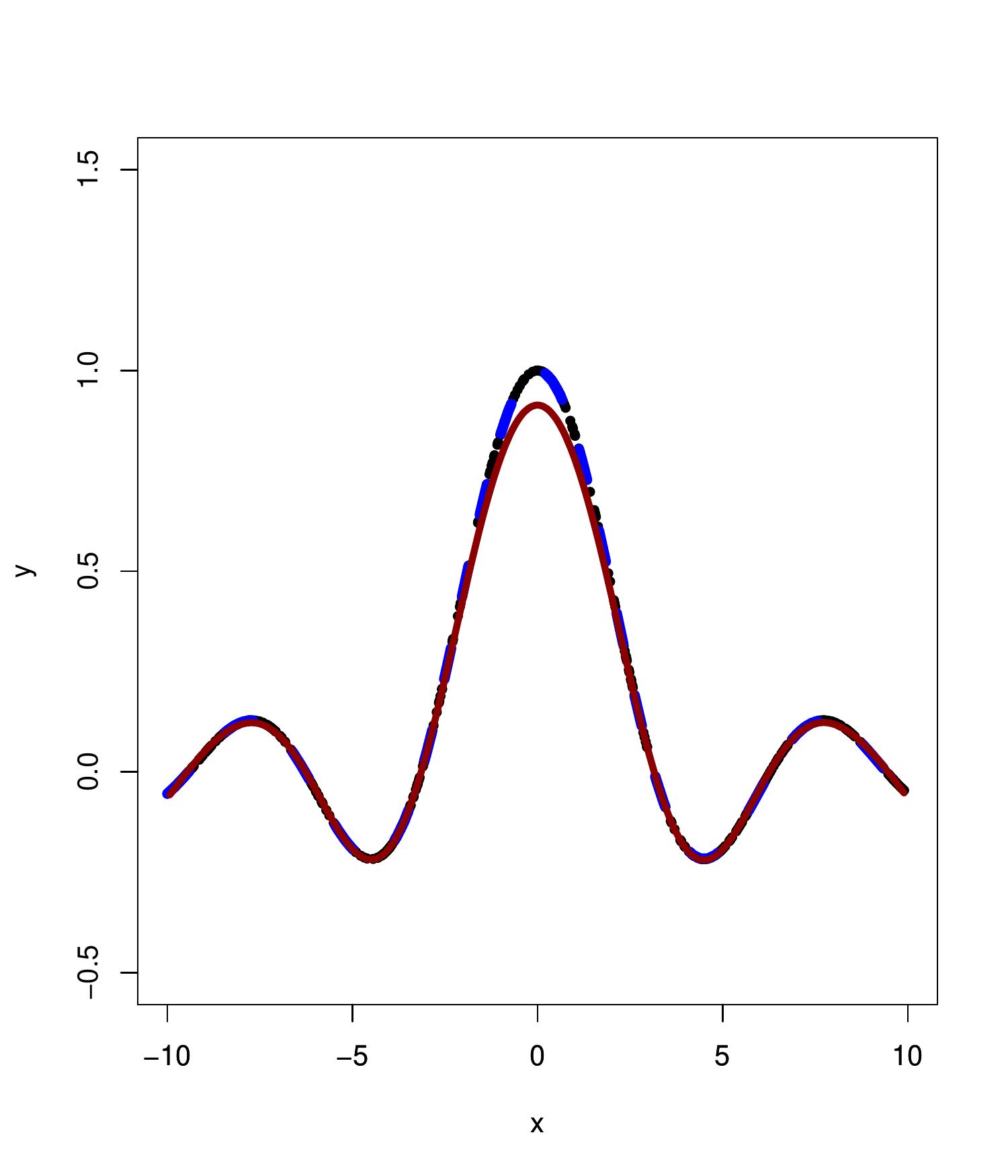}}
\subfloat[noise level 10\%]{\label{fg:c3-cuda-test-000_500_f1_b0_m4_i180_n10_h10}
 \includegraphics[width=0.3\textwidth]{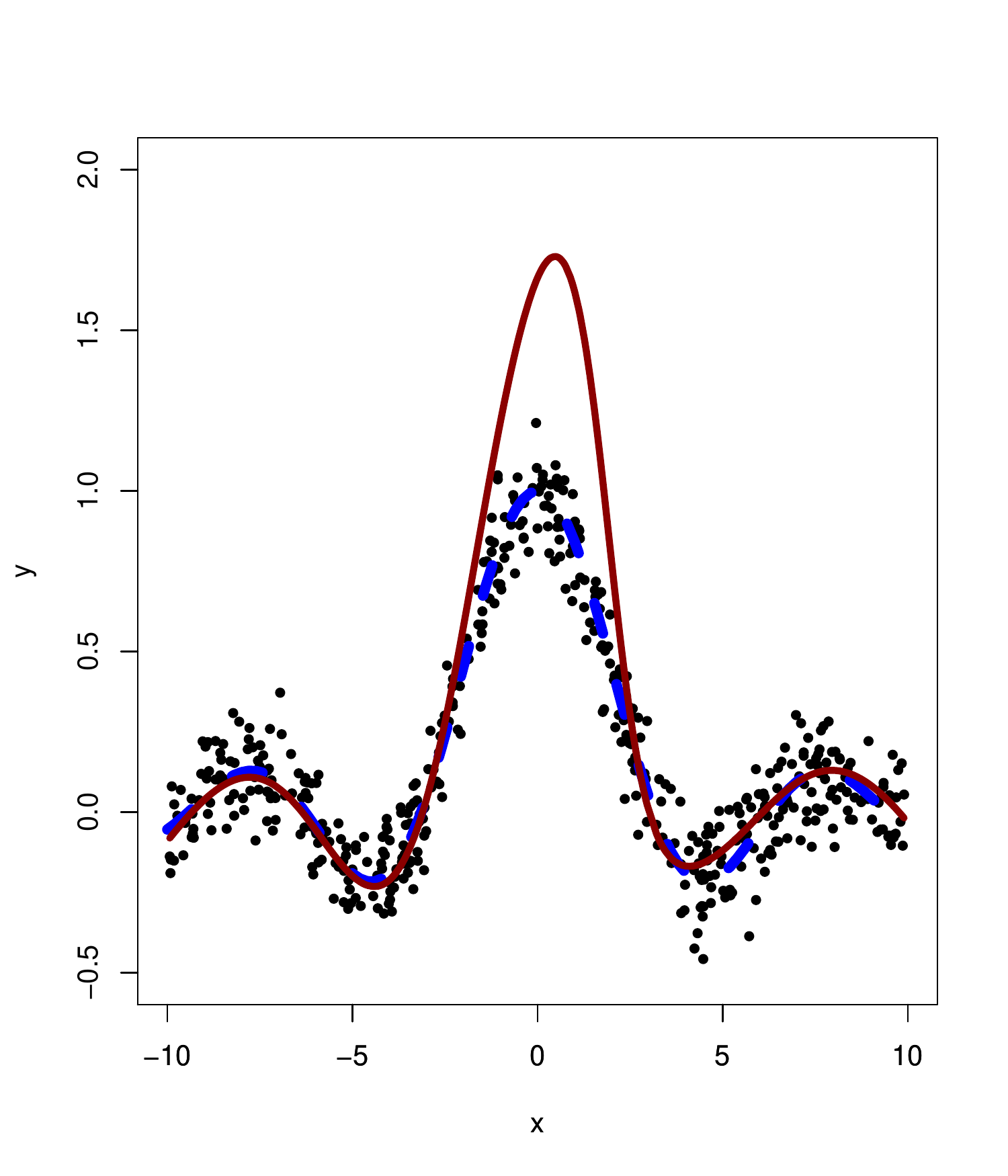}}
\subfloat[noise level 20\%]{\label{fg:c3-cuda-test-000_500_f1_b0_m4_i180_n20_h10}
 \includegraphics[width=0.3\textwidth]{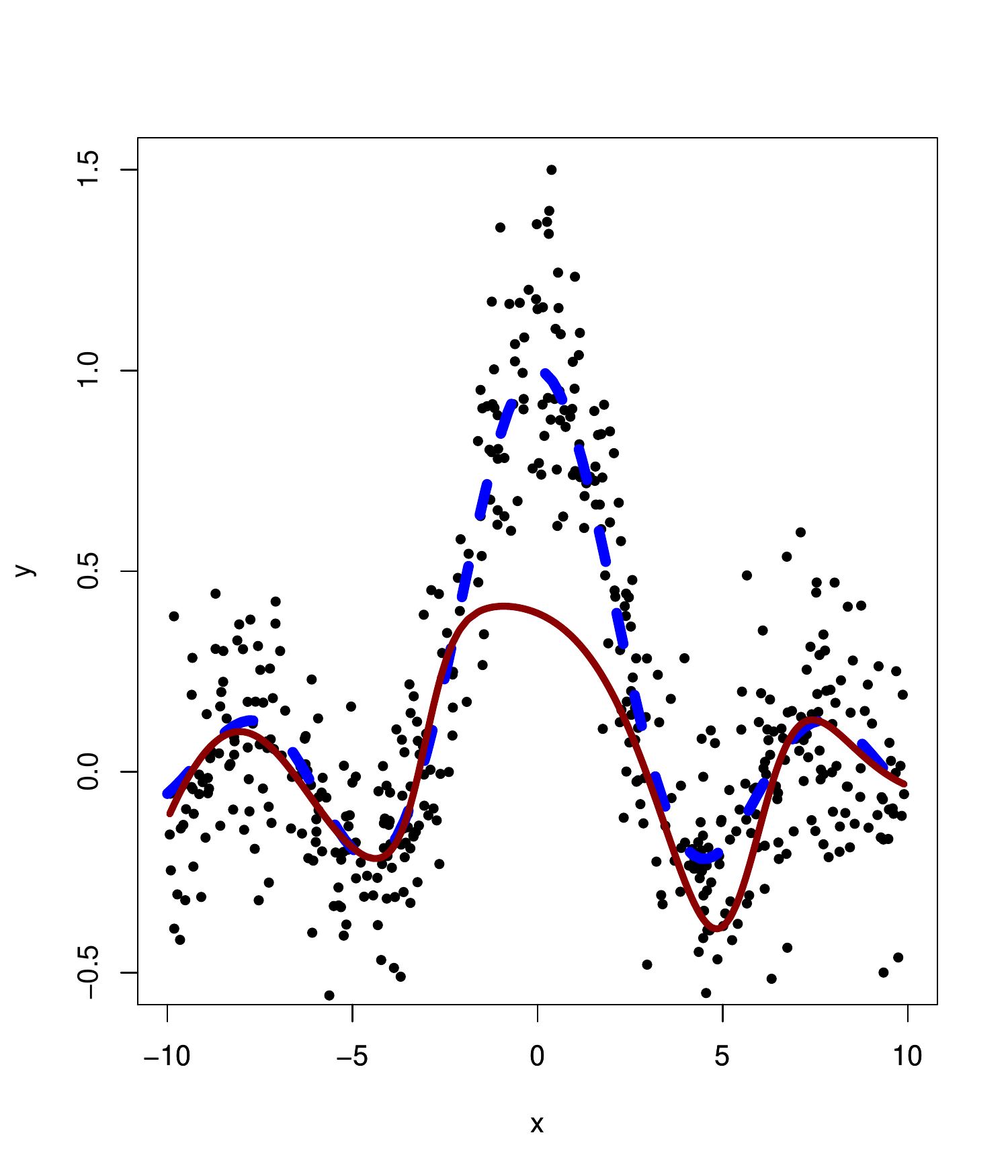}}
 \\
 \subfloat[noise level 0\%]{\label{fg:c3-d1_020_500_p2_f1_b0_m1_i200_n0_h5_cuda_nnet_test_predictions}
 \includegraphics[width=0.3\textwidth]{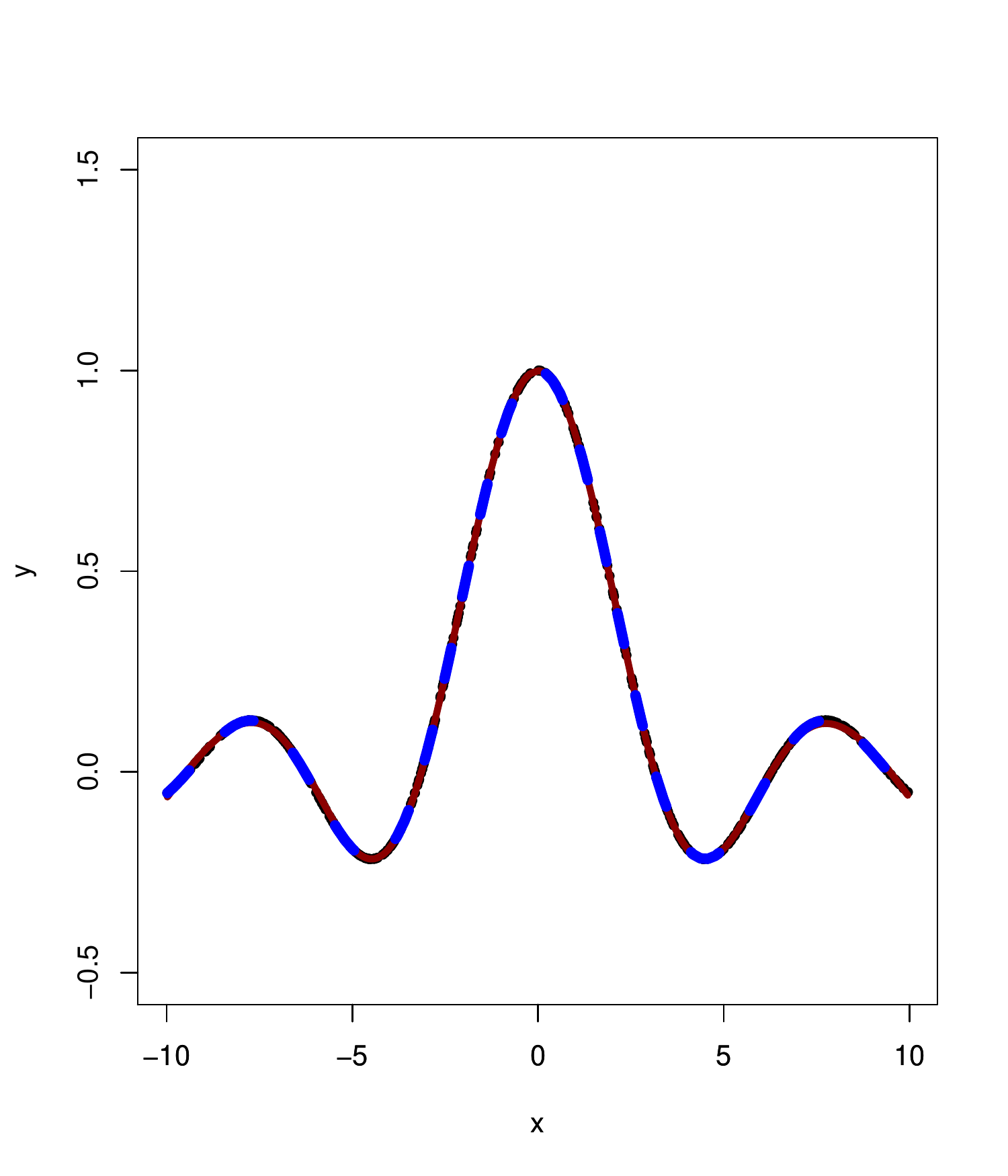}}
\subfloat[noise level 10\%]{\label{fg:c3-d1_020_500_p2_f1_b0_m1_i200_n10_h5_cuda_nnet_test_predictions}
 \includegraphics[width=0.3\textwidth]{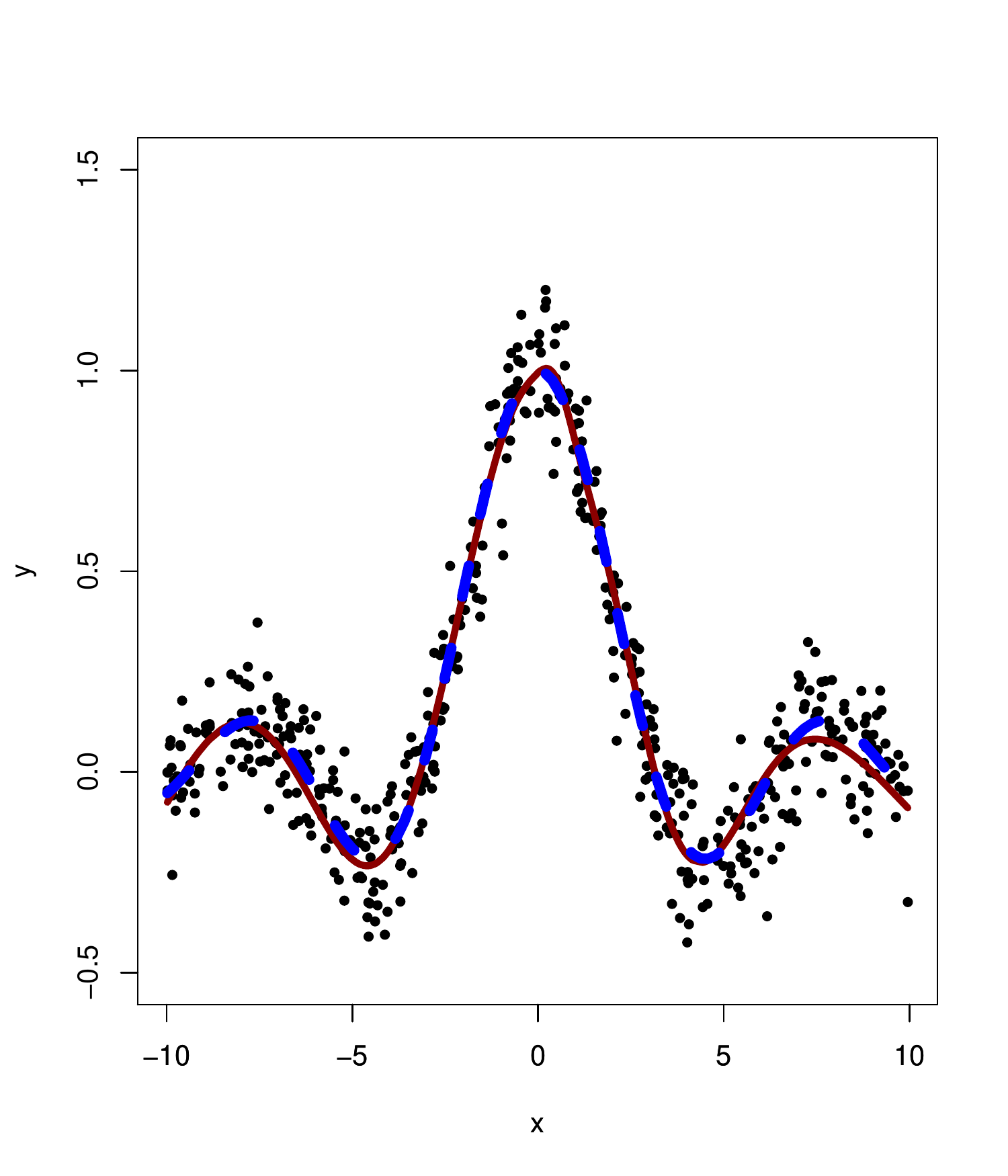}}
\subfloat[noise level 20\%]{\label{fg:c3-d1_020_500_p2_f1_b0_m1_i200_n20_h5_cuda_nnet_test_predictions}
 \includegraphics[width=0.3\textwidth]{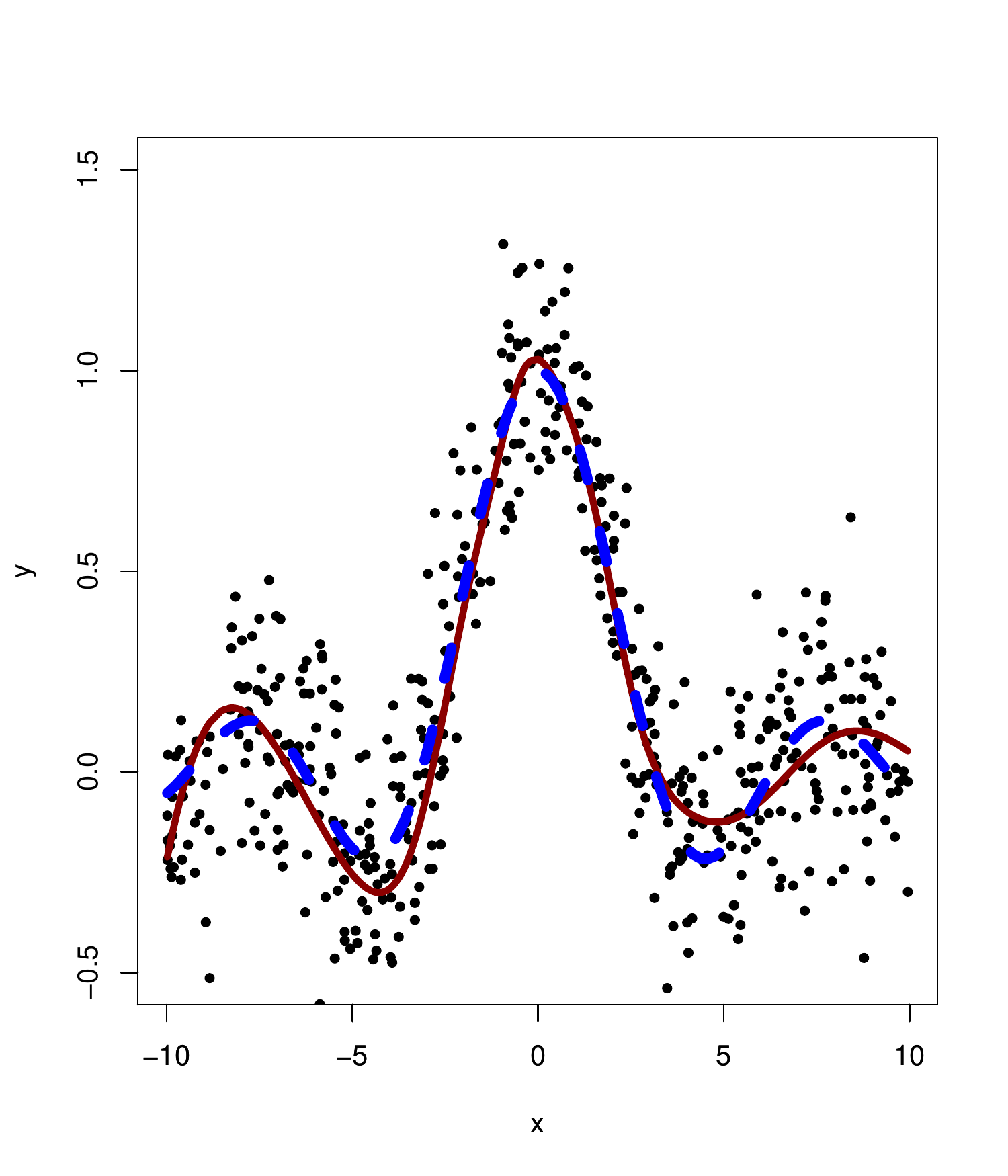}}
\caption{Data set 3, training with PCLTS with $C=1$, $B=0$ (LTS criterion, top) and  $C=8$, $B=8$ (bottom).
The proportion of outliers is $\delta=0.2$, and $n=500$.  Red solid curve is the ANN prediction, the blue dashed curve is the (noiseless) test data, and black dots are training data. The outliers are not shown. We observe that LTS treats some valid data as outliers (near the center of the graph).}
\label{fg:ds3train}

\end{figure}

\begin{figure}
\centering
\subfloat[noise level 0\%]{\label{fg:c10-d1_020_500_p2_f8_b8_m4_i200_n0_h10_cuda_nnet_test_predictions}
 \includegraphics[width=0.3\textwidth]{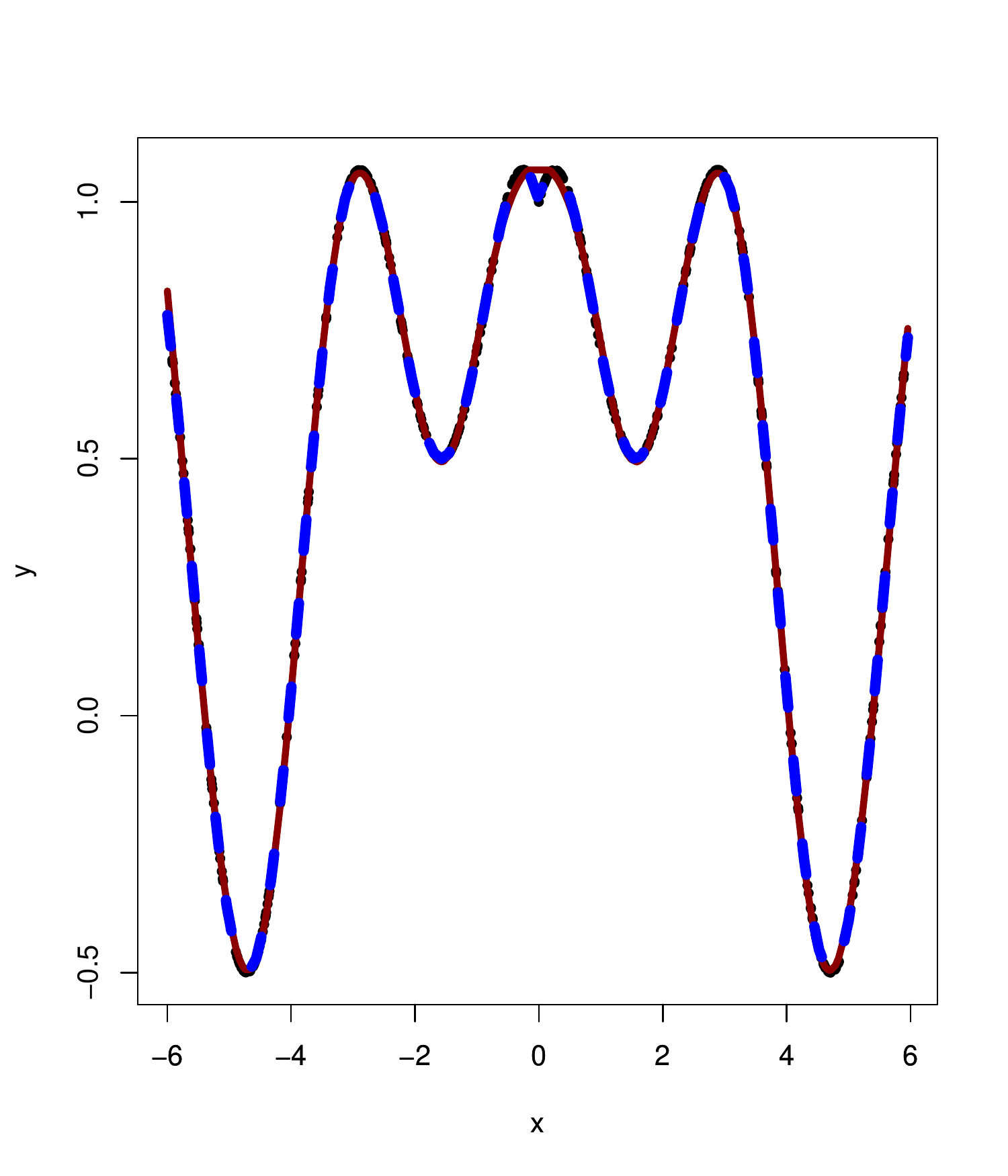}}
\subfloat[noise level 10\%]{\label{fg:c10-d1_020_500_p2_f8_b8_m4_i200_n10_h10_cuda_nnet_test_predictions}
 \includegraphics[width=0.3\textwidth]{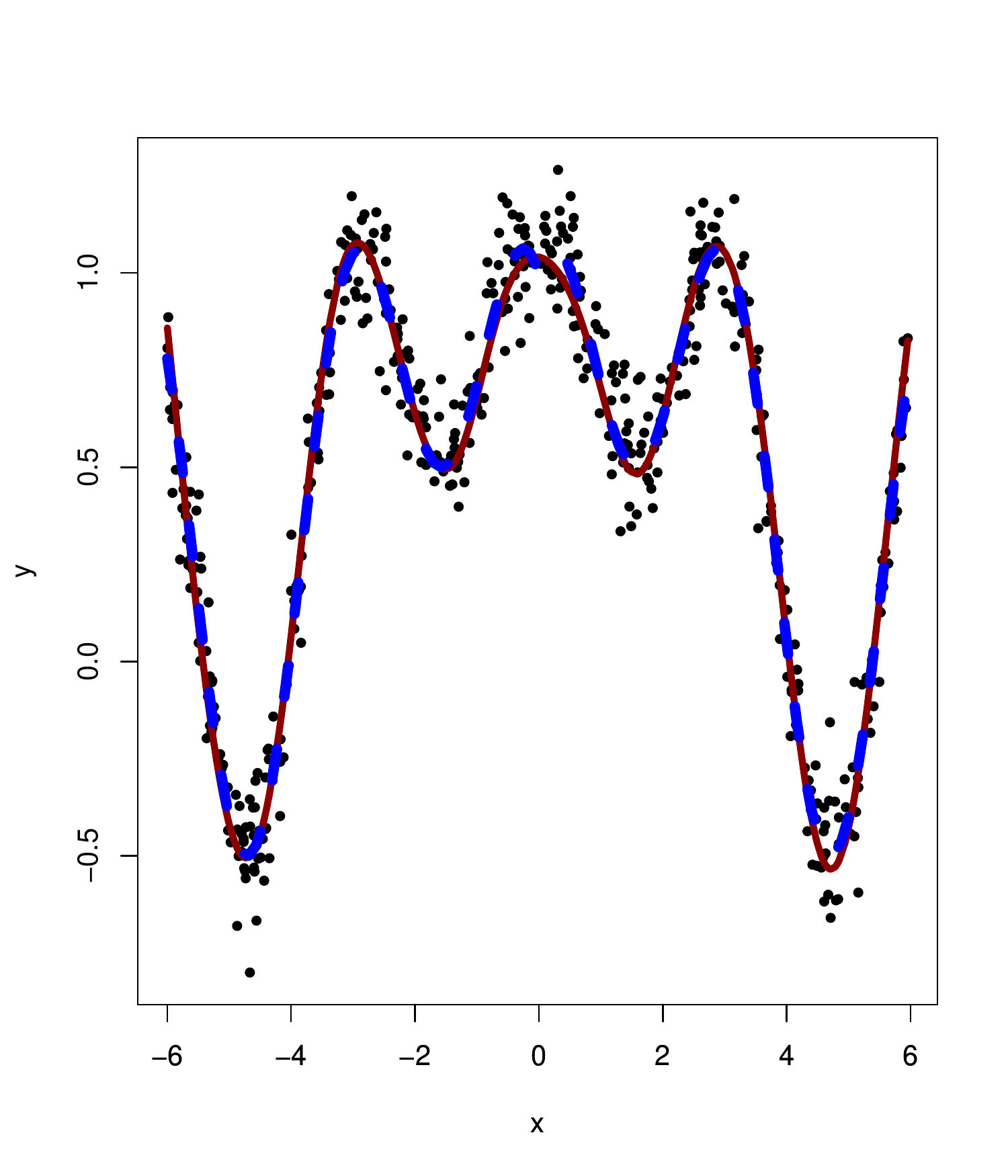}}
\subfloat[noise level 20\%]{\label{fg:c10-d1_020_500_p2_f8_b8_m4_i200_n20_h10_cuda_nnet_test_predictions}
 \includegraphics[width=0.3\textwidth]{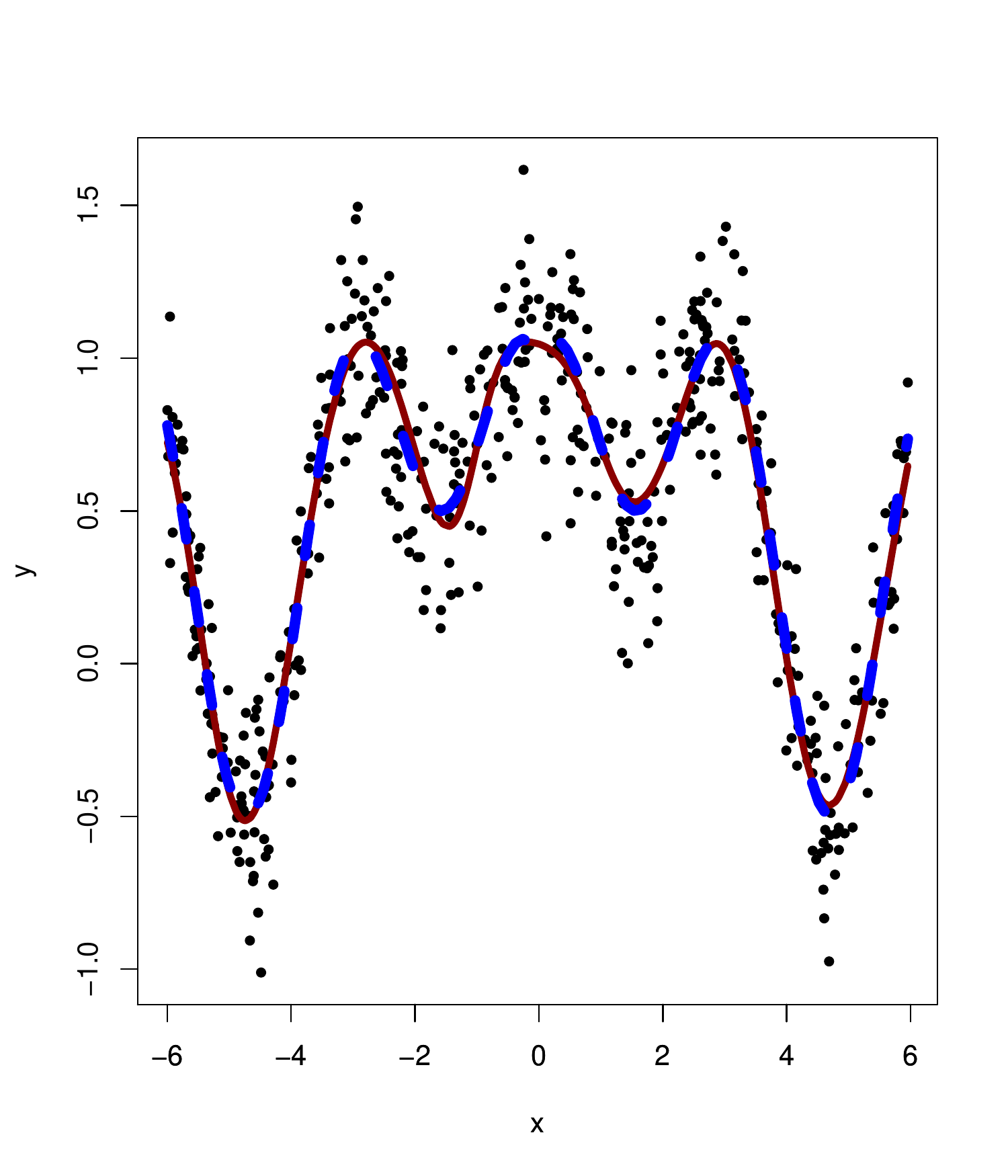}}
\\
\subfloat[noise level 0\%]{\label{fg:c10-d1_020_500_p2_f8_b8_m4_i200_n0_h5_cuda_nnet_test_predictions}
 \includegraphics[width=0.3\textwidth]{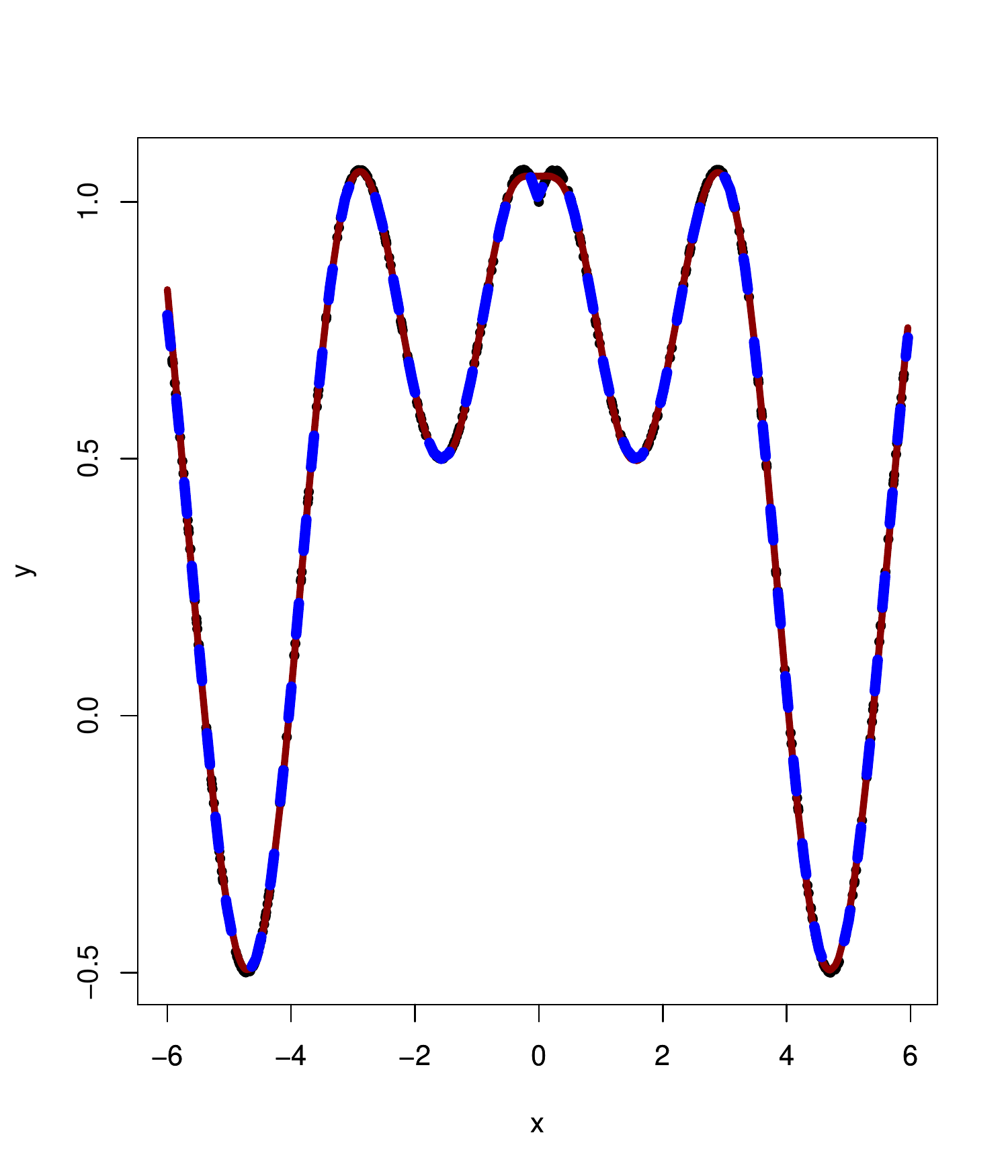}}
\subfloat[noise level 10\%]{\label{fg:c10-d1_020_500_p2_f8_b8_m4_i200_n10_h5_cuda_nnet_test_predictions}
 \includegraphics[width=0.3\textwidth]{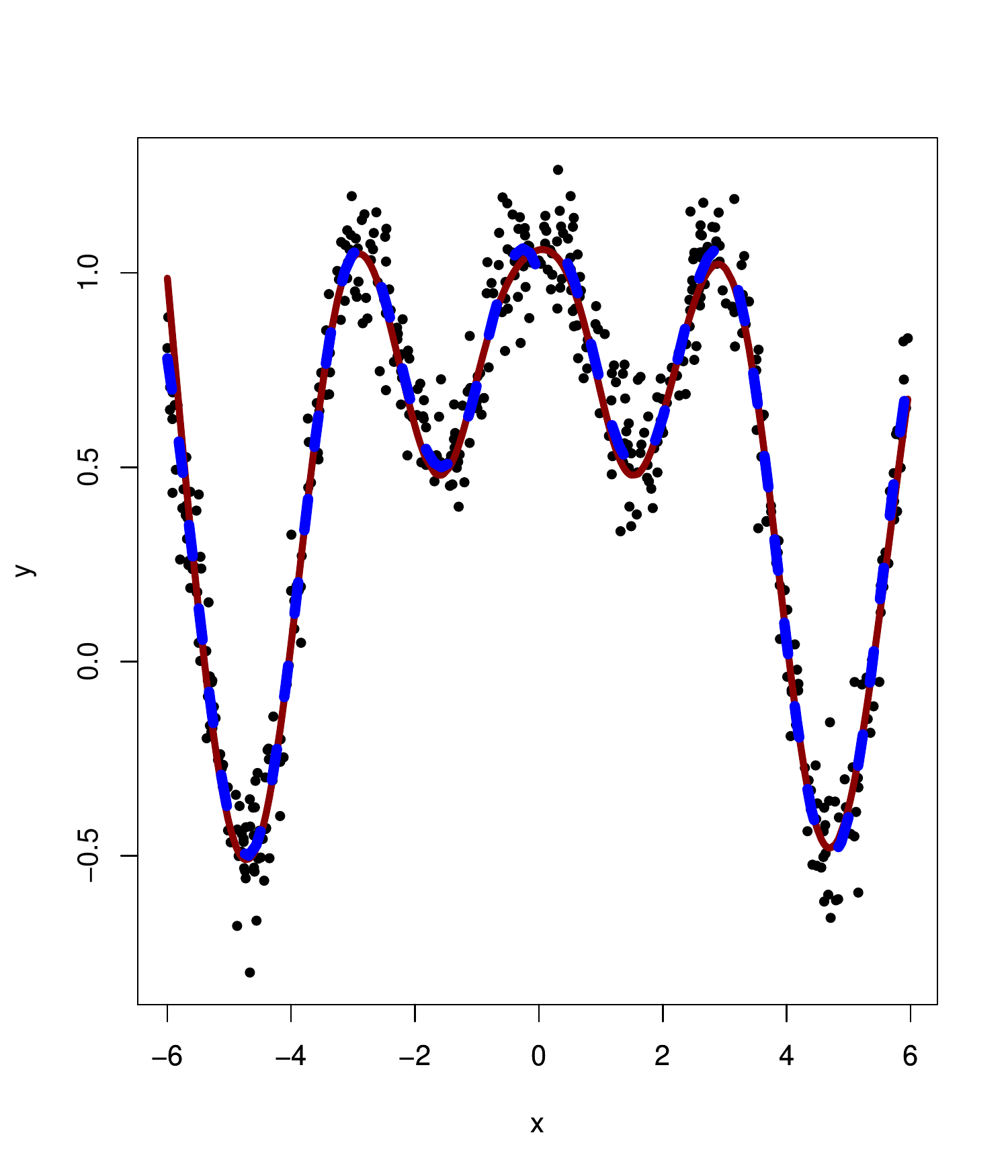}}
\subfloat[noise level 20\%]{\label{fg:c10-d1_020_500_p2_f8_b8_m4_i200_n20_h5_cuda_nnet_test_predictions}
 \includegraphics[width=0.3\textwidth]{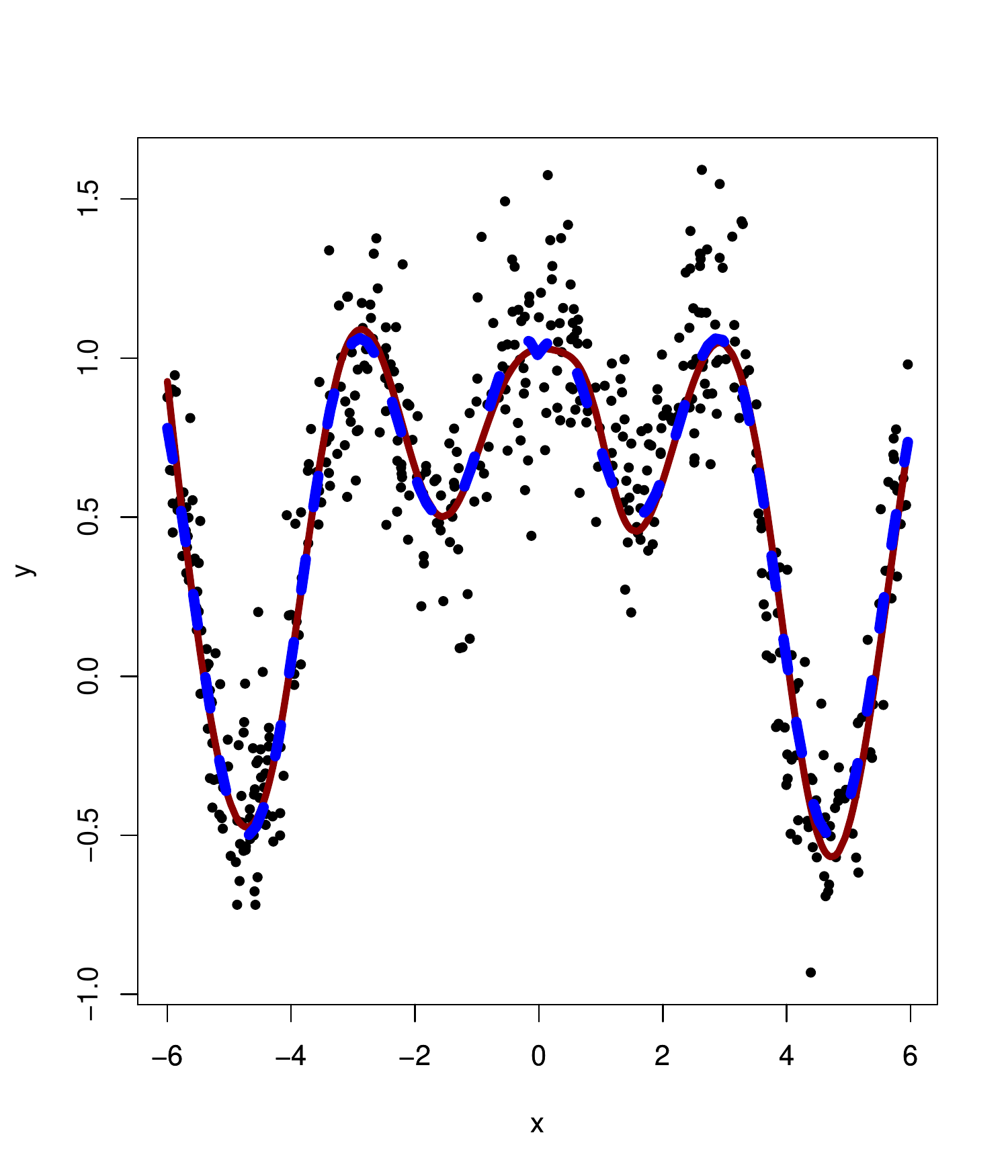}}
\caption{Data set 10, $\delta=0.2$, $n=500$. PCLTS ($C=8$, $B=8$) correctly predicts the model, using 10 (top row) and 5 (bottom row) hidden neurons.}
\label{fg:ds10f8b8h10}
\end{figure}

\begin{figure}[ph!]
\centering
\subfloat[noise level 0\%]{\label{fg:c10-d1_020_500_p2_f1_b0_m4_i200_n0_h10_cuda_nnet_test_predictions}
 \includegraphics[width=0.3\textwidth]{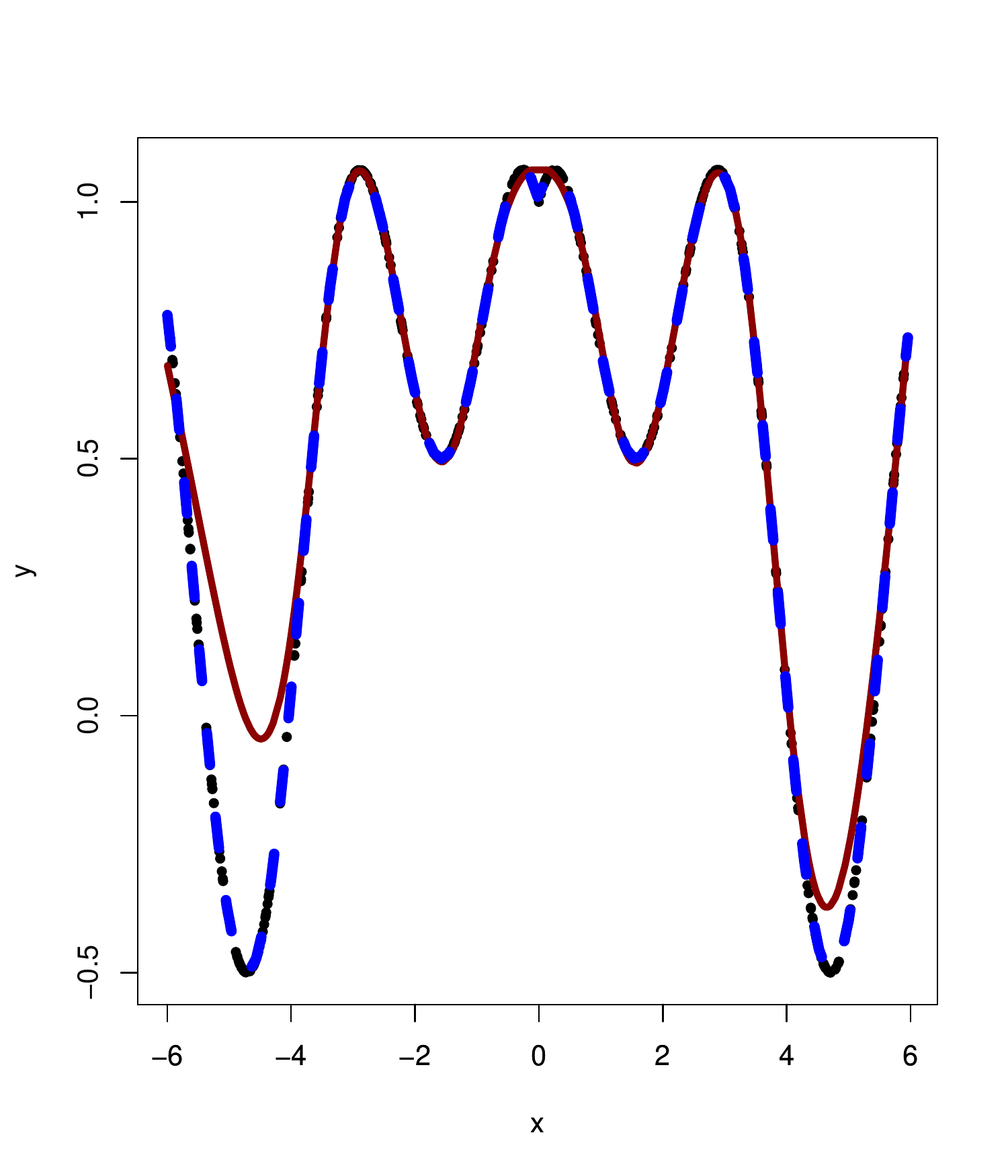}}
\subfloat[noise level 10\%]{\label{fg:c10-d1_020_500_p2_f1_b0_m4_i200_n10_h10_cuda_nnet_test_predictions}
 \includegraphics[width=0.3\textwidth]{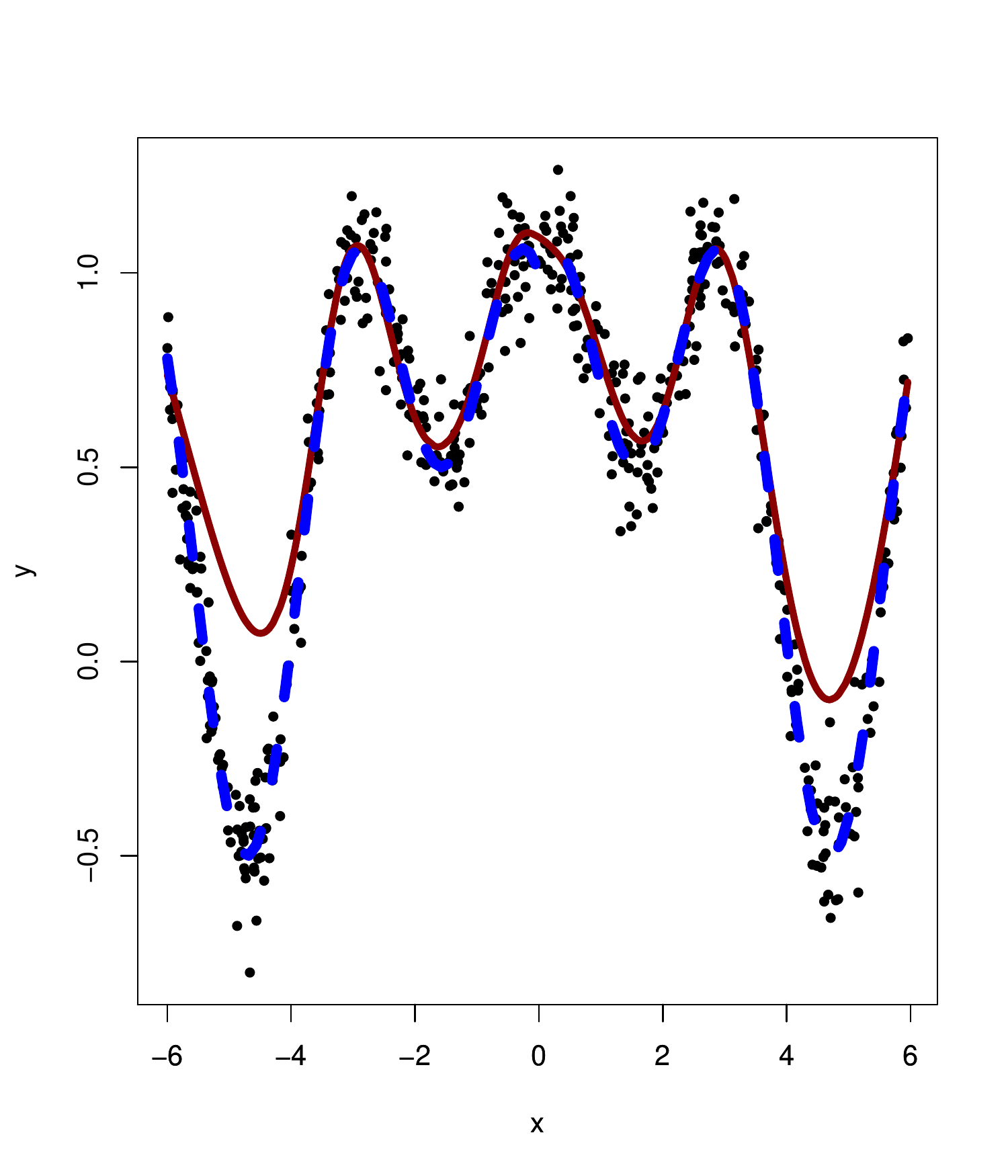}}
\subfloat[noise level 20\%]{\label{fg:c10-d1_020_500_p2_f1_b0_m4_i200_n20_h10_cuda_nnet_test_predictions}
 \includegraphics[width=0.3\textwidth]{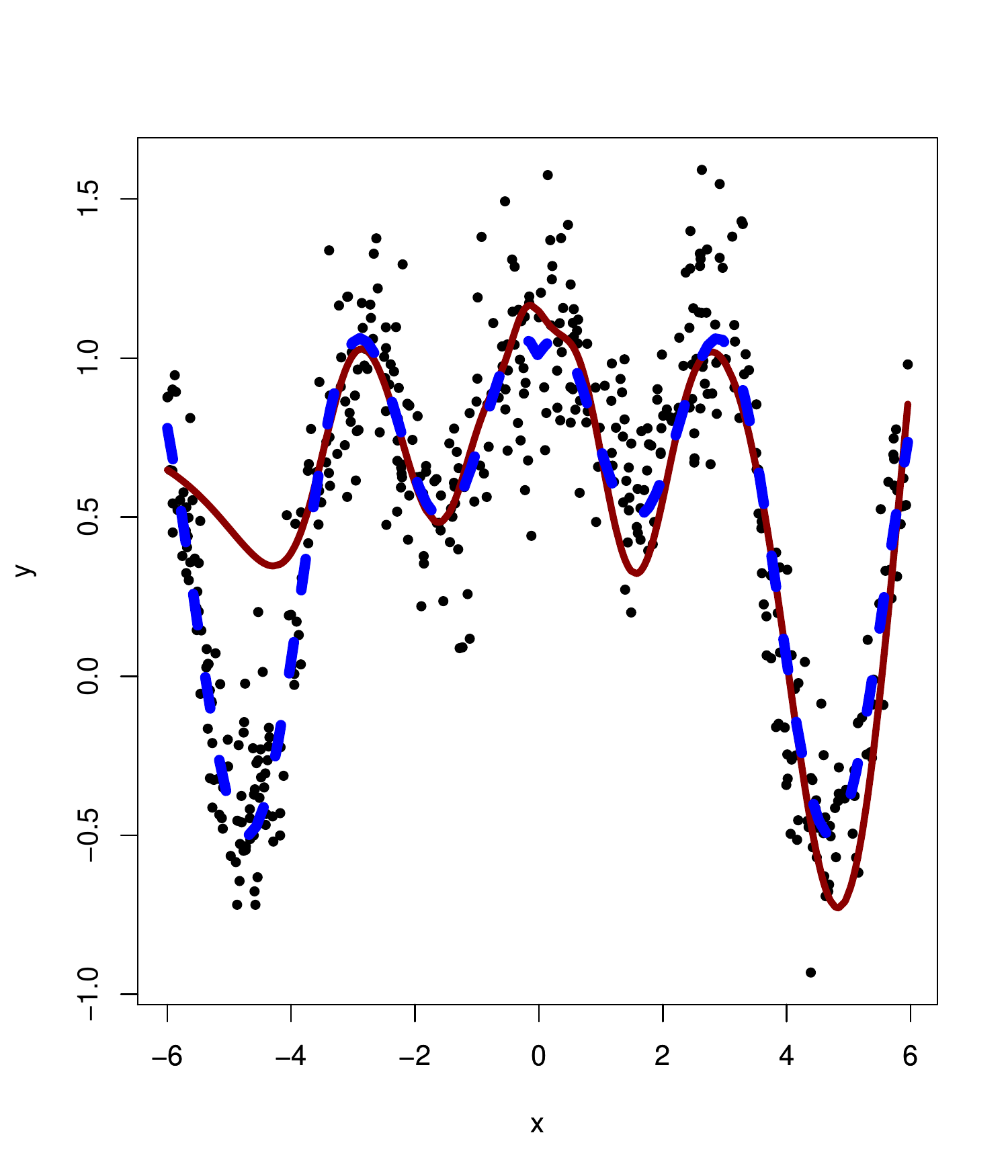}}
\\
\subfloat[noise level 0\%]{\label{fg:c10-d1_020_500_p2_f1_b0_m4_i200_n0_h5_cuda_nnet_test_predictions}
 \includegraphics[width=0.3\textwidth]{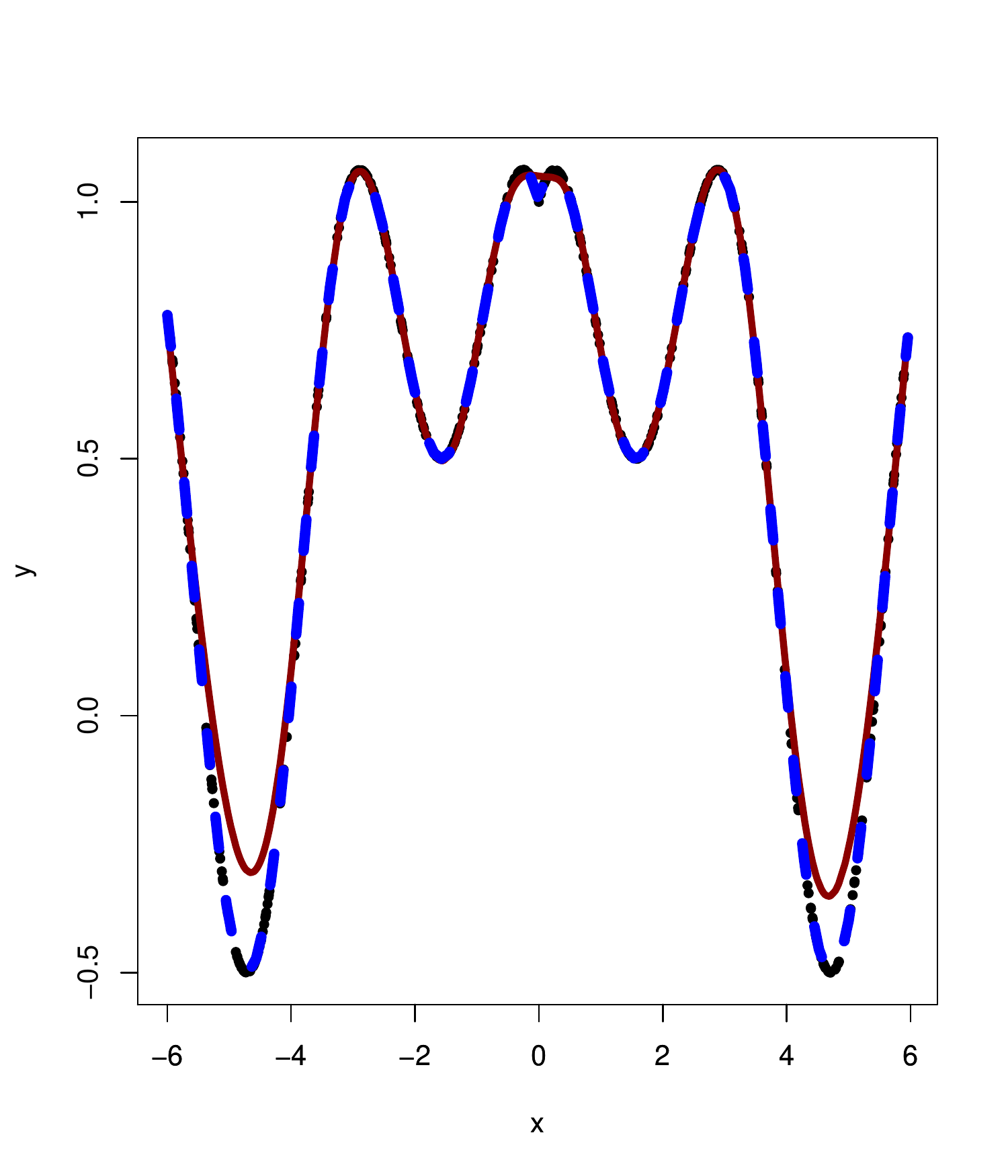}}
\subfloat[noise level 10\%]{\label{fg:c10-d1_020_500_p2_f1_b0_m4_i200_n10_h5_cuda_nnet_test_predictions}
 \includegraphics[width=0.3\textwidth]{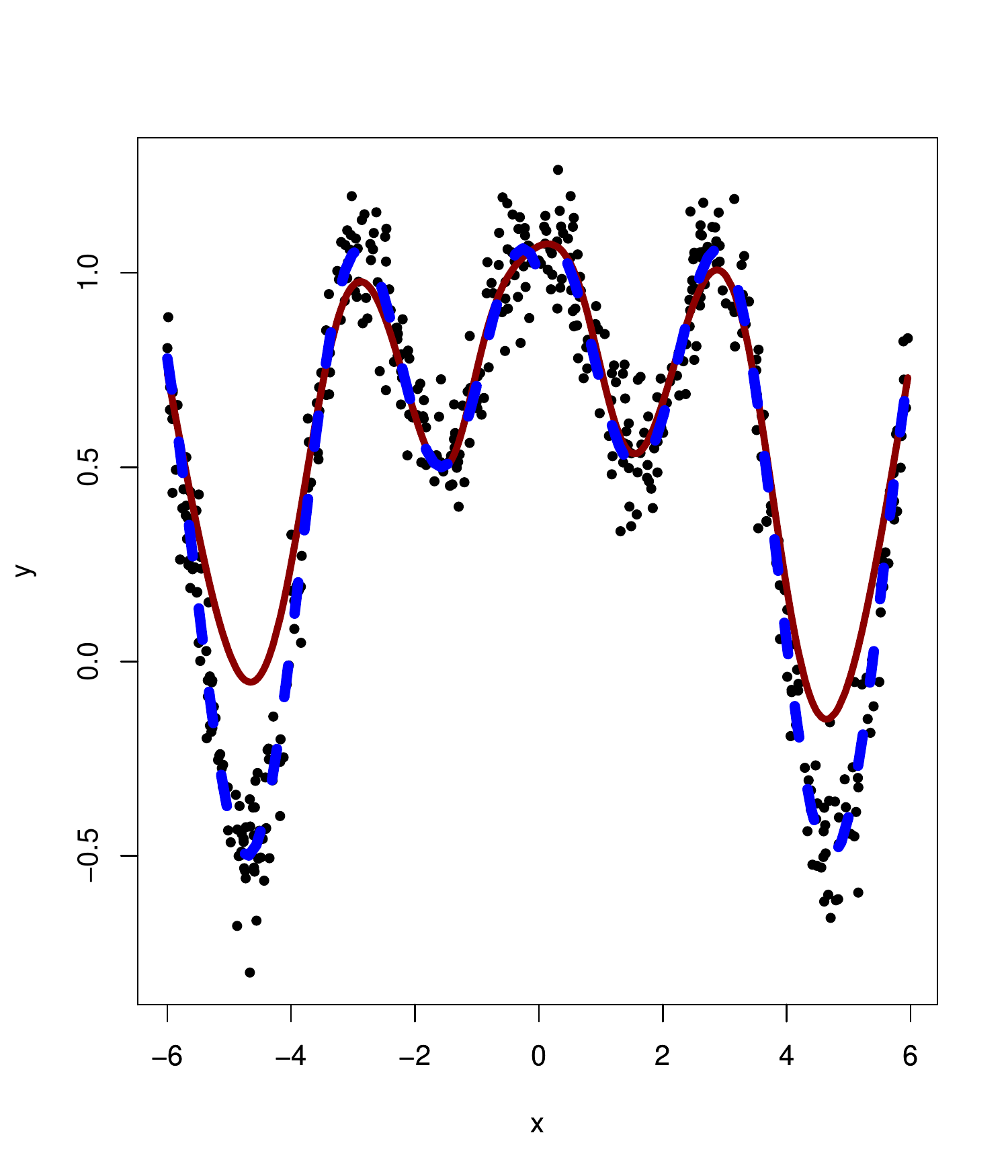}}
\subfloat[noise level 20\%]{\label{fg:c10-d1_020_500_p2_f1_b0_m4_i200_n20_h5_cuda_nnet_test_predictions}
 \includegraphics[width=0.3\textwidth]{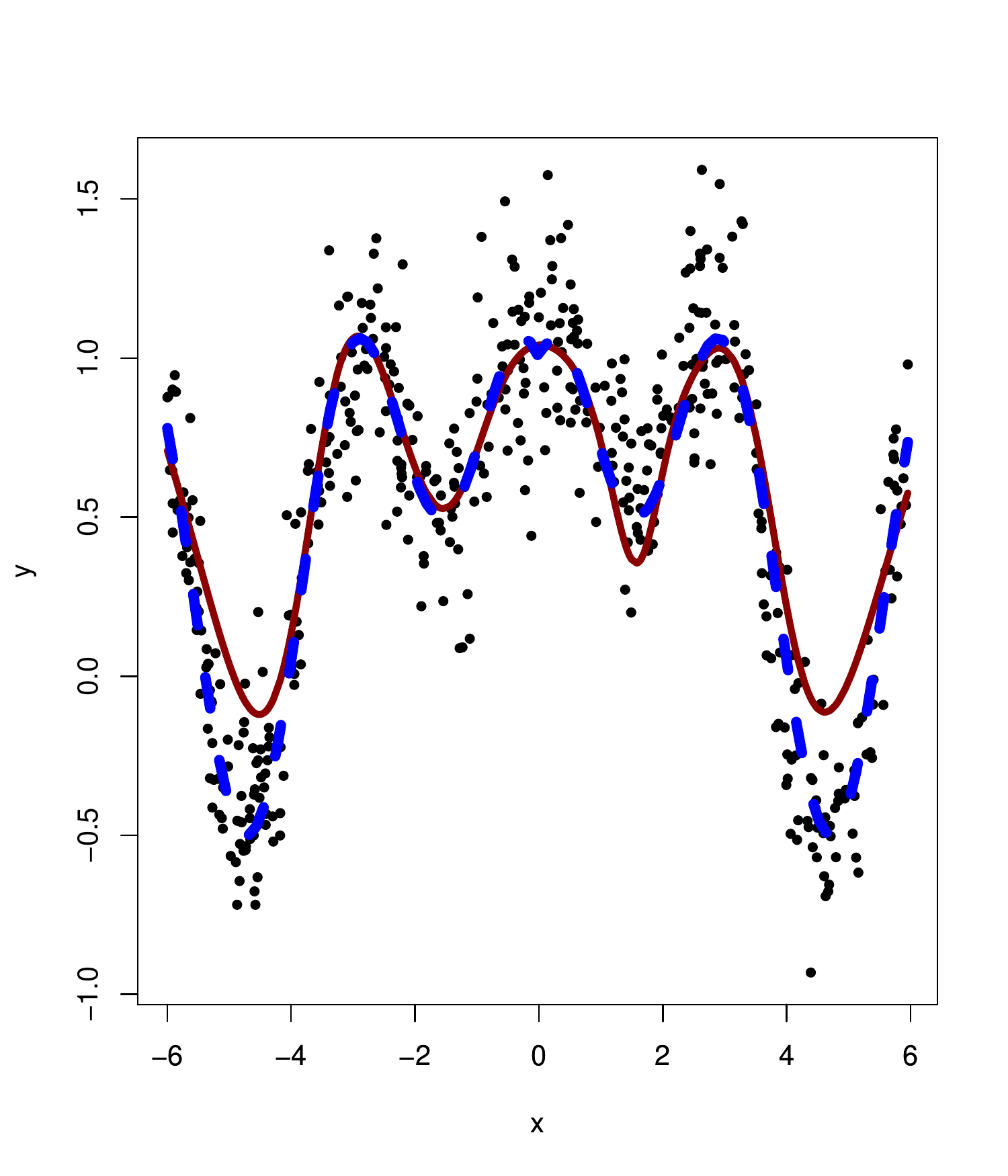}}
\caption{Data set 10, $\delta=0.2$, $n=500$. LTS ($C=1$, $B=0$, 10 hidden neurons (top), and 5 hidden neurons (bottom)) removed valid data near the two lowest minima of the model function, and hence did not achieve good accuracy. }
\label{fg:ds10f1b0h10}
\end{figure}

Let us now look at the figures. As we mentioned, Figures \ref{fg:ds1nnet} illustrate the inability of the standard ANNs to predict the correct model, and the ability of the robust training based on LTS and PCLTS to remove the outliers. Figures \ref{fg:ds1trainPCLTS} and \ref{fg:ds3train} (see also Figures  \ref{fg:ds10f8b8h10} and\ref{fg:ds10f1b0h10}) illustrate the difference between the LTS and PCLTS criteria. Note that using LTS (top row), the cusp of the graph of the model function near the origin is lost. The reason is that LTS criterion treated the data near the cusp as outliers. This is of course undesirable. The modified criterion PCLTS (bottom row)  gave a much better (in fact, nearly perfect) prediction. Hence PCLTS achieved its aim to block unnecessary removal of data, by imposing a penalty for every removal.
In Figures \ref{fg:ds10f8b8h10},\ref{fg:ds10f1b0h10} the same effect is visible at both ends of the domain of the data.
Thus we conclude that the addition of the penalty term in PCLTS criterion was justified.

As far as the optimisation methods are concerned, we noticed that the method DSO reliably detected the outliers, and its CPU was significantly smaller than that of the other two methods. CPU of random start with NEWUOA was on average six times higher, and CPU of random start with DFBM was twelve to fifteen times higher than that of DSO. On a few occasions random start with NEWUOA or DFBM have failed to deliver a model that identified the outliers correctly. Therefore we concluded that for our objective, DSO was the most suitable method among the alternatives we studied, and we
 concentrated on studying this method in greater detail. The tables in the Appendix give CPU and RMSE for the DSO method.

\newpage

\section{Conclusion} \label{sec_concl}

In this paper we investigated robust training of ANNs and detection of outliers in the data from two perspectives. First, we illustrated the inability of the classical least squares criterion to produce
regression models not sensitive to outliers in the data. We investigated the use of the least trimmed criterion and found that it delivered robust regression models. We benchmarked various derivative free optimisation algorithms for optimising the LTS criterion and found that the method based on dynamical systems optimisation (DSO) was superior to several alternative methods. A contributing factor here is relatively large number of variables, which ranged from ten to one hundred. We produced a hybrid algorithm, combining detection of outliers by optimising the LTS criterion, their removal and subsequent fine tuning of the ANN by backpropagation.

Second, we investigated an undesirable feature of the LTS criterion, unjustified removal of valid data, and subsequent loss of accuracy. We introduced a new criterion, called  Penalised CLTS, which imposes a penalty for removing the data. By optimising PCLTS criterion, we achieved a high degree of accuracy of the resulting regression model. Our method pinpoints and  filters out the outliers reliably, and its computational cost is reasonable.


\bibliographystyle{plain}
\bibliography{median}

\newpage
\section*{Appendix}


\renewcommand{\baselinestretch}{1}
\begin{table}[ph!]
\begin{center}
\caption{Data Set 1, RMSE scores
for test function $h(\vect x) = ||\vect x||^{2/3}$.
Values in bold indicate breakdown of the method.}
\label{t:rmse1}
{\tiny
\begin{tabular}{|r|r|r|r|r|r|r|r|r|r|}
\hline\hline
 Input  & Proportion  & Size of the& \multicolumn{2}{|c|}{Noise level =0} & \multicolumn{2}{|c|}{Noise level=0.1} & \multicolumn{2}{|c|}{Noise level=0.2} \\
dimension   &  of outliers  & data set $n$   & pclts   & nnet   & pclts   & nnet   & pclts   & nnet \\
\hline\hline
1  & 0  & 100  & 0.021  & 0.020  & 0.035  & 0.035  & 0.071  & 0.071 \\
1  & 0  & 500  & 0.011  & 0.011  & 0.016  & 0.016  & 0.042  & 0.042 \\
1  & 0  & 5000  & 0.006  & 0.006  & 0.008  & 0.008  & 0.010  & 0.009 \\
\hline
1  & 0.2  & 100  & 0.013  &{\bf 1109.050  }& 0.035  &{\bf 2151.400  }& 0.076  &{\bf 2221.849}\\
1  & 0.2  & 500  & 0.012  &{\bf 1626.452  }& 0.018  &{\bf 1447.013  }& 0.044  &{\bf 1018.946}\\
1  & 0.2  & 5000  & 0.007  &{\bf 1020.319  }& 0.419  &{\bf 1742.355  }& 0.014  &{\bf 1545.335}\\
\hline\hline
\end{tabular}
} 
\end{center}
\end{table}
\renewcommand{\baselinestretch}{2}

\renewcommand{\baselinestretch}{1}
\begin{table}[ph!]
\begin{center}
\caption{Data Set 2, RMSE scores for test function
$h(\vect x) = x_1 e^{||\vect x||}$.
Values in bold indicate breakdown of the method. }
\label{t:rmse2}
{\tiny
\begin{tabular}{|r|r|r|r|r|r|r|r|r|r|}
\hline\hline
 Input  & Proportion  & Size of the& \multicolumn{2}{|c|}{Noise level =0} & \multicolumn{2}{|c|}{Noise level=0.1} & \multicolumn{2}{|c|}{Noise level=0.2} \\
dimension   &  of outliers  & data set $n$   & pclts   & nnet   & pclts   & nnet   & pclts   & nnet \\
\hline\hline
1 &  0 &  100 &  0.013 &  0.012 &  0.048 &  0.034 &  0.061 &  0.061  \\
1 &  0 &  500 &  0.001 &  0.002 &  0.016 &  0.014 &  0.030 &  0.030  \\
1 &  0 &  5000 &  0.000 &  0.000 &  0.005 &  0.005 &  0.011 &  0.010  \\
\hline
1 &  0.2 &  100 &  0.010 & {\bf  2148.754 }&  0.033 &  2.417 &  0.080 & {\bf  2400.194  }\\
1 &  0.2 &  500 &  0.001 & {\bf  2602.706 }&  0.020 & {\bf  2604.085 }&  0.031 & {\bf  2276.562  }\\
1 &  0.2 &  5000 &  0.000 & {\bf  2334.587 }&  0.006 & {\bf  2887.469 }&  0.013 & {\bf  2337.995  }\\
\hline\hline
\end{tabular}
} 
\end{center}
\end{table}
\renewcommand{\baselinestretch}{2}

\renewcommand{\baselinestretch}{1}
\begin{table}[ph!]
\begin{center}
\caption{Data Set 3, RMSE scores for test function
$h(\vect x)=sin(||\vect x||)/||\vect x||$.
Values in bold indicate breakdown of the method. }
\label{t:rmse3}
{\tiny
\begin{tabular}{|r|r|r|r|r|r|r|r|r|r|}
\hline\hline
 Input  & Proportion  & Size of the& \multicolumn{2}{|c|}{Noise level =0} & \multicolumn{2}{|c|}{Noise level=0.1} & \multicolumn{2}{|c|}{Noise level=0.2} \\
dimension   &  of outliers  & data set $n$   & pclts   & nnet   & pclts   & nnet   & pclts   & nnet \\
\hline\hline
1 &  0 &  100 &  0.004 &  0.004 &  0.024 &  0.024 &  0.087 &  0.087  \\
1 &  0 &  500 &  0.002 &  0.002 &  0.019 &  0.021 &  0.041 &  0.042  \\
1 &  0 &  5000 &  0.002 &  0.002 &  0.007 &  0.007 &  0.011 &  0.011  \\
\hline
1 &  0.2 &  100 &  0.005 & {\bf  2350.122 }&  0.049 & {\bf  2621.706 }&  0.147 & {\bf  3816.212  }\\
1 &  0.2 &  500 &  0.003 & {\bf  2331.880 }&  0.019 & {\bf  2326.221 }&  0.044 & {\bf  1604.963  }\\
1 &  0.2 &  5000 &  0.002 & {\bf  2341.197 }&  0.007 & {\bf  1742.731 }&  0.012 & {\bf  1800.171  }\\
\hline\hline
\end{tabular}
} 
\end{center}
\end{table}
\renewcommand{\baselinestretch}{2}

\renewcommand{\baselinestretch}{1}
\begin{table}[ph!]
\begin{center}
\caption{Data Set 4, RMSE scores for test function
$h(\vect x) =
\sin(\frac{5}{m}\sum^m_{i=1}x_i)
\,\mbox{acos}(\frac{1}{m}\sum^m_{i=1}x_i)
\cos(\frac{3}{m}\sum^m_{i=1}x_i-2/n)$.
Values in bold indicate breakdown of the method.
 }
\label{t:rmse4}
{\tiny
\begin{tabular}{|r|r|r|r|r|r|r|r|r|r|}
\hline\hline
 Input  & Proportion  & Size of the& \multicolumn{2}{|c|}{Noise level =0} & \multicolumn{2}{|c|}{Noise level=0.1} & \multicolumn{2}{|c|}{Noise level=0.2} \\
dimension   &  of outliers  & data set $n$   & pclts   & nnet   & pclts   & nnet   & pclts   & nnet \\
\hline\hline
1 &  0 &  100 &  0.018 &  0.018 &  0.039 &  0.039 &  0.051 &  0.051  \\
1 &  0 &  500 &  0.009 &  0.010 &  0.021 &  0.020 &  0.040 &  0.043  \\
1 &  0 &  5000 &  0.004 &  0.004 &  0.005 &  0.005 &  0.007 &  0.007  \\
\hline
1 &  0.2 &  100 &  0.029 & {\bf  2536.255 }&  0.051 & {\bf  2281.526 }&  0.102 & {\bf  2113.818  }\\
1 &  0.2 &  500 &  0.011 & {\bf  2349.493 }&  0.022 & {\bf  2460.206 }&  0.023 & {\bf  2480.789  }\\
1 &  0.2 &  5000 &  0.005 & {\bf  2392.121 }&  0.008 & {\bf  2311.305 }&  0.015 & {\bf  2314.949  }\\
\hline\hline
\end{tabular}
} 
\end{center}
\end{table}
\renewcommand{\baselinestretch}{2}

\renewcommand{\baselinestretch}{1}
\begin{table}[ph!]
\begin{center}
\caption{Data Set 5, RMSE scores for test function
$h(\vect x)=\sin(10\pi
||\vect x||)+\sin(20\pi||\vect x||)$.
Values in bold indicate breakdown of the method.  }
\label{t:rmse5}
{\tiny
\begin{tabular}{|r|r|r|r|r|r|r|r|r|r|}
\hline\hline
 Input  & Proportion  & Size of the& \multicolumn{2}{|c|}{Noise level =0} & \multicolumn{2}{|c|}{Noise level=0.1} & \multicolumn{2}{|c|}{Noise level=0.2} \\
dimension   &  of outliers  & data set $n$   & pclts   & nnet   & pclts   & nnet   & pclts   & nnet \\
\hline\hline
1 &  0 &  100 &  0.114 &  0.109 &  0.075 &  0.056 &  0.148 &  0.146  \\
1 &  0 &  500 &  0.012 &  0.012 &  0.014 &  0.027 &  0.040 &  0.037  \\
1 &  0 &  5000 &  0.007 &  0.007 &  0.008 &  0.012 &  0.011 &  0.011  \\
\hline
1 &  0.2 &  100 &  0.129 & {\bf  2324.811 }&  0.133 & {\bf  2504.451 }&  0.178 & {\bf  1986.435  }\\
1 &  0.2 &  500 &  0.015 & {\bf  2376.646 }&  0.027 & {\bf  2306.402 }&  0.042 & {\bf  2510.584  }\\
1 &  0.2 &  5000 &  0.007 & {\bf  2139.114 }&  0.008 & {\bf  2391.442 }&  0.017 & {\bf  2139.109  }\\
\hline\hline
\end{tabular}
} 
\end{center}
\end{table}
\renewcommand{\baselinestretch}{2}

\renewcommand{\baselinestretch}{1}
\begin{table}[ph!]
\begin{center}
\caption{Data Set 6, RMSE scores for
$h(\vect x) = (x_1^2-x_2^2+x_3^2-x_4^2+\cdots)
\sin(0.5(x_1+x_3+\cdots))$.
Values in bold indicate breakdown of the method. }
\label{t:rmse6}
{\tiny
\begin{tabular}{|r|r|r|r|r|r|r|r|r|r|}
\hline\hline
 Input  & Proportion  & Size of the& \multicolumn{2}{|c|}{Noise level =0} & \multicolumn{2}{|c|}{Noise level=0.1} & \multicolumn{2}{|c|}{Noise level=0.2} \\
dimension   &  of outliers  & data set $n$   & pclts   & nnet   & pclts   & nnet   & pclts   & nnet \\
\hline\hline
1 &  0 &  100 &  0.011 &  0.011 &  0.032 &  0.031 &  0.116 &  0.129  \\
1 &  0 &  500 &  0.033 &  0.004 &  0.047 &  0.013 &  0.032 &  0.029  \\
1 &  0 &  5000 &  0.012 &  0.001 &  0.049 &  0.005 &  0.016 &  0.013  \\
\hline
1 &  0.2 &  100 &  0.086 & {\bf  1428.515 }&  0.205 & {\bf  1244.780 }&  0.139 & {\bf  2210.010  }\\
1 &  0.2 &  500 &  0.004 & {\bf  2050.179 }&  0.020 & {\bf  2343.837 }&  0.041 & {\bf  2274.445  }\\
1 &  0.2 &  5000 &  0.001 & {\bf  2131.853 }&  0.007 & {\bf  1964.621 }&  0.011 & {\bf  1938.302  }\\
\hline\hline
\end{tabular}
} 
\end{center}
\end{table}
\renewcommand{\baselinestretch}{2}

\renewcommand{\baselinestretch}{1}
\begin{table}[ph!]
\begin{center}
\caption{Data Set 7, RMSE scores for test function
$h(\vect x) = \frac{\sin(x_1+x_3+\cdots)}{x_1+x_3+\cdots}
\frac{\sin(x_2+x_4+\cdots)}{x_2+x_4+\cdots}$.
Values in bold indicate breakdown of the method. }
\label{t:rmse7}
{\tiny
\begin{tabular}{|r|r|r|r|r|r|r|r|r|r|}
\hline\hline
 Input  & Proportion  & Size of the& \multicolumn{2}{|c|}{Noise level =0} & \multicolumn{2}{|c|}{Noise level=0.1} & \multicolumn{2}{|c|}{Noise level=0.2} \\
dimension   &  of outliers  & data set $n$   & pclts   & nnet   & pclts   & nnet   & pclts   & nnet \\
\hline
\hline
1 &  0 &  100 &  0.009 &  0.009 &  0.054 &  0.054 &  0.106 &  0.113  \\
1 &  0 &  500 &  0.002 &  0.002 &  0.016 &  0.016 &  0.049 &  0.047  \\
1 &  0 &  5000 &  0.000 &  0.000 &  0.005 &  0.005 &  0.013 &  0.013  \\
\hline
1 &  0.2 &  100 &  0.008 & {\bf  1765.182 }&  0.050 & {\bf  1758.972 }&  0.055 & {\bf  2102.820  }\\
1 &  0.2 &  500 &  0.002 & {\bf  2188.202 }&  0.016 & {\bf  2113.608 }&  0.042 & {\bf  2097.431  }\\
1 &  0.2 &  5000 &  0.001 & {\bf  1826.183 }&  0.004 & {\bf  1980.115 }&  0.010 & {\bf  2050.203  }\\
\hline\hline
\end{tabular}
} 
\end{center}
\end{table}
\renewcommand{\baselinestretch}{2}

\renewcommand{\baselinestretch}{1}
\begin{table}[ph!]
\begin{center}
\caption{Data Set 8, RMSE scores for test function
$h(\vect x) = 0.2(x_1x_2\cdots x_m)+
1.2\sin(||\vect x||^2)$.
Values in bold indicate breakdown of the method. }
\label{t:rmse8}
{\tiny
\begin{tabular}{|r|r|r|r|r|r|r|r|r|r|}
\hline\hline
 Input  & Proportion  & Size of the& \multicolumn{2}{|c|}{Noise level =0} & \multicolumn{2}{|c|}{Noise level=0.1} & \multicolumn{2}{|c|}{Noise level=0.2} \\
dimension   &  of outliers  & data set $n$   & pclts   & nnet   & pclts   & nnet   & pclts   & nnet \\
\hline
\hline
1 &  0 &  100 &  0.041 &  0.012 &  0.095 &  0.042 &  0.057 &  0.048  \\
1 &  0 &  500 &  0.004 &  0.004 &  0.013 &  0.013 &  0.038 &  0.038  \\
1 &  0 &  5000 &  0.001 &  0.001 &  0.007 &  0.004 &  0.012 &  0.012  \\
\hline
1 &  0.2 &  100 &  0.014 & {\bf  3128.316 }&  0.044 & {\bf  3952.802 }&  0.076 & {\bf  3533.059  }\\
1 &  0.2 &  500 &  0.009 & {\bf  2259.867 }&  0.031 & {\bf  2258.475 }&  0.070 & {\bf  2098.067  }\\
1 &  0.2 &  5000 &  0.001 & {\bf  2187.879 }&  0.006 & {\bf  2303.391 }&  0.017 & {\bf  2237.420  }\\
\hline\hline
\end{tabular}
} 
\end{center}
\end{table}
\renewcommand{\baselinestretch}{2}

\renewcommand{\baselinestretch}{1}
\begin{table}[ph!]
\begin{center}
\caption{Data Set 9, RMSE scores for test function
$h(\vect x) = \max\{ \exp(-10 x_1^2),
\exp(-50x_2^2), 1.25\exp(-5(||\vect x||^2))\}$.
Values in bold indicate breakdown of the method. }
\label{t:rmse9}
{\tiny
\begin{tabular}{|r|r|r|r|r|r|r|r|r|r|}
\hline\hline
 Input  & Proportion  & Size of the& \multicolumn{2}{|c|}{Noise level =0} & \multicolumn{2}{|c|}{Noise level=0.1} & \multicolumn{2}{|c|}{Noise level=0.2} \\
dimension   &  of outliers  & data set $n$   & pclts   & nnet   & pclts   & nnet   & pclts   & nnet \\
\hline
\hline
2 &  0 &  100 &  0.184 &  0.141 &  0.206 &  0.153 &  0.308 &  0.309  \\
2 &  0 &  500 &  0.079 &  0.075 &  0.099 &  0.088 &  0.106 &  0.102  \\
2 &  0 &  5000 &  0.059 &  0.059 &  0.064 &  0.065 &  0.067 &  0.066  \\
\hline
2 &  0.2 &  100 &  0.184 & {\bf  5173.392 }&  0.250 & {\bf  2889.054 }&  0.306 & {\bf  3264.975  }\\
2 &  0.2 &  500 &  0.062 & {\bf  2448.795 }&  0.078 & {\bf  2744.082 }&  0.111 & {\bf  2458.956  }\\
2 &  0.2 &  5000 &  0.059 & {\bf  2224.556 }&  0.056 & {\bf  2061.706 }&  0.059 & {\bf  2159.524  }\\
\hline\hline
\end{tabular}
} 
\end{center}
\end{table}
\renewcommand{\baselinestretch}{2}

\renewcommand{\baselinestretch}{1}
\begin{table}[ph!]
\begin{center}
\caption{Data Set 10, RMSE scores for test function
$h(\vect x) = 0.5 ||\vect x|| \sin(||\vect x||) + \cos^2(||\vect x||)$.
Values in bold indicate breakdown of the method. }
\label{t:rmse10}
{\tiny
\begin{tabular}{|r|r|r|r|r|r|r|r|r|r|}
\hline\hline
 Input  & Proportion  & Size of the& \multicolumn{2}{|c|}{Noise level =0} & \multicolumn{2}{|c|}{Noise level=0.1} & \multicolumn{2}{|c|}{Noise level=0.2} \\
dimension   &  of outliers  & data set $n$   & pclts   & nnet   & pclts   & nnet   & pclts   & nnet \\
\hline
\hline
1 &  0 &  100 &  0.012 &  0.012 &  0.038 &  0.037 &  0.115 &  0.116  \\
1 &  0 &  500 &  0.011 &  0.011 &  0.020 &  0.020 &  0.042 &  0.041  \\
1 &  0 &  5000 &  0.008 &  0.008 &  0.015 &  0.011 &  0.015 &  0.014  \\
\hline
1 &  0.2 &  100 &  0.014 & {\bf  3057.407 }&  0.043 & {\bf  3137.274 }&  0.116 & {\bf  3190.071  }\\
1 &  0.2 &  500 &  0.010 & {\bf  1692.090 }&  0.028 & {\bf  1778.629 }&  0.035 & {\bf  1845.963  }\\
1 &  0.2 &  5000 &  0.008 & {\bf  2209.299 }&  0.009 & {\bf  1780.260 }&  0.015 & {\bf  2128.146  }\\
\hline\hline
\end{tabular}
} 
\end{center}
\end{table}
\renewcommand{\baselinestretch}{2}

\renewcommand{\baselinestretch}{1}
\begin{table}[ph!]
\begin{center}
\caption{Data Set 11, RMSE scores for
prediction of the electricity consumption.
Values in bold indicate breakdown of the method.}
\label{t:rmse11}
{\small
\begin{tabular}{|r|r|r|r|r|r|r|}
\hline\hline
Proportion of outliers & PCLTS train & PCLTS test & nnet train & nnet test & CPU PCLTS & CPU nnet \\
\hline\hline
0  &  0.07377254  &  0.07669042  &  0.07662421  &  0.07942267  &  732.852  &  4.9196  \\
0.2  &  0.0732535  &  0.0768083  & {\bf  1721.716  }& {\bf  1071.685  }&  1313.783  &  5.0199  \\
0.4  &  0.07328527  &  0.07779184  & {\bf  3395.867  }& {\bf  2798.297  }&  800.7642  &  1.9977  \\
0.5  &  0.07314826  &  0.07771504  & {\bf  3726.083  }& {\bf  3783.546  }&  1585.115  &  1.3546  \\
\hline\hline
\end{tabular}
} 
\end{center}
\end{table}
\renewcommand{\baselinestretch}{2}

\renewcommand{\baselinestretch}{1}
\begin{table}[ph!]
\begin{center}
\caption{Data Sets 12, RMSE scores for
prediction of the cold water consumption.
Values in bold indicate breakdown of the method. }
\label{t:rmse12}
{\small
\begin{tabular}{|r|r|r|r|r|r|r|}
\hline\hline
Proportion of outliers & PCLTS train & PCLTS test & nnet train & nnet test & CPU PCLTS & CPU nnet \\
\hline\hline
0  &  0.03352069  &  0.03799872  &  0.04224156  &  0.04277332  &  1551.998  &  5.5072  \\
0.2  &  0.03314315  &  0.03812518  & {\bf  1504.155  }& {\bf  736.1127  }&  879.4155  &  3.6342  \\
0.4  &  0.03574897  &  0.03810132  & {\bf  3112.657  }& {\bf  2725.164  }&  676.9802  &  2.4515  \\
0.5  &  0.03614024  &  0.03821851  & {\bf  3522.750  }& {\bf  3630.974  }&  960.6753  &  1.6193  \\
\hline\hline
\end{tabular}
} 
\end{center}
\end{table}
\renewcommand{\baselinestretch}{2}

\renewcommand{\baselinestretch}{1}
\begin{table}[ph!]
\begin{center}
\caption{Data Sets 13, RMSE scores for
prediction of the hot water consumption.
Values in bold indicate breakdown of the method.}
\label{t:rmse13}
{\small
\begin{tabular}{|r|r|r|r|r|r|r|}
\hline\hline
Proportion of outliers & PCLTS train & PCLTS test & nnet train & nnet test & CPU PCLTS & CPU nnet \\
\hline\hline
0  &  0.05732137  &  0.05889331  &  0.05732247  &  0.05874144  &  582.1506  &  5.6062  \\
0.2  &  0.05633721  &  0.05898911  & {\bf  1966.307  }& {\bf  1136.745  }&  1025.413  &  4.9611  \\
0.4  &  0.055599  &  0.05912662  & {\bf  3335.348  }& {\bf  2747.096  }&  604.2534  &  3.0936  \\
0.5  &  0.05553557  &  0.05952139  & {\bf  1685.271  }& {\bf  1570.897  }&  1657.796  &  6.5511  \\
\hline\hline
\end{tabular}
} 
\end{center}
\end{table}
\renewcommand{\baselinestretch}{2}

\renewcommand{\baselinestretch}{1}
\begin{table}[ph!]
\begin{center}
\caption{Data Sets 14, RMSE scores for
prediction of heart disease.
Values in bold indicate breakdown of the method. }
\label{t:rmse14}
{\small
\begin{tabular}{|r|r|r|r|r|r|r|}
\hline\hline
Proportion of outliers & PCLTS train & PCLTS test & nnet train & nnet test & CPU PCLTS & CPU nnet \\
\hline\hline
0  &  0.2139786  &  0.2517887  &  0.2118541  &  0.2665598  &  922.7474  &  0.9273  \\
0.2  &  0.2190235  &  0.2697414  & {\bf  2651.938  }& {\bf  1682.814  }&  1284.666  &  0.4501  \\
0.4  &  0.2071636  &  0.292656  & {\bf  2926.835  }& {\bf  2752.914  }&  636.9798  &  0.371  \\
0.5  &  0.2002703  &  0.3038773  & {\bf  2934.247  }& {\bf  3429.113  }&  448.0119  &  0.7827  \\
\hline\hline
\end{tabular}
} 
\end{center}
\end{table}
\renewcommand{\baselinestretch}{2}


\renewcommand{\baselinestretch}{1}
\setlength{\tabcolsep}{1pt}
\noindent
$\:$\hskip -2cm
\begin{table}[ph!]
\caption{CPU times for PCLTS and standard ANN training using nnet, in seconds.}
\label{t:all-cpu-large1}
{\tiny
\begin{tabular}{||r|r|r||r|r||r|r||r|r||r|r||r|r||r|r||r|r||r|r||r|r||r|r||}
\hline\hline
    & Out- &   &\multicolumn{2}{|c||}{Data set 1}&\multicolumn{2}{|c||}{Data set 2}&\multicolumn{2}{|c||}{Data set 3}&\multicolumn{2}{|c||}{Data set 4}&\multicolumn{2}{|c||}{Data set 5}&\multicolumn{2}{|c||}{Data set 6}&\multicolumn{2}{|c||}{Data set 7}&\multicolumn{2}{|c||}{Data set 8}&\multicolumn{2}{|c||}{Data set 9}&\multicolumn{2}{|c||}{Data set 10}\\
Dim & liers & Size & PCLTS & nnet & PCLTS & nnet & PCLTS & nnet & PCLTS & nnet & PCLTS & nnet & PCLTS & nnet & PCLTS & nnet & PCLTS & nnet & PCLTS & nnet & PCLTS & nnet\\
\hline\hline
1 & 0 & 100 & 496.75 & 0.04 & 219.32 & 0.04 & 252.36 & 0.07 & 250.17 & 0.07 & 135.70 & 0.10 & 212.33 & 0.05 & 392.40 & 0.07 & 1455.44 & 0.09 &  &  & 304.92 & 0.08\\
1 & 0 & 500 & 679.66 & 0.21 & 242.90 & 0.30 & 308.30 & 0.43 & 192.44 & 0.38 & 173.20 & 0.55 & 154.51 & 0.23 & 246.51 & 0.38 & 57.93 & 0.40 &  &  & 231.32 & 0.41\\
1 & 0 & 5000 & 898.40 & 4.26 & 702.94 & 2.66 & 511.65 & 3.42 & 242.22 & 4.29 & 535.68 & 7.37 & 96.90 & 2.85 & 578.60 & 3.19 & 275.26 & 4.70 &  &  & 213.27 & 5.55\\
1 & 0.2 & 100 & 2034.51 & 0.04 & 224.57 & 0.03 & 352.25 & 0.07 & 132.74 & 0.04 & 129.44 & 0.07 & 155.90 & 0.04 & 2035.75 & 0.04 & 365.34 & 0.07 &  &  & 2567.31 & 0.07\\
1 & 0.2 & 500 & 321.40 & 0.18 & 137.18 & 0.18 & 301.95 & 0.26 & 199.56 & 0.24 & 94.57 & 0.40 & 298.49 & 0.28 & 832.34 & 0.26 & 431.11 & 0.30 &  &  & 1230.57 & 0.37\\
1 & 0.2 & 5000 & 84.55 & 2.64 & 276.78 & 2.55 & 41.44 & 2.72 & 238.21 & 4.04 & 40.72 & 5.00 & 295.15 & 3.21 & 204.76 & 3.05 & 1130.11 & 4.28 &  &  & 435.72 & 5.36\\
1 & 0.4 & 100 & 285.64 & 0.02 & 299.08 & 0.03 & 523.99 & 0.06 & 214.12 & 0.03 & 138.93 & 0.06 & 357.91 & 0.03 & 289.11 & 0.03 & 274.10 & 0.05 &  &  & 551.48 & 0.06\\
1 & 0.4 & 500 & 98.66 & 0.22 & 357.07 & 0.14 & 748.60 & 0.27 & 164.57 & 0.20 & 29.09 & 0.22 & 96.85 & 0.17 & 274.20 & 0.15 & 211.31 & 0.21 &  &  & 305.79 & 0.28\\
1 & 0.4 & 5000 & 230.00 & 3.26 & 91.97 & 1.60 & 754.08 & 2.02 & 239.76 & 3.20 & 41.82 & 4.30 & 223.28 & 2.82 & 106.37 & 2.12 & 156.27 & 2.79 &  &  & 483.86 & 3.56\\
1 & 0.5 & 100 & 101.90 & 0.03 & 248.73 & 0.02 & 1444.12 & 0.04 & 427.61 & 0.03 & 77.07 & 0.02 & 1459.09 & 0.02 & 241.20 & 0.05 & 146.64 & 0.03 &  &  & 288.88 & 0.04\\
1 & 0.5 & 500 & 68.67 & 0.16 & 201.84 & 0.12 & 1628.66 & 0.19 & 146.05 & 0.13 & 86.54 & 0.27 & 108.41 & 0.12 & 105.26 & 0.22 & 127.66 & 0.15 &  &  & 456.14 & 0.17\\
1 & 0.5 & 5000 & 173.49 & 1.29 & 270.08 & 1.96 & 129.77 & 2.06 & 189.34 & 2.11 & 99.70 & 3.42 & 164.30 & 1.96 & 124.92 & 1.54 & 238.84 & 2.13 &  &  & 370.20 & 2.60\\
2 & 0 & 100 & 547.53 & 0.08 & 412.31 & 0.06 & 427.92 & 0.04 & 132.38 & 0.07 & 221.30 & 0.13 & 80.43 & 0.08 & 620.42 & 0.07 & 2200.18 & 0.06 & 112.95 & 0.05 & 449.78 & 0.04\\
2 & 0 & 500 & 79.05 & 0.37 & 264.78 & 0.33 & 364.42 & 0.37 & 690.49 & 0.44 & 316.23 & 0.66 & 488.03 & 0.44 & 438.66 & 0.51 & 419.43 & 0.45 & 61.14 & 0.30 & 517.92 & 0.44\\
2 & 0 & 5000 & 562.98 & 4.62 & 355.03 & 3.97 & 275.83 & 2.27 & 269.44 & 7.13 & 174.75 & 4.99 & 424.71 & 7.66 & 693.21 & 6.33 & 393.48 & 7.59 & 100.49 & 3.11 & 770.06 & 4.28\\
2 & 0.2 & 100 & 475.15 & 0.08 & 747.54 & 0.07 & 215.25 & 0.04 & 1096.95 & 0.06 & 87.32 & 0.12 & 423.30 & 0.09 & 583.42 & 0.05 & 1029.96 & 0.07 & 1245.96 & 0.04 & 374.18 & 0.04\\
2 & 0.2 & 500 & 389.48 & 0.39 & 544.45 & 0.37 & 226.29 & 0.19 & 876.37 & 0.41 & 53.82 & 0.49 & 1633.23 & 0.45 & 519.59 & 0.33 & 1292.09 & 0.36 & 124.16 & 0.19 & 381.84 & 0.24\\
2 & 0.2 & 5000 & 336.52 & 4.99 & 479.05 & 3.46 & 300.22 & 3.14 & 344.91 & 4.91 & 58.31 & 5.20 & 780.18 & 5.60 & 570.18 & 3.87 & 380.49 & 5.73 & 91.68 & 2.76 & 247.78 & 3.81\\
2 & 0.4 & 100 & 335.69 & 0.06 & 222.57 & 0.04 & 714.58 & 0.04 & 378.30 & 0.05 & 257.28 & 0.07 & 725.68 & 0.04 & 311.52 & 0.04 & 809.63 & 0.04 & 306.61 & 0.04 & 456.34 & 0.04\\
2 & 0.4 & 500 & 195.75 & 0.29 & 540.35 & 0.24 & 787.95 & 0.16 & 409.78 & 0.27 & 152.62 & 0.42 & 549.01 & 0.33 & 456.96 & 0.24 & 410.08 & 0.25 & 446.09 & 0.19 & 187.06 & 0.24\\
2 & 0.4 & 5000 & 177.14 & 2.51 & 320.08 & 2.69 & 76.12 & 1.16 & 1217.63 & 4.27 & 66.60 & 6.02 & 697.88 & 4.34 & 501.40 & 3.02 & 373.52 & 4.39 & 570.35 & 2.10 & 325.98 & 2.24\\
2 & 0.5 & 100 & 31.61 & 0.04 & 33.41 & 0.04 & 1032.92 & 0.03 & 33.72 & 0.04 & 30.98 & 0.05 & 31.05 & 0.05 & 33.62 & 0.04 & 31.61 & 0.04 & 31.28 & 0.05 & 24.40 & 0.04\\
2 & 0.5 & 500 & 32.02 & 0.21 & 35.95 & 0.20 & 1160.10 & 0.12 & 35.92 & 0.20 & 33.82 & 0.34 & 36.77 & 0.19 & 36.35 & 0.16 & 36.12 & 0.21 & 36.50 & 0.13 & 35.32 & 0.17\\
2 & 0.5 & 5000 & 39.86 & 2.22 & 42.74 & 2.09 & 1215.49 & 1.83 & 40.80 & 3.70 & 41.78 & 3.21 & 42.36 & 2.88 & 42.91 & 2.17 & 36.56 & 3.37 & 41.03 & 1.58 & 41.27 & 1.93\\
3 & 0 & 100 & 1318.37 & 0.07 & 392.70 & 0.07 & 195.01 & 0.06 & 912.20 & 0.09 & 350.33 & 0.13 & 91.23 & 0.08 & 126.94 & 0.04 & 117.87 & 0.05 & 41.68 & 0.05 & 739.10 & 0.05\\
3 & 0 & 500 & 194.17 & 0.40 & 583.70 & 0.36 & 138.41 & 0.21 & 786.80 & 0.43 & 305.35 & 0.51 & 304.91 & 0.92 & 104.59 & 0.35 & 379.41 & 0.37 & 96.02 & 0.35 & 722.90 & 0.31\\
3 & 0 & 5000 & 529.31 & 5.82 & 461.61 & 4.78 & 20.66 & 3.85 & 211.56 & 7.38 & 246.76 & 6.01 & 261.91 & 8.09 & 237.43 & 5.00 & 415.79 & 5.12 & 132.59 & 2.83 & 839.92 & 4.67\\
3 & 0.2 & 100 & 793.12 & 0.07 & 1050.84 & 0.07 & 255.37 & 0.05 & 740.36 & 0.07 & 982.26 & 0.09 & 500.63 & 0.07 & 770.20 & 0.05 & 907.44 & 0.05 & 1608.91 & 0.05 & 422.54 & 0.05\\
3 & 0.2 & 500 & 1065.14 & 0.35 & 1254.05 & 0.28 & 453.35 & 0.21 & 1553.27 & 0.36 & 239.91 & 0.53 & 482.26 & 0.40 & 1071.81 & 0.26 & 763.10 & 0.28 & 873.80 & 0.23 & 403.16 & 0.18\\
3 & 0.2 & 5000 & 421.41 & 4.51 & 498.22 & 4.46 & 206.14 & 2.52 & 568.44 & 5.20 & 94.76 & 2.26 & 355.08 & 6.57 & 630.81 & 4.41 & 734.31 & 3.82 & 228.76 & 2.86 & 267.18 & 2.63\\
3 & 0.4 & 100 & 468.57 & 0.05 & 570.71 & 0.06 & 652.32 & 0.04 & 653.88 & 0.06 & 553.02 & 0.08 & 636.94 & 0.10 & 803.92 & 0.05 & 504.09 & 0.04 & 500.96 & 0.08 & 467.94 & 0.05\\
3 & 0.4 & 500 & 451.09 & 0.28 & 626.09 & 0.19 & 722.05 & 0.17 & 415.71 & 0.31 & 1208.96 & 0.28 & 501.07 & 0.35 & 611.85 & 0.14 & 454.96 & 0.17 & 622.12 & 0.18 & 411.77 & 0.17\\
3 & 0.4 & 5000 & 428.67 & 4.25 & 1136.44 & 2.92 & 650.55 & 1.69 & 607.85 & 5.04 & 946.19 & 2.86 & 588.62 & 5.09 & 256.92 & 2.27 & 621.70 & 2.73 & 378.00 & 2.65 & 350.33 & 3.66\\
3 & 0.5 & 100 & 44.66 & 0.06 & 50.29 & 0.08 & 1454.08 & 0.05 & 46.91 & 0.05 & 48.04 & 0.08 & 50.77 & 0.05 & 51.69 & 0.03 & 45.41 & 0.04 & 37.79 & 0.05 & 57.23 & 0.04\\
3 & 0.5 & 500 & 53.94 & 0.19 & 54.31 & 0.13 & 1266.73 & 0.16 & 54.70 & 0.29 & 53.96 & 0.18 & 50.72 & 0.23 & 53.27 & 0.16 & 53.91 & 0.14 & 54.77 & 0.15 & 53.21 & 0.11\\
3 & 0.5 & 5000 & 62.21 & 2.20 & 51.09 & 3.07 & 1429.54 & 1.46 & 64.06 & 3.21 & 63.55 & 2.79 & 61.65 & 4.83 & 65.36 & 2.16 & 64.41 & 2.48 & 64.99 & 1.88 & 65.35 & 1.63\\
5 & 0 & 100 & 950.31 & 0.10 & 601.42 & 0.12 & 351.57 & 0.07 & 571.34 & 0.13 & 304.94 & 0.13 & 831.21 & 0.08 & 123.34 & 0.07 & 1752.19 & 0.06 & 706.40 & 0.07 & 412.43 & 0.06\\
5 & 0 & 500 & 1263.60 & 0.49 & 500.25 & 0.37 & 933.98 & 0.37 & 641.15 & 0.51 & 331.82 & 0.48 & 346.82 & 1.13 & 888.89 & 0.35 & 259.11 & 0.34 & 199.08 & 0.36 & 1012.63 & 0.18\\
5 & 0 & 5000 & 1466.39 & 7.52 & 148.12 & 5.08 & 514.01 & 2.98 & 732.56 & 8.63 & 369.66 & 5.60 & 276.46 & 15.78 & 232.57 & 5.77 & 364.63 & 7.52 & 422.92 & 3.82 & 1125.25 & 5.05\\
5 & 0.2 & 100 & 1335.14 & 0.10 & 4147.40 & 0.10 & 347.67 & 0.07 & 423.78 & 0.09 & 2887.74 & 0.12 & 888.39 & 0.08 & 783.25 & 0.06 & 803.19 & 0.11 & 726.47 & 0.09 & 860.12 & 0.05\\
5 & 0.2 & 500 & 1085.54 & 0.35 & 1477.85 & 0.29 & 301.62 & 0.18 & 996.00 & 0.50 & 3563.79 & 0.44 & 1052.70 & 0.64 & 1031.60 & 0.23 & 1082.71 & 0.40 & 1662.99 & 0.39 & 823.96 & 0.18\\
5 & 0.2 & 5000 & 608.99 & 6.11 & 196.59 & 2.78 & 44.11 & 3.54 & 1001.02 & 5.73 & 228.28 & 6.58 & 805.67 & 7.14 & 2823.58 & 3.52 & 724.28 & 6.16 & 580.61 & 4.20 & 657.57 & 3.66\\
5 & 0.4 & 100 & 914.85 & 0.13 & 1247.45 & 0.07 & 758.12 & 0.13 & 718.13 & 0.09 & 1253.16 & 0.11 & 771.06 & 0.13 & 937.43 & 0.12 & 853.46 & 0.10 & 717.50 & 0.09 & 799.58 & 0.15\\
5 & 0.4 & 500 & 1536.03 & 0.31 & 3057.03 & 0.23 & 778.53 & 0.16 & 880.26 & 0.37 & 170.57 & 0.25 & 769.97 & 0.64 & 1001.51 & 0.19 & 513.09 & 0.22 & 676.37 & 0.22 & 665.52 & 0.13\\
5 & 0.4 & 5000 & 1166.16 & 2.49 & 2278.71 & 2.09 & 910.46 & 1.65 & 706.08 & 5.11 & 151.84 & 3.57 & 1061.37 & 5.73 & 709.54 & 2.50 & 427.18 & 2.67 & 682.45 & 2.76 & 382.72 & 2.47\\
5 & 0.5 & 100 & 87.57 & 0.12 & 77.93 & 0.14 & 1514.40 & 0.11 & 76.39 & 0.11 & 87.78 & 0.14 & 80.79 & 0.16 & 85.70 & 0.16 & 89.51 & 0.16 & 87.77 & 0.14 & 117.76 & 0.14\\
5 & 0.5 & 500 & 91.37 & 0.22 & 87.67 & 0.20 & 1499.13 & 0.15 & 92.04 & 0.29 & 92.31 & 0.22 & 90.37 & 0.43 & 88.81 & 0.15 & 91.63 & 0.17 & 84.36 & 0.20 & 92.60 & 0.11\\
5 & 0.5 & 5000 & 110.58 & 2.70 & 104.39 & 1.73 & 164.34 & 1.34 & 113.49 & 4.12 & 105.37 & 3.17 & 110.37 & 7.23 & 115.36 & 2.12 & 112.16 & 2.84 & 108.93 & 2.27 & 109.26 & 2.79\\
10 & 0 & 100 & 2757.84 & 0.35 & 4326.25 & 0.25 & 107.44 & 0.40 & 1183.30 & 0.36 & 639.67 & 0.38 & 1071.85 & 0.46 & 3277.33 & 0.38 & 3625.45 & 0.42 & 2152.67 & 0.44 & 1855.01 & 0.46\\
10 & 0 & 500 & 2590.14 & 0.71 & 1095.09 & 0.36 & 206.05 & 0.31 & 1080.42 & 0.90 & 967.88 & 0.69 & 1363.92 & 0.98 & 759.82 & 0.38 & 1186.38 & 0.37 & 641.53 & 0.38 & 1384.22 & 0.29\\
10 & 0 & 5000 & 3645.12 & 13.33 & 309.22 & 4.99 & 272.86 & 3.35 & 1139.06 & 10.41 & 1026.37 & 8.65 & 2049.83 & 12.09 & 724.34 & 5.54 & 797.15 & 3.46 & 627.56 & 7.13 & 1261.36 & 3.81\\
10 & 0.2 & 100 & 923.73 & 0.35 & 708.45 & 0.22 & 515.29 & 0.38 & 1298.31 & 0.30 & 5686.52 & 0.34 & 1460.27 & 0.41 & 934.70 & 0.37 & 1381.85 & 0.35 & 1149.68 & 0.27 & 1483.71 & 0.40\\
10 & 0.2 & 500 & 1026.55 & 0.72 & 471.09 & 0.19 & 967.80 & 0.26 & 1101.19 & 0.68 & 726.47 & 0.48 & 1396.33 & 0.45 & 829.65 & 0.28 & 1119.06 & 0.29 & 715.50 & 0.42 & 1561.98 & 0.27\\
10 & 0.2 & 5000 & 1235.41 & 2.58 & 2145.38 & 2.95 & 1172.67 & 2.81 & 1068.90 & 9.90 & 635.84 & 6.29 & 1609.09 & 10.02 & 577.08 & 3.38 & 1391.16 & 2.63 & 882.57 & 6.88 & 1300.23 & 3.13\\
10 & 0.4 & 100 & 1127.56 & 0.23 & 1164.25 & 0.16 & 725.52 & 0.25 & 1406.01 & 0.19 & 3755.16 & 0.26 & 1445.14 & 0.37 & 1258.34 & 0.24 & 1450.73 & 0.33 & 1519.86 & 0.26 & 1865.13 & 0.29\\
10 & 0.4 & 500 & 1073.37 & 0.36 & 901.27 & 0.21 & 769.89 & 0.18 & 1245.44 & 0.54 & 950.66 & 0.43 & 1351.65 & 0.41 & 1284.80 & 0.24 & 1385.30 & 0.31 & 1179.80 & 0.44 & 1484.50 & 0.23\\
10 & 0.4 & 5000 & 1660.87 & 2.91 & 669.08 & 2.72 & 83.95 & 1.79 & 1309.16 & 8.91 & 627.09 & 3.41 & 1581.53 & 4.79 & 806.69 & 2.31 & 1282.66 & 2.24 & 1359.75 & 3.91 & 2162.42 & 2.18\\
10 & 0.5 & 100 & 177.63 & 0.19 & 173.75 & 0.11 & 1788.45 & 0.21 & 161.28 & 0.15 & 177.12 & 0.20 & 178.72 & 0.33 & 176.50 & 0.25 & 186.67 & 0.33 & 178.34 & 0.20 & 145.08 & 0.29\\
10 & 0.5 & 500 & 189.69 & 0.32 & 187.19 & 0.17 & 1092.62 & 0.19 & 188.99 & 0.42 & 192.57 & 0.29 & 187.90 & 0.35 & 195.14 & 0.20 & 192.07 & 0.81 & 188.00 & 0.23 & 191.50 & 0.16\\
10 & 0.5 & 5000 & 231.75 & 4.74 & 232.26 & 1.87 & 1345.41 & 2.32 & 232.99 & 6.69 & 198.40 & 4.14 & 229.13 & 5.63 & 201.72 & 2.30 & 239.00 & 3.37 & 217.86 & 2.84 & 229.90 & 1.55\\
\hline\hline
\end{tabular}
} 
\end{table}
\setlength{\tabcolsep}{6pt}
\clearpage



\clearpage

\clearpage

\begin{figure}
\centering
\subfloat[noise level 0\%]{\label{fg:c2-d1_020_500_p2_f8_b8_m4_i200_n0_h10_cuda_test}
 \includegraphics[width=0.3\textwidth]{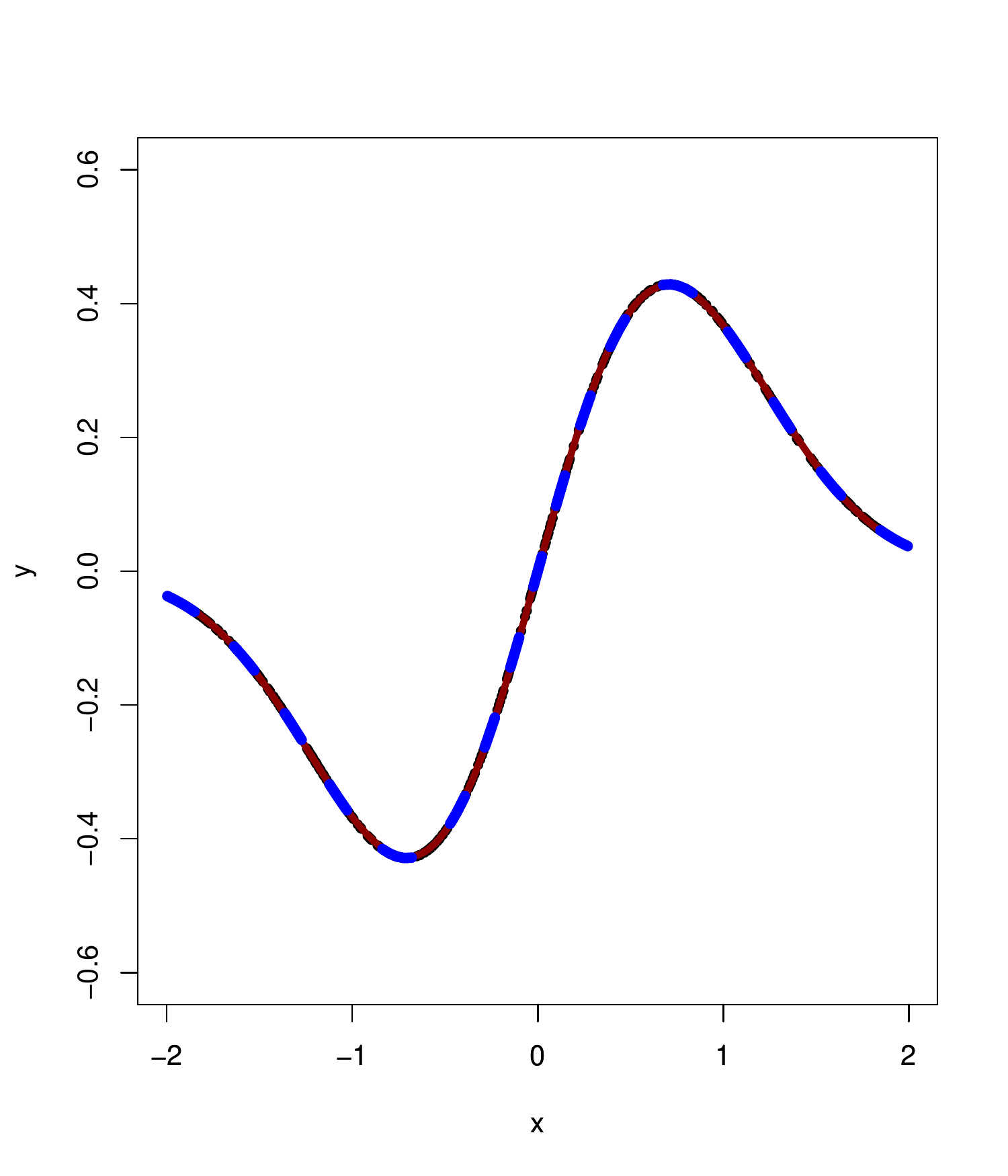}}
\subfloat[noise level 10\%]{\label{fg:c2-d1_020_500_p2_f8_b8_m4_i200_n10_h10_cuda_test}
 \includegraphics[width=0.3\textwidth]{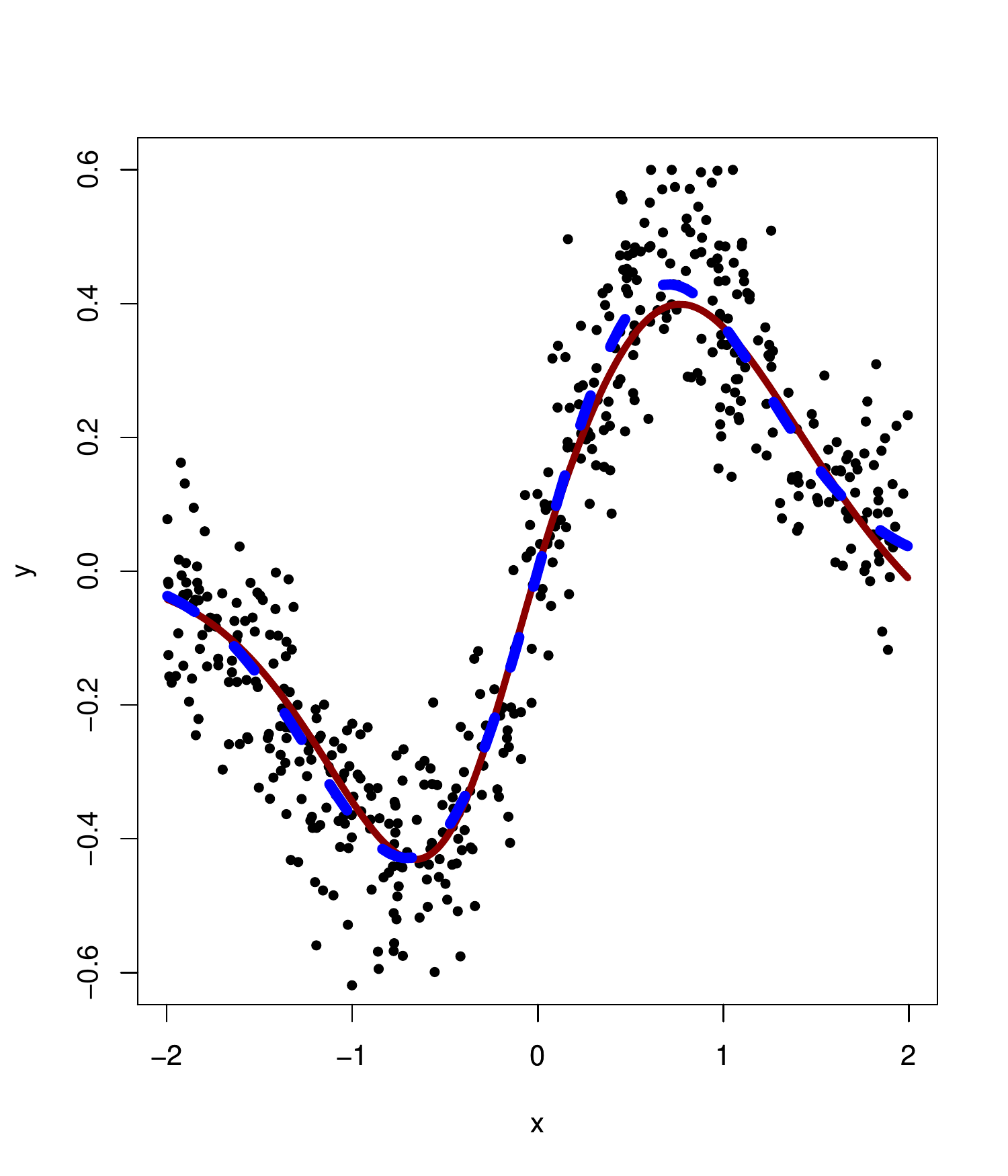}}
\subfloat[noise level 20\%]{\label{fg:c2-d1_020_500_p2_f8_b8_m4_i200_n20_h10_cuda_test}
 \includegraphics[width=0.3\textwidth]{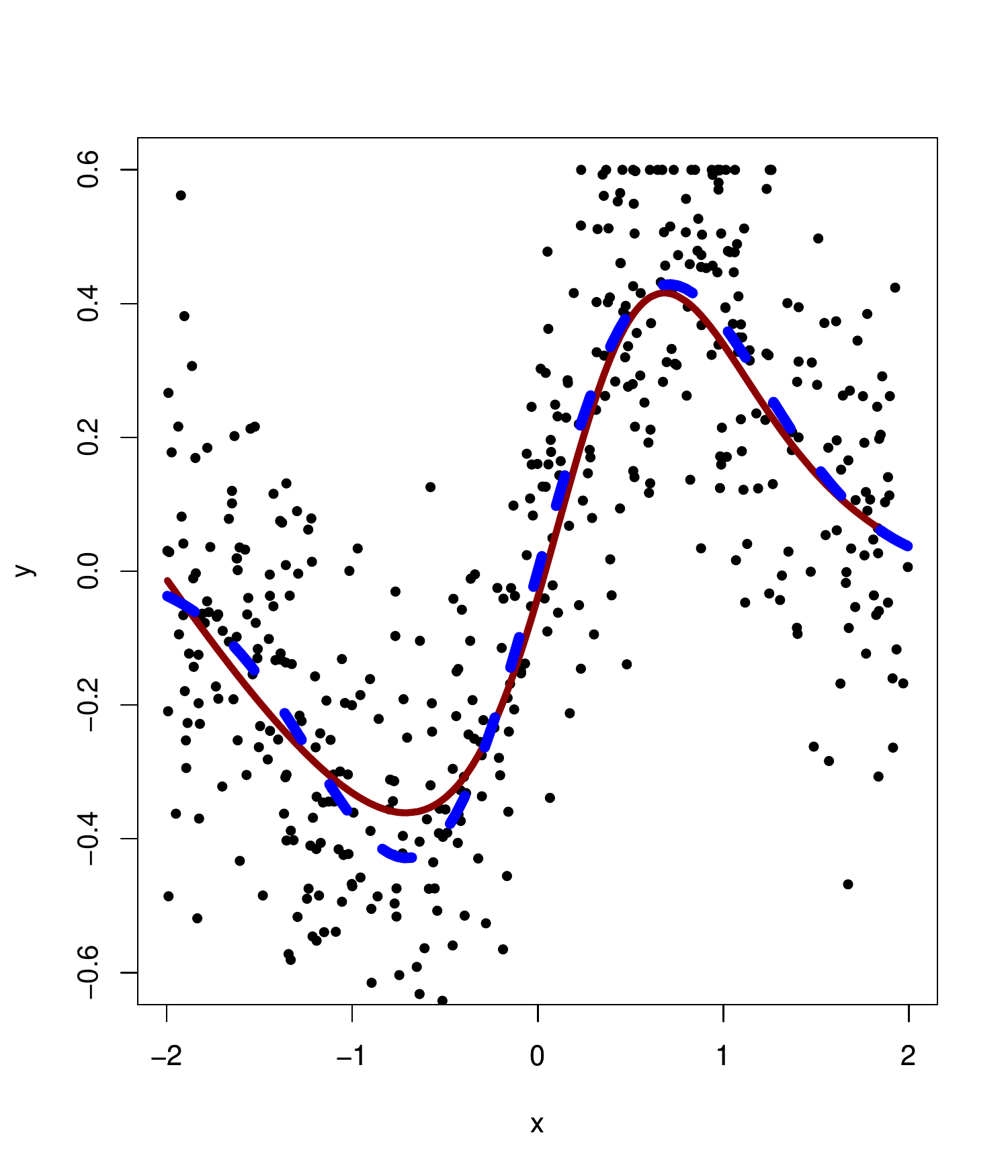}}
\caption{Data set 2. The proportion of outliers is $\delta=0.2$, and $n=500$. Gaussian noise in the data is 0, 0.1 and 0.2 respectively. Red solid curve is the PCLTS  prediction, the blue dashed curve is the (noiseless) test data, and black dots are training data. The outliers are not shown.}
\label{fg:ds2}
\end{figure}

\begin{figure}
\centering
\subfloat[noise level 0\% ]{\label{fg:c4-d1_020_500_p2_f8_b8_m4_i200_n0_h10_cuda_nnet_test_predictions}
 \includegraphics[width=0.3\textwidth]{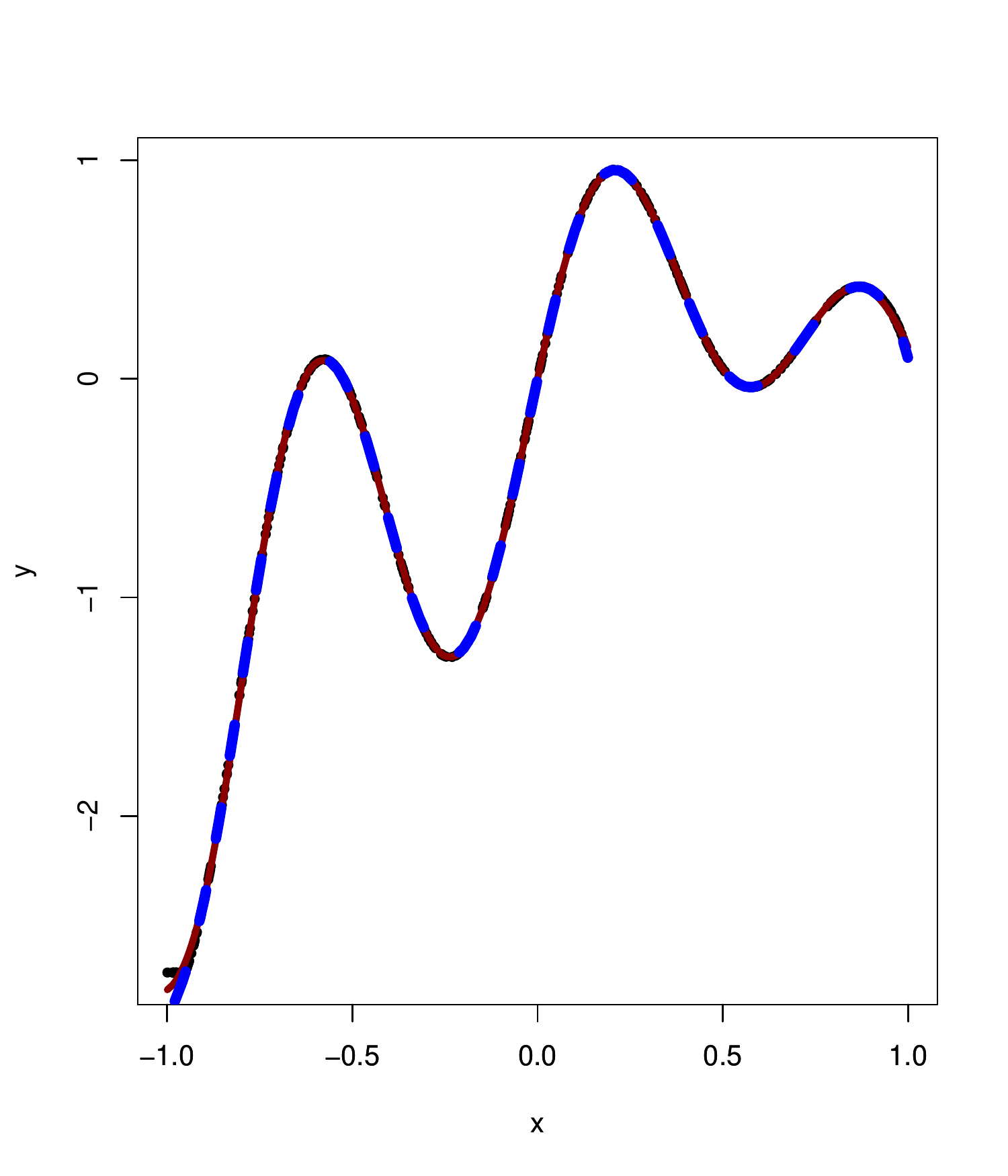}}
\subfloat[noise level 10\%]{\label{fg:c4-020_500_p2_f8_b8_m4_i200_n10_h10_cuda_nnet_short_train_predictions}
 \includegraphics[width=0.3\textwidth]{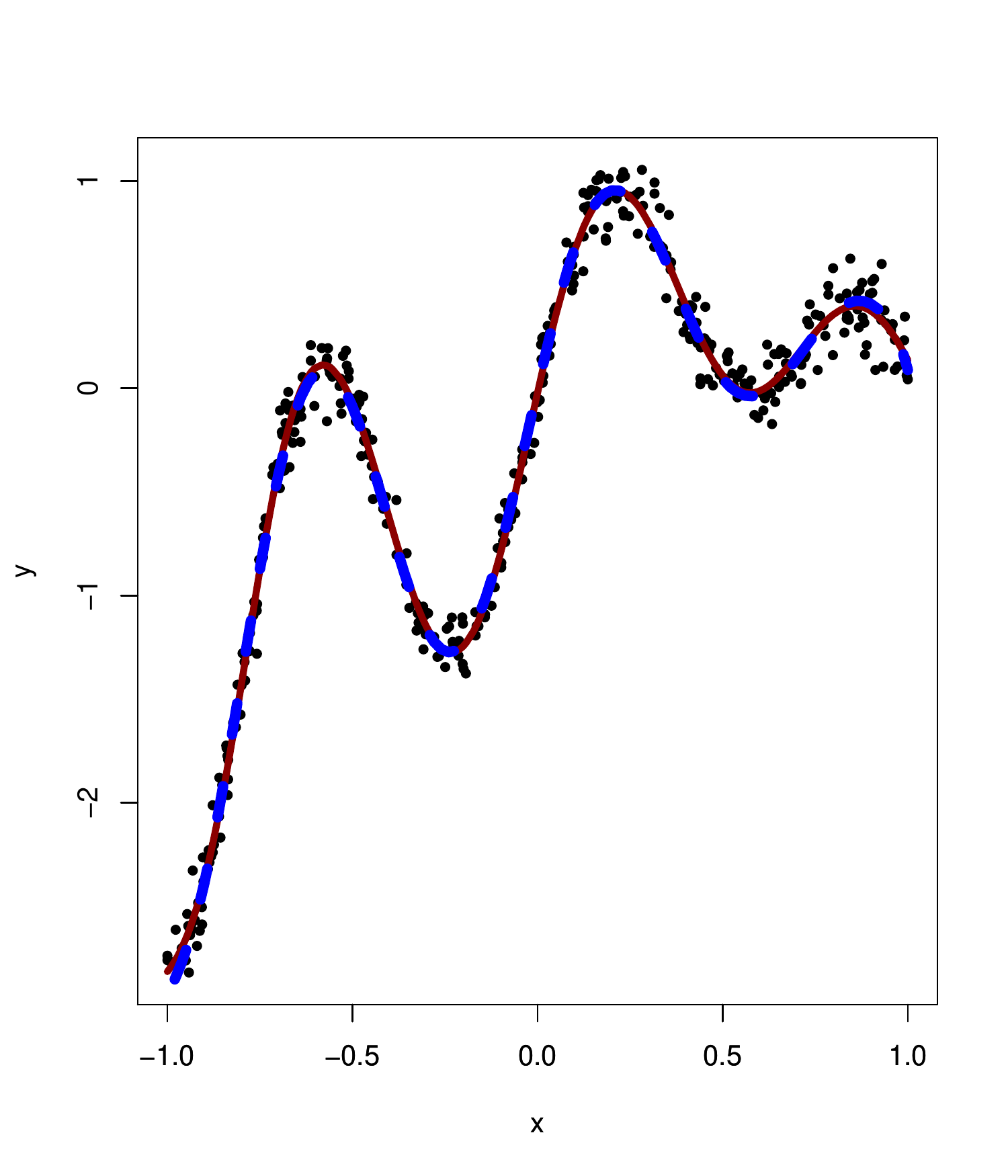}}
\subfloat[noise level 20\%]{\label{fg:c4-020_500_p2_f8_b8_m4_i200_n20_h10_cuda_nnet_short_train_predictions}
 \includegraphics[width=0.3\textwidth]{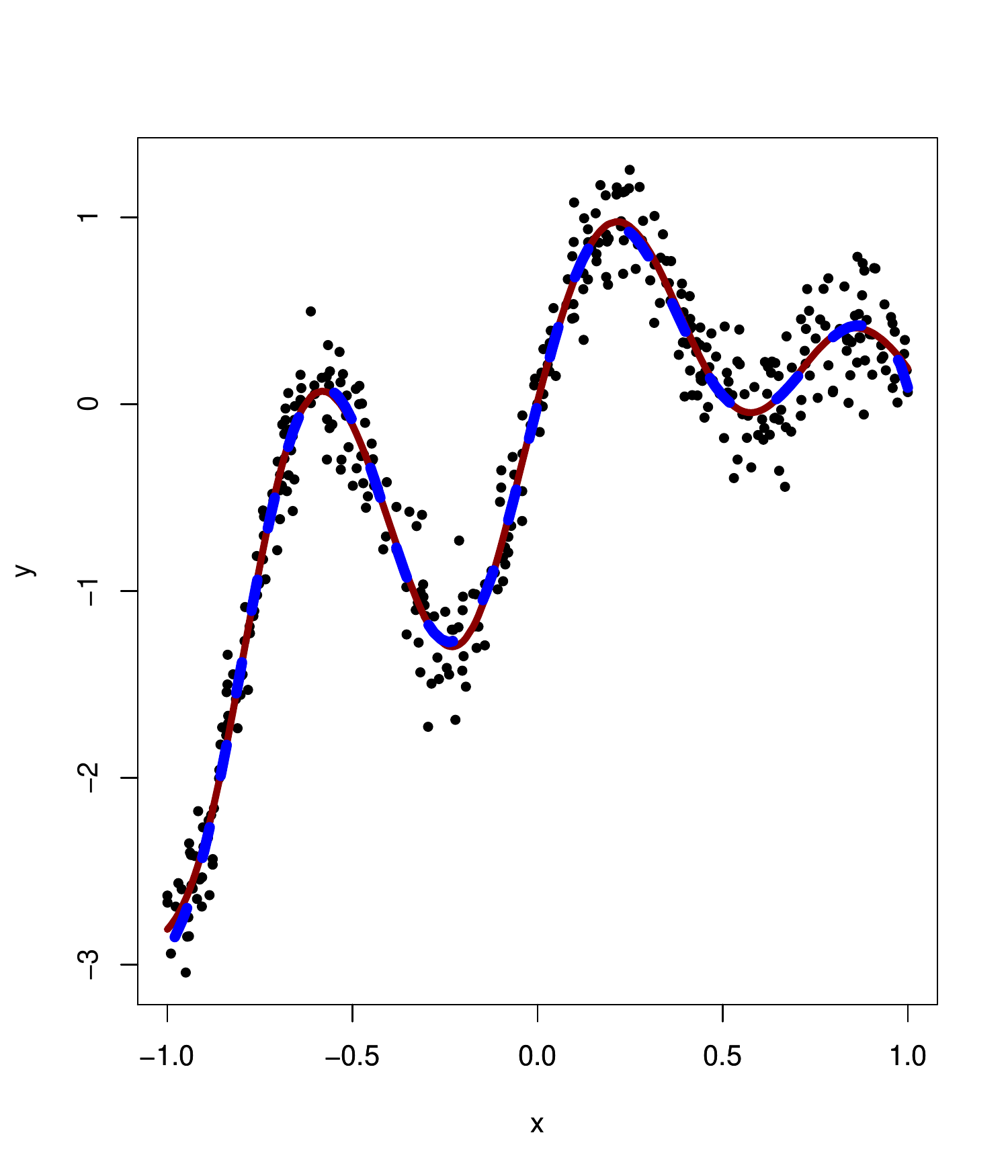}}
\caption{Data set 4, $\delta=0.2$, $n=500$. PCLTS ($C=8$, $B=8$) correctly predicts the model.}
\label{fg:ds4}
\end{figure}

\begin{figure}
\centering
\subfloat[noise level 0\%]{\label{fg:c5-d1_020_500_p2_f8_b8_m4_i200_n0_h10_cuda_nnet_test_predictions}
 \includegraphics[width=0.3\textwidth]{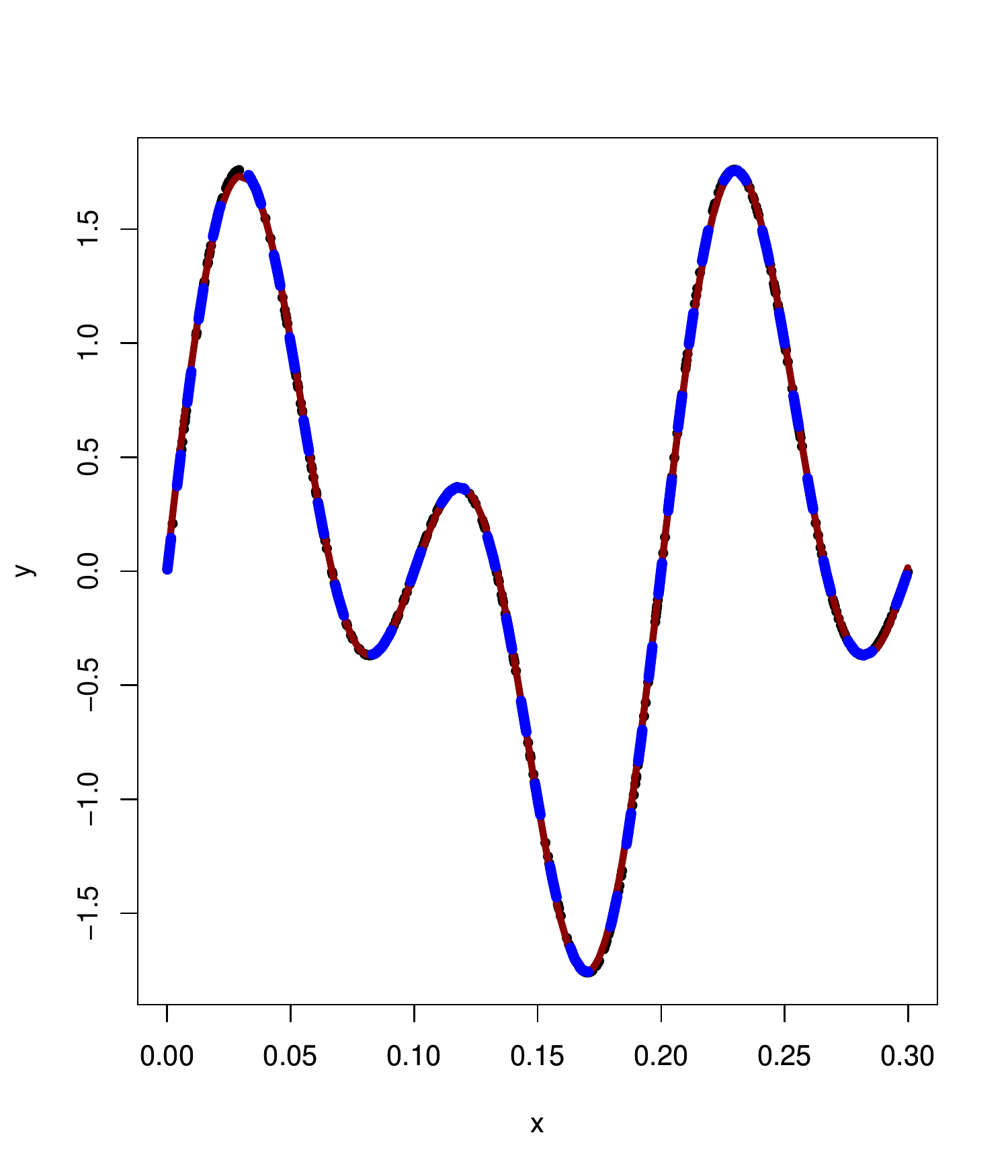}}
\subfloat[noise level 10\% ]{\label{fg:c5-d1_020_500_p2_f8_b8_m4_i200_n10_h10_cuda_nnet_short_train_predictions}
 \includegraphics[width=0.3\textwidth]{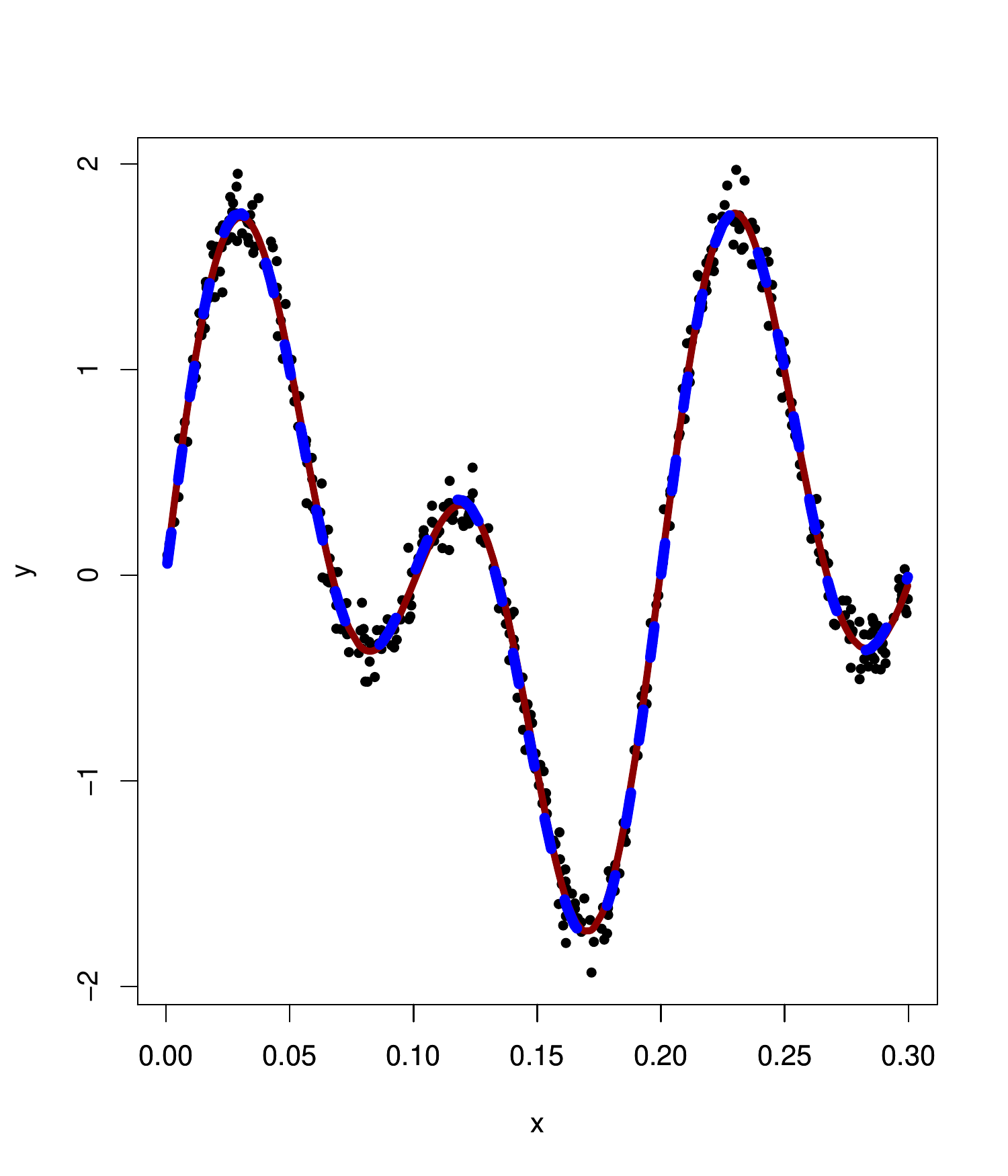}}
\subfloat[noise level 20\%]{\label{fg:c5-d1_020_500_p2_f8_b8_m4_i200_n20_h10_cuda_nnet_short_train_predictions}
 \includegraphics[width=0.3\textwidth]{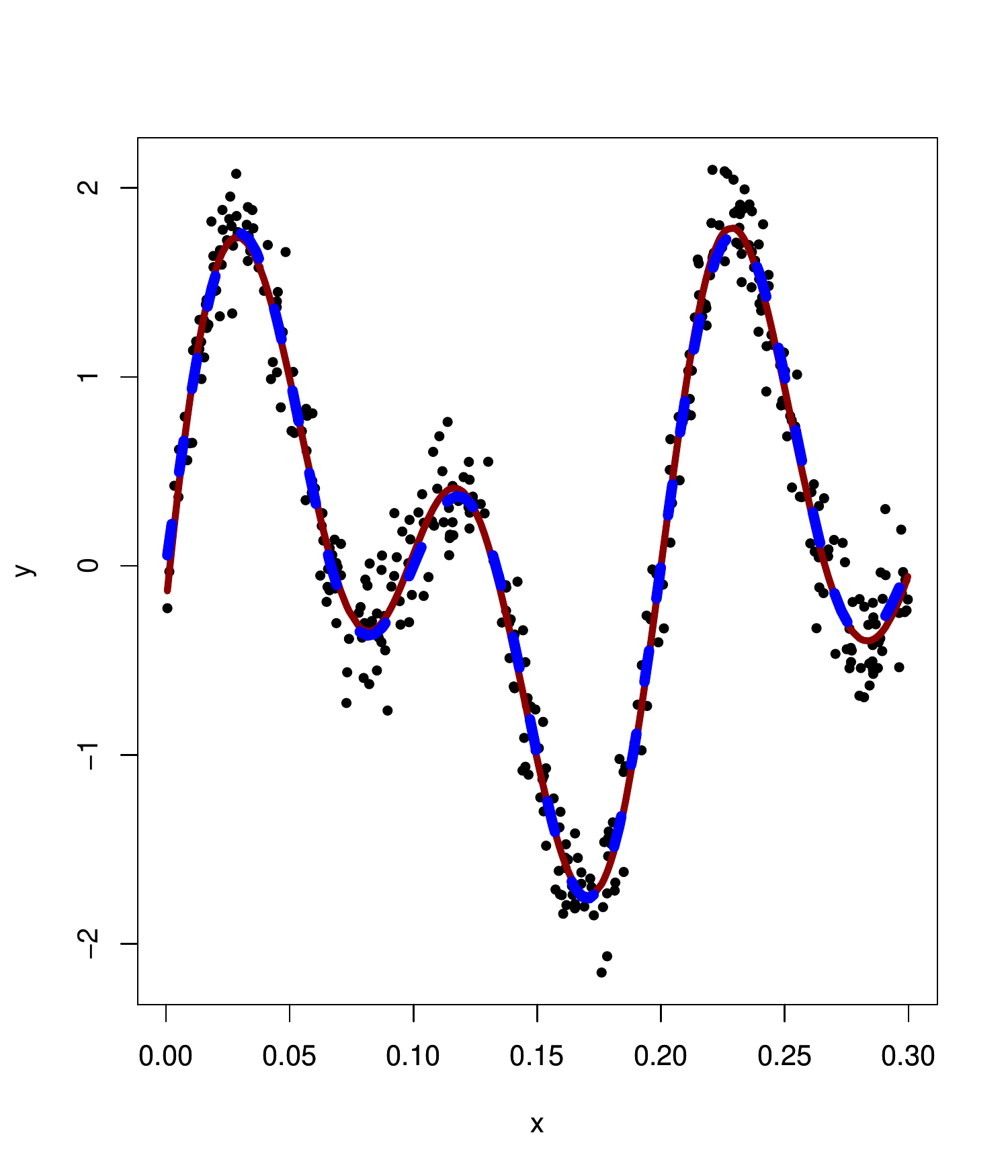}}
\caption{Data set 5, $\delta=0.2$, $n=500$. PCLTS ($C=8$, $B=8$) correctly predicts the model.}
\label{fg:ds5}
\end{figure}

\begin{figure}[ph!]
\centering
\subfloat[noise level 0\% ]{\label{fg:c6-d1_020_500_p2_f8_b8_m4_i200_n0_h10_cuda_nnet_test_predictions}
 \includegraphics[width=0.3\textwidth]{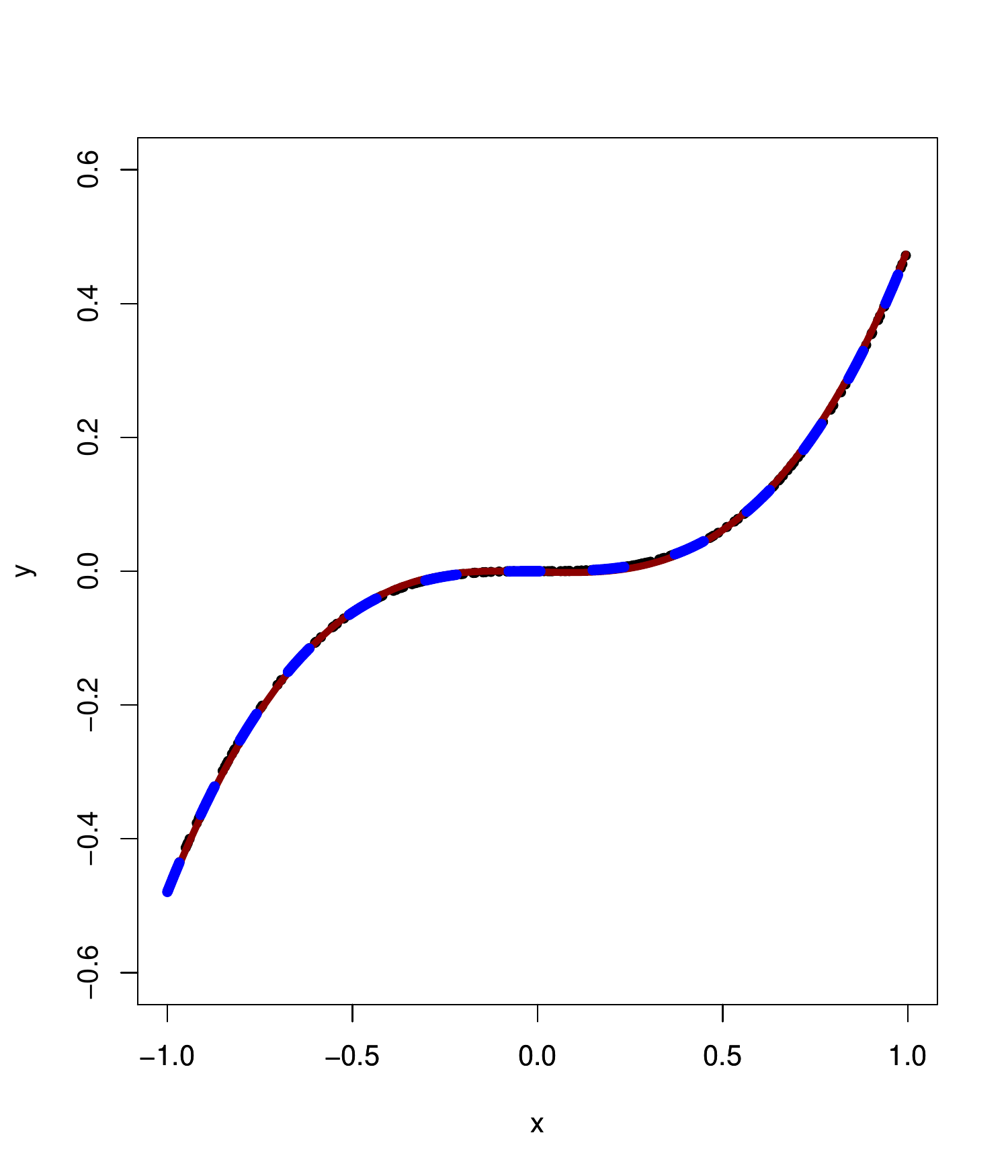}}
\subfloat[noise level 10\%]{\label{fg:c6-d1_020_500_p2_f8_b8_m4_i200_n10_h10_cuda_nnet_test_predictions}
 \includegraphics[width=0.3\textwidth]{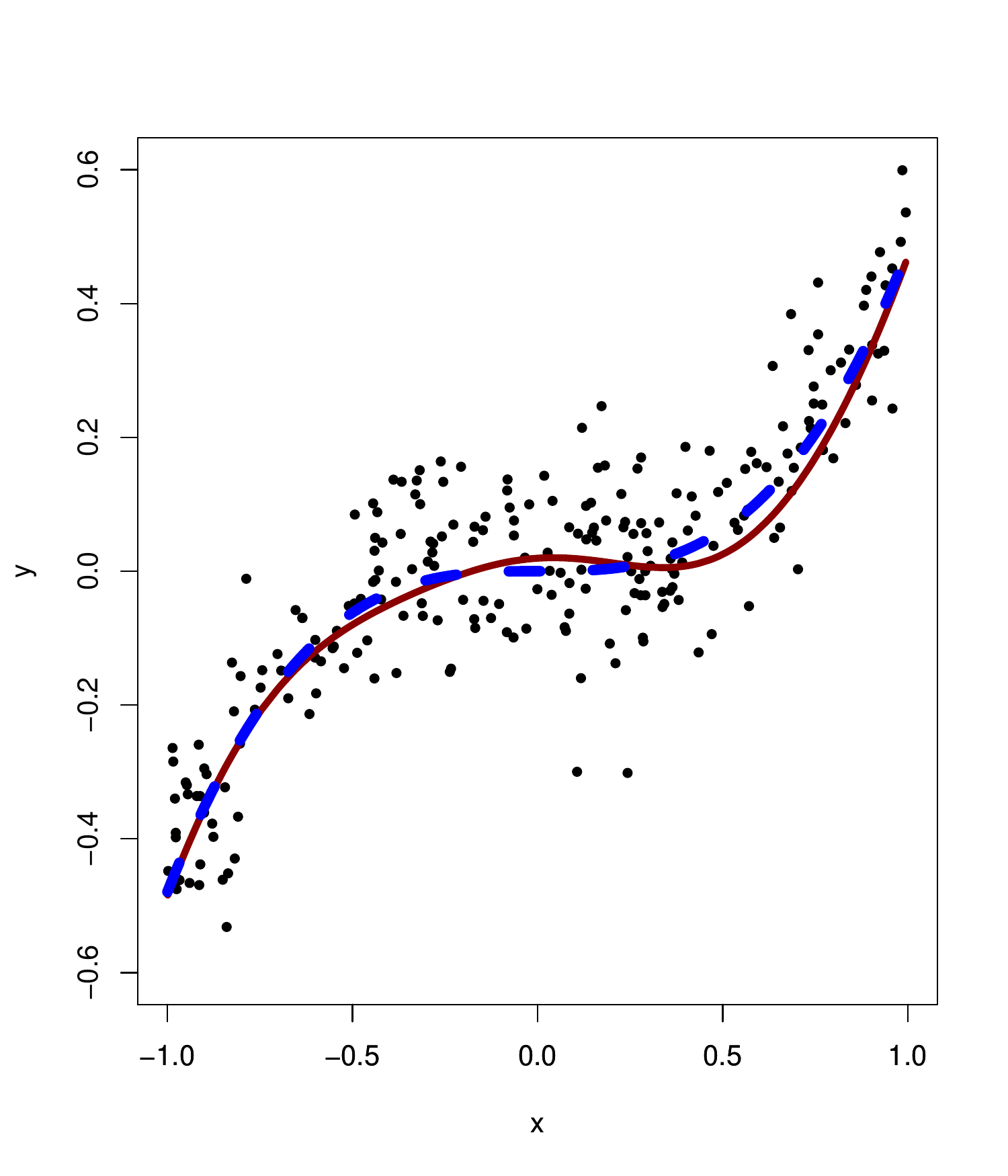}}
\subfloat[noise level 20\%]{\label{fg:c6-d1_020_500_p2_f8_b8_m4_i200_n20_h10_cuda_nnet_test_predictions}
 \includegraphics[width=0.3\textwidth]{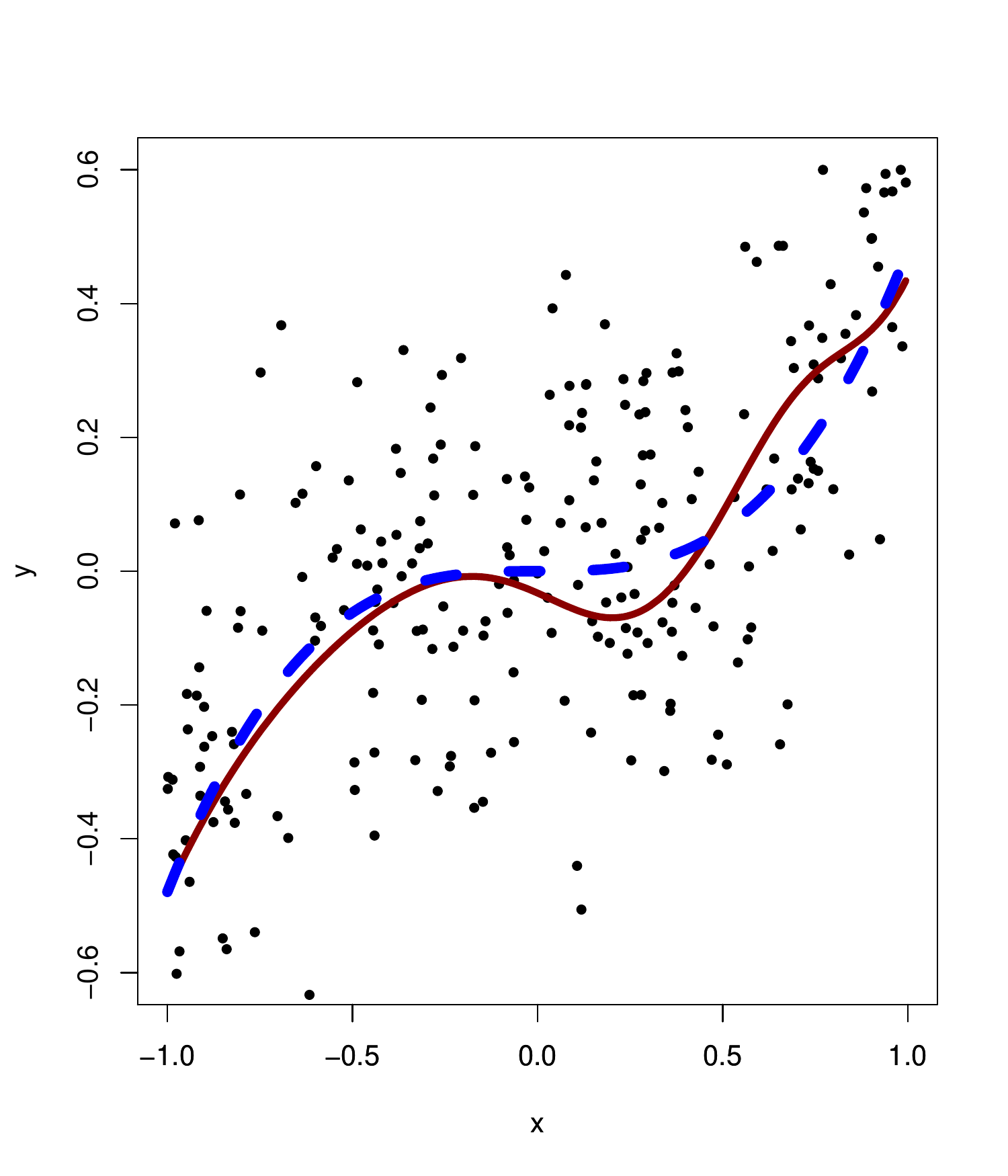}}
\caption{Data set 6, $\delta=0.2$, $n=500$. PCLTS ($C=8$, $B=8$) correctly predicts the model.}
\label{fg:ds6}
\end{figure}

\begin{figure}
\centering
\subfloat[noise level 0\%]{\label{fg:c7-d1_020_500_p2_f8_b8_m4_i200_n0_h10_cuda_nnet_test_predictions}
 \includegraphics[width=0.3\textwidth]{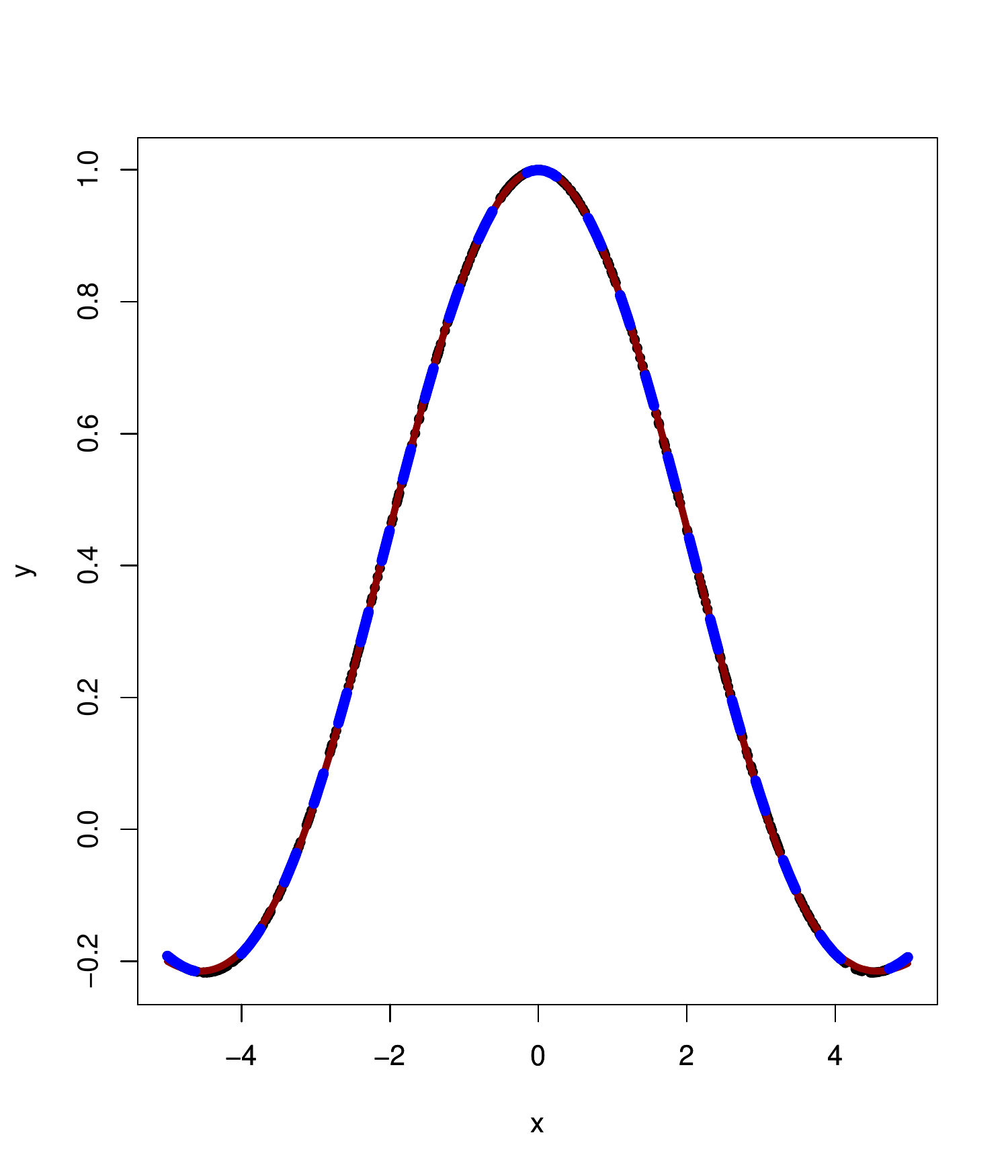}}
\subfloat[noise level 10\%]{\label{fg:c7-d1_020_500_p2_f8_b8_m4_i200_n10_h10_cuda_nnet_test_predictions}
 \includegraphics[width=0.3\textwidth]{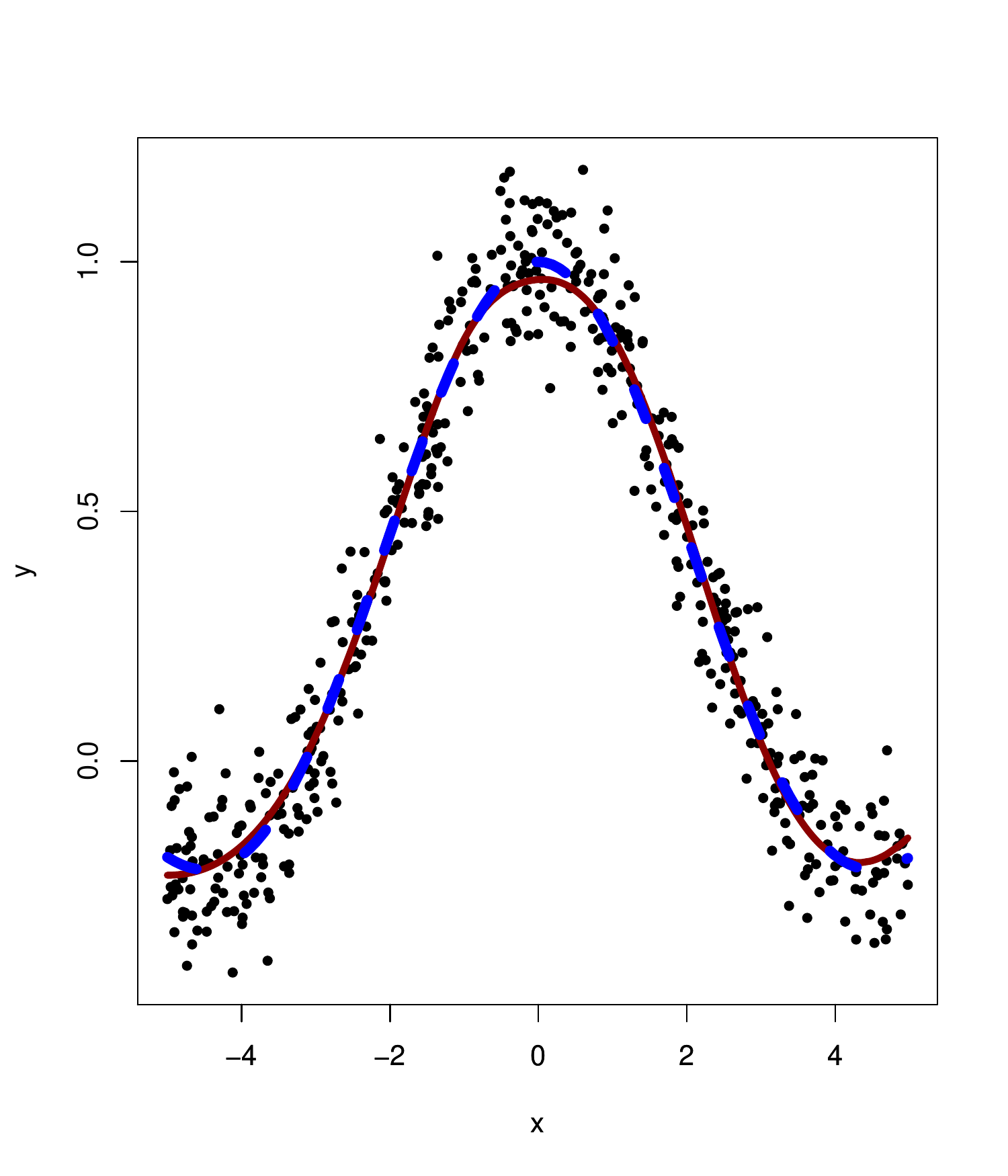}}
\subfloat[noise level 20\%]{\label{fg:c7-d1_020_500_p2_f8_b8_m4_i200_n20_h10_cuda_nnet_test_predictions}
 \includegraphics[width=0.3\textwidth]{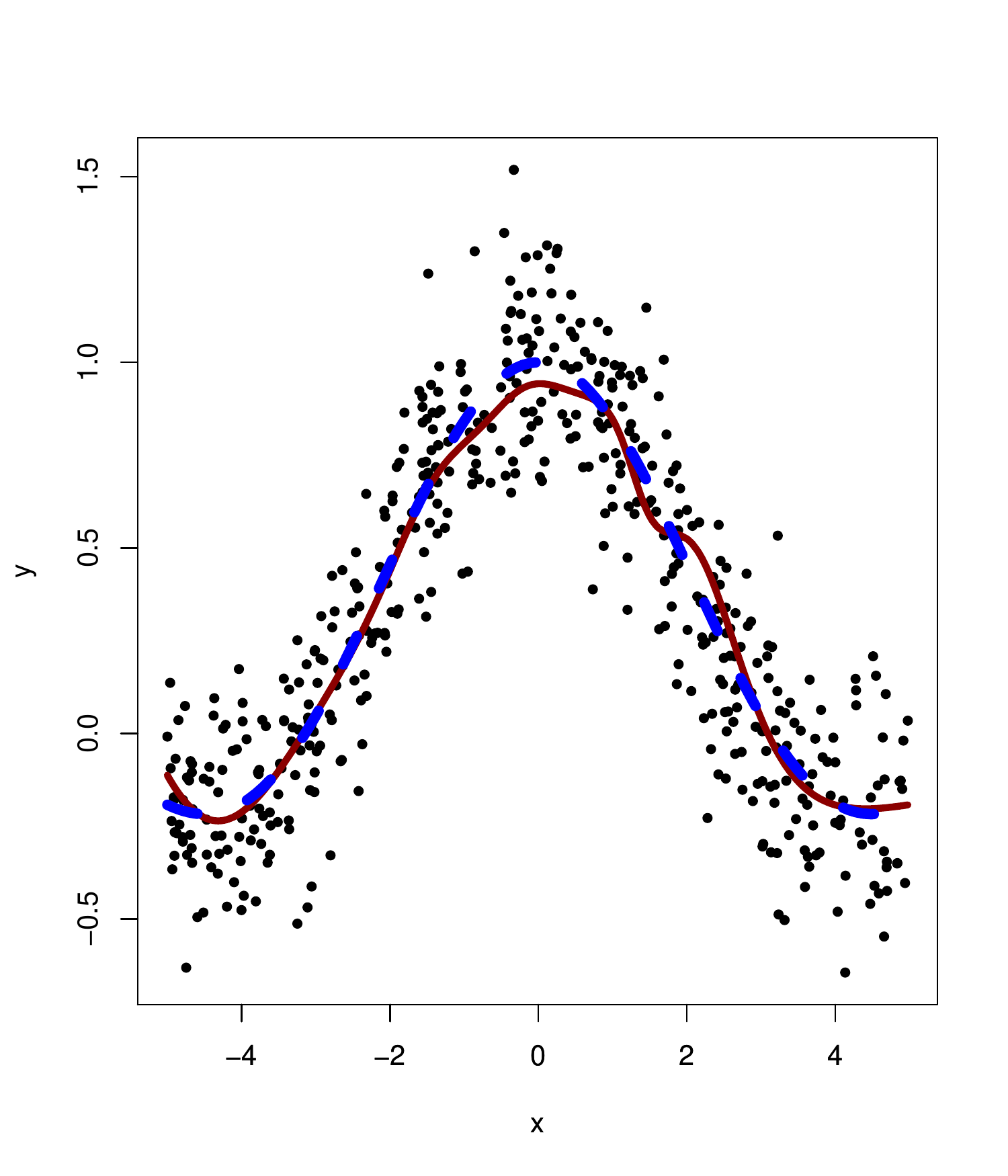}}
\caption{Data set 7, $\delta=0.2$, $n=500$. PCLTS ($C=8$, $B=8$) correctly predicts the model.}
\label{fg:ds7}
\end{figure}

\end{document}